\documentclass[11pt]{article}
\usepackage[margin=2.0cm]{geometry}
\usepackage[utf8]{inputenc}
\usepackage[T1]{fontenc}
\usepackage[english]{babel}
\usepackage{amsmath,amsfonts,amssymb,amsthm,mathtools}
\usepackage{bm}
\usepackage{graphicx}
\usepackage{float}
\usepackage{booktabs}
\usepackage[authoryear,longnamesfirst]{natbib}
\usepackage[colorlinks=true,linkcolor=blue,citecolor=blue,urlcolor=blue]{hyperref}
\allowdisplaybreaks\hbadness=10000
\theoremstyle{plain}
\newtheorem{theorem}{Theorem}
\newtheorem{proposition}{Proposition}
\newtheorem{corollary}[proposition]{Corollary}
\newtheorem{lemma}{Lemma}
\theoremstyle{definition}
\newtheorem{remark}{Remark}

\newcommand{\E}{\mathbb{E}}
\newcommand{\PP}{\mathbb{P}}
\newcommand{\Var}{\mathrm{Var}}
\newcommand{\Cov}{\mathrm{Cov}}
\newcommand{\ind}{\mathbf{1}}
\newcommand{\Sph}{\mathbb{S}}
\newcommand{\GOE}{\mathrm{GOE}}
\newcommand{\lmax}{\lambda_{\max}}
\newcommand{\lmin}{\lambda_{\min}}

\newcommand{\ba}{\mathbf{a}}
\newcommand{\Mb}{\overline{M}}
\newcommand{\doi}[1]{\href{https://doi.org/#1}{\nolinkurl{https://doi.org/#1}}}
\renewcommand{\paragraph}[1]{\par\smallskip\noindent\textit{#1.}\hspace{0.5em}}

\title{The $t$-spacing test for a spiked tensor}
\usepackage{authblk}
\author[1]{Jean-Marc Aza\"{i}s}
\author[2]{Federico Dalmao}
\author[3,4]{Yohann De Castro\thanks{Corresponding author. E-mail: \texttt{yohann.de-castro@ec-lyon.fr}. ORCID: 0000-0002-9008-7474.}}
\affil[1]{{\small Institut de Math\'ematiques de Toulouse\\ Universit\'e de Toulouse, France.}}
\affil[2]{{\small DMEL, CENUR Litoral Norte\\ Universidad de la Rep\'ublica, Salto, Uruguay.}}
\affil[3]{{\small Institut Camille Jordan, CNRS UMR 5208\\ \'Ecole Centrale Lyon, France.}}
\affil[4]{{\small Institut Universitaire de France~(IUF)}}
\date{November 27, 2026}

\begin{document}
\maketitle

\begin{abstract}
We observe one symmetric tensor $Y=\lambda\,e^{\otimes k}+\sigma W$ of odd order $k$ on $\mathbb R^{d}$, with
Gaussian noise $W$, and test $\lambda=0$. The observation does not determine the noise level $\sigma$, so a test
that uses it needs external information. The $t$-spacing test of Aza\"is, Dalmao, and De Castro needs none and is
exact at every $\sigma$, $k$ and $d$, in closed form. We study its power. At every $\lambda$ we bound the
probability that the test misses a signal of amplitude $\lambda$ by a noncentral Student distribution function. At
fixed $(k,d)$ this probability decays exponentially in~$\lambda^{2}$ as $\lambda$ grows, and we identify the rate.
The rate is below the lower bound we prove for the rate of the known-variance spacing test. It tends to this bound
as~$d$ grows. As $d$ grows at fixed $k$, the power tends to one when $\lambda/\sigma$ is a large enough multiple of
$\sqrt d$. The test also returns an estimator of the
direction~$e$, the maximiser of the profile likelihood. It is asymptotically minimax as $\lambda$ grows, with the sharp constant
$(d-1)/(k\lambda^{2})$ for the risk $\E\sin^{2}\angle(\hat e,e)$: a Bayesian Cram\'er--Rao lower bound valid at
every $\lambda$, and an upper bound with the same leading term.
\end{abstract}


\smallskip\noindent\textbf{MSC2020 subject classifications.} Primary 62F03, 60G60; secondary 62C20, 15B52

\section{Introduction}\label{sec:intro}

Let $Y=\lambda\,e^{\otimes k}+\sigma W$ be an observed symmetric tensor of odd order~$k\ge3$ on $\mathbb R^{d}$, where
$e\in\Sph^{d-1}$ is an unknown direction, $\lambda\ge0$ an unknown amplitude, $\sigma>0$ an unknown noise level,
and $W$ a Gaussian symmetric tensor normalised so that~$X(t)=\langle W,t^{\otimes k}\rangle$ is the
Kostlan--Shub--Smale field, $\E X(s)X(t)=\langle s,t\rangle^{k}$. Testing $\lambda=0$ against $\lambda>0$ is the
detection problem of tensor principal component analysis \citep{montanari2014statistical,arous2019landscape}; the
maximiser of the profile likelihood~$Z(t)=\langle Y,t^{\otimes k}\rangle$ is the associated estimator of $e$.

\medskip

\begin{figure}[t]
\centering
\includegraphics[width=\textwidth]{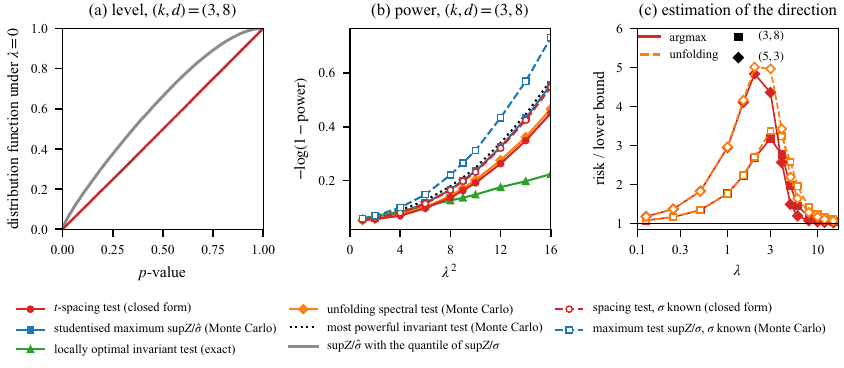}
\caption{The $t$-spacing test against the competitors of Section~\ref{sec:num}, $10^{5}$ replicates per point.
(a) Under $\lambda=0$, the distribution function of each $p$-value. Only $\sup Z/\hat\sigma$ with the quantile of
$\sup Z/\sigma$ is not uniform: it rejects with probability $0.097$ at the nominal level $0.05$. (b) Power at level
$0.05$, on a logarithmic scale. The $t$-spacing test is below the simulated most powerful invariant test by at most
$0.08$. (c) Risk $\E\sin^{2}\angle(\hat e,e)$ of two estimators of the direction, divided by the lower bound of
Theorem~\ref{thm:estim}. The ratio tends to one at both ends. For the maximiser of the profile likelihood it is
within $10\%$ of one at $\lambda=10$.}
\label{fig:intro}
\end{figure}

Exact tests are hard to come by here. The null law of $\sup Z$ has no closed form. A scale-free statistic, such as
$\sup Z$ divided by an estimate of $\sigma$, can be calibrated by simulation, but each replicate costs one full
optimisation of a non-convex landscape, for every pair $(k,d)$. Comparing it with the quantile of $\sup Z/\sigma$
is not valid: at $(k,d)=(3,8)$ this test rejects too often (Figure~\ref{fig:intro}(a)). Asymptotic calibrations are of little use at
the dimensions where the model is interesting. The spacing test and the $t$-spacing test of \citet{azais2024second},
related to the Kac--Rice test of \citet{taylor2016inference}, solve the calibration problem in closed form: they
are uniformly distributed on $(0,1)$ under the null, at every dimension and every degree, and the second one does
not use $\sigma$. They need no simulation, and the maximiser of $Z$ that they compute is the estimator of $e$. The power
of the spacing test was studied for least-angle regression by \citet{azais2018power}; for these tests in the tensor
model nothing was known on the power.
Figure~\ref{fig:intro} shows the picture.

\paragraph{Why the $t$-test, and not the known-variance one} In tensor principal component analysis the noise
level is often unknown in practice \citep{xia2022inference}. It is tempting to regard the $t$-spacing test as the
fallback one uses when $\sigma$ is unavailable. In this model it is the other way round.
Section~\ref{sub:rank} shows that the kernel $\langle s,t\rangle^{k}$ has finite rank
$m=\binom{d+k-1}{k}$: the whole path of $\sigma X$ determines exactly $m$ Gaussian coordinates and nothing more,
so $\sigma$ is not a function of the observation and the known-variance test requires information from outside the
data. The $t$-spacing test needs no such information. It estimates the noise level with $\kappa=m-d$ degrees of
freedom. In simulations its power is below that of the known-variance test, by an amount that decreases as $\kappa$
grows (Sections~\ref{sec:low} and~\ref{sec:num}).

\paragraph{Contributions} The test is invariant under rotations and the estimator is equivariant. So the power of
the test and the risk of the estimator are the same at every direction $e$.
\begin{enumerate}
\item \emph{Power.} At every $\lambda$, Corollary~\ref{cor:closed} bounds the probability that the $t$-spacing test
      misses a signal of amplitude $\lambda$ by a noncentral Student distribution function. This probability decays
      exponentially in $\lambda^{2}$. We identify the rate $E$ at every $(k,d)$ and every level
      (Theorem~\ref{thm:exponent}). Its constant is a maximum over the spectra of the normalised Hessian, computed for
      $d=3$ and bracketed for $d\ge4$ in Remark~\ref{rem:bracket}. The lower bound on the rate follows from
      Corollary~\ref{cor:closed}. As $d$ grows at
      fixed $k$, the power tends to one when $\lambda/\sigma$ is a large enough multiple of $\sqrt d$
      (Theorem~\ref{thm:sqrtd}). The bound in closed form, Corollary~\ref{cor:closed}, is informative only on the
      larger scale $\lambda/\sigma\asymp d^{k/2}$ (Remark~\ref{rem:dim}).
\item \emph{Minimax estimation.} The test also returns an estimator of the direction $e$, the
      argmax of the profile likelihood. At fixed $(k,d)$ it is asymptotically minimax as $\lambda\to\infty$ for the
      squared sine of the angle to $e$, with the sharp constant $(d-1)/(k\lambda^{2})$ (Theorem~\ref{thm:estim},
      Corollary~\ref{cor:minimax}). Theorem~\ref{thm:estim} is non-asymptotic: a lower bound valid at every
      $\lambda$, and an upper bound with an explicit remainder which controls the probability that the argmax of a
      non-convex landscape with many critical points is far from $e$. The constant is the classical one: at fixed
      $(k,d)$ the model is a curved exponential family in $e$ with Fisher information $k\lambda^{2}\,\mathrm{Id}$ on
      the tangent space, and the small-noise theory of the maximum likelihood estimator gives this constant for its
      limit law. The difficulty is not the local limit but the global landscape.
\item \emph{Benchmarks.} With $\sigma$ known, Corollary~\ref{cor:Sclosed} bounds the probability of a miss for the
      spacing test, and gives the lower bound $k/(2(3k-2))$ on its rate. The rate $E$ of the $t$-spacing test is
      below this bound at every $(k,d)$ (Remark~\ref{rem:high}) and tends to it as $d\to\infty$
      (Remark~\ref{rem:dim}). With~$\sigma$ known, the test on the maximum of the profile likelihood has the rate
      $1/2$. At high signal the Bayes-optimal test, which averages the likelihood over directions, differs from it
      by a level shift of order $d\log\lambda/\lambda$ and has the same rate (Proposition~\ref{prop:laplace},
      Remark~\ref{rem:high}). At vanishing signal the locally optimal test uses quadratic forms, not the $t$-test.
      At $(k,d)=(3,3)$, with~$\sigma$ known, it is above the maximum test only where both have power below $0.19$
      at level~$0.05$. It gives no estimator of $e$ (Section~\ref{sec:low}).
\end{enumerate}

\paragraph{Growing dimension} In tensor principal component analysis the order $k$ is fixed and the dimension $d$ is
large. The noise is then of order $\sqrt d$: $\sup X/\sqrt d$ tends almost surely to a constant, about $1.657$ for
$k=3$ \citep{subag2017complexity,jagannath2020statistical}. On this scale, with $\lambda/\sigma=\bar\lambda\sqrt d$,
the $t$-spacing test has power tending to one when $\bar\lambda$ is above a constant $C_{k}$, with $C_{3}=7.70$
(Theorem~\ref{thm:sqrtd}). For every odd $k$, $C_{k}\le1.72\cdot2^{k/2}\sqrt k$, and~$C_{k}/(2^{k/2}\sqrt k)\to1$ as
$k\to\infty$. So this constant grows exponentially with the order. It is a sufficient value of $\bar\lambda$, not a
threshold, and for $k=3$ it is far above the two values quoted below. We know of no other theorem on the power of spacing tests for a non-convex likelihood
when the dimension grows. For $k=3$ and a fixed $\bar\lambda<1.12$, no test of level $\alpha$ has power larger than
$\alpha+o(1)$ at every $e$ \citep[Theorem~3.2]{perry2020statistical}. The maximum test, which rejects for large
$\sup Z/\sigma$ with $\sigma$ known, has power tending to one when $\bar\lambda>1.21$
\citep{jagannath2020statistical}. We do not know the threshold of the $t$-spacing test (Remark~\ref{rem:dim}).

\paragraph{Related work} Tensor methods recast several latent variable models as the low-rank decomposition of
an empirical moment tensor of order $k\ge3$, which removes the identifiability obstructions of the matrix case
\citep{anandkumar2014tensor}; spiked tensor PCA, introduced by \citet{montanari2014statistical}, isolates the
low-rank-plus-noise structure of that problem. The matrix analogue is the Baik--Ben Arous--P\'ech\'e transition
\citep{baik2005phase}. \citet{zhang2018tensor} give the minimax rates of low-rank tensor estimation, and their
computational limits, as the dimension grows; Theorem~\ref{thm:estim} is at fixed $(k,d)$ and gives the
constant. On the Bayesian side, \citet{lesieur2017statistical} compute the asymptotic mutual information and the
minimum mean square error for a spike with independent entries drawn from a known prior, at fixed
signal-to-noise ratio as $d\to\infty$. \citet{lesieur2017statistical} and the works recalled in the previous paragraph
normalise the noise differently from us; Appendix~\ref*{app:norm} of the Supplementary Material translates their
thresholds into our variables. The critical points of the profile likelihood are counted in the spin-glass literature,
by \citet{auffinger2013} for the pure noise field and by \citet{arous2019landscape} for the spiked one. The
noise field $X$ is studied in \citet{kss}, from which the present paper takes, restated without proof in the
Supplementary Material, the layer decomposition of the helix process and a concentration bound
for the norm of the Gaussian Orthogonal Ensemble. The high-dimensional dynamics of Langevin and online
stochastic gradient descent for the multi-spiked model are characterised by
\citet{benarous2024langevin,benarous2024sgd}; these address the computational side of the problem,
complementary to the statistical guarantees here.

\paragraph{Notation} Throughout the paper $k\ge3$ is odd and $d\ge3$. $\Sph^{d-1}$ is the unit sphere of $\mathbb R^{d}$, $M=\Sph^{d-1}$ the parameter manifold,
$\dim M=d-1=:N$. We write $\rho^{2}=k/(2(k-1))$, $\varphi$ and $\Phi$ for the standard Gaussian density and
distribution function, $\lmax$ and $\lmin$ for the largest and smallest eigenvalues, $\|r\|_{F}$ and $\|r\|_{\rm op}$ for the Frobenius norm and the largest absolute eigenvalue
of a symmetric matrix $r$, and $\GOE(n)$ for the Gaussian Orthogonal Ensemble normalised by
$\Var(v^{\top}Gv)=1$ for unit $v$. Multi-indices $\ba\in\mathbb N^{d}$ with $|\ba|=k$ index the monomials
$t^{\ba}$, and $\binom k\ba=k!/\ba!$. For a function on $M$, $\nabla$ and $\nabla^{2}$ denote its gradient and
its Hessian on $M$.

\section{The model, and why the noise level must be estimated}\label{sec:model}

\subsection{The spiked tensor and its profile likelihood}

Write $\psi_{t}=t^{\otimes k}$, so that
\[\langle\psi_{s},\psi_{t}\rangle=\langle s,t\rangle^{k},\qquad \|\psi_{t}\|=1 ,\]
and $Z(t)=\langle Y,\psi_{t}\rangle=\lambda\langle e,t\rangle^{k}+\sigma X(t)$. The model is the Gaussian shift
$Y=\lambda\psi_{e}+\sigma W$ with a shift of norm exactly $\lambda$, so the signal-to-noise ratio is $\lambda/\sigma$.

\subsection{The second maximum and the normalised Hessian}

Take $\sigma=1$ in this subsection; the general case follows by homogeneity, $\lambda_{1}$, $\lambda_{2}$ and
$\Omega$ below being of degree one in $(\lambda,\sigma)$. The test is built on three statistics of the profile
likelihood: its maximum and the maximiser, which is the estimator of $e$; the second maximum, defined through a
regressed field; and the normalised Hessian at the maximiser:
\begin{subequations}\label{eq:stats}
\begin{align}
\lambda_{1}&=\max_{t\in M}Z(t)=Z(\hat t),\qquad \hat t=\operatorname*{arg\,max}_{t\in M}Z(t),\label{eq:argmax}\\
Z^{|t}(s)&=\frac{Z(s)-\langle s,t\rangle^{k}Z(t)-\langle s,t\rangle^{k-1}\langle s,\nabla Z(t)\rangle}{1-\langle s,t\rangle^{k}},
\qquad s\ne t,\label{eq:helix}\\
\lambda_{2}&=\sup_{s\ne\hat t}Z^{|\hat t}(s),\label{eq:lam2}\\
\Omega&=\Omega_{\hat t},\qquad \Omega_{t}=\frac{\nabla^{2}Z(t)+k\,Z(t)\,\mathrm{Id}}{k}\quad\text{on the tangent space at $t$.}\label{eq:Omega}
\end{align}
\end{subequations}
The maximiser is almost surely unique (Appendix~\ref*{app:kr} of the Supplementary Material). The second statistic looks at what remains of $Z$
once the information carried by the maximiser is removed: for a fixed $t$, \eqref{eq:helix} regresses $Z(s)$ on
the value $Z(t)$ and the gradient $\nabla Z(t)$ and divides the residual by $1-\langle s,t\rangle^{k}$. This is the
\emph{helix process} of \citet{azais2024second}; under the null it is independent of $(Z(t),\nabla Z(t))$. The
\emph{second maximum}~\eqref{eq:lam2} is its largest value at the maximiser. The
\emph{normalised Hessian}~\eqref{eq:Omega} of \citet{azais2024second} is the Hessian of $Z$ at the maximiser corrected by its value.

\medskip

At a critical point $t$ of $Z$ with value $x=Z(t)$ and $c=\langle e,t\rangle$, the Hessian is
$\nabla^{2}Z(t)=\sqrt{2k(k-1)}\,(G_{\lambda}-\rho x\,\mathrm{Id})$ with $G\sim\GOE(N)$ and
\begin{subequations}\label{eq:hessian}
\begin{align}
G_{\lambda}&=G+\beta_{2}(c)\,e_{\perp}e_{\perp}^{\top},\qquad
\beta_{2}(c)=\lambda\binom k2^{1/2}c^{k-2}(1-c^{2}),\qquad e_{\perp}=\frac{e-ct}{\sqrt{1-c^{2}}},\label{eq:spike}\\
\intertext{a rank-one spiked Gaussian Orthogonal Ensemble (Appendix~\ref*{app:layers}); $\beta_{2}(c)$ is the spike
of the Hessian at $t$. The shift $\rho x$ cancels in $\Omega_{t}$, and}
\Omega_{t}&=\frac{G_{\lambda}}{\rho},
\qquad \det\big(u\,\mathrm{Id}-\Omega_{t}\big)=\rho^{-N}\det\big(\rho u\,\mathrm{Id}-G_{\lambda}\big).\label{eq:omega}
\end{align}
\end{subequations}
\emph{The determinant in the spacing statistic is the determinant in the Kac--Rice integrand}, the formula that
counts the critical points of $Z$ (Appendix~\ref*{app:kr}). Two
order relations follow. At a local maximum $\nabla^{2}Z\prec0$, so $\lmax(G_{\lambda})<\rho x$, that is
$\lmax(\Omega)<\lambda_{1}$; and the limit of the helix as $s\to\hat t$ along a unit tangent direction $v$ is
$v^{\top}\Omega v$, so $\lambda_{2}\ge\lmax(\Omega)$. Since $\nabla Z(\hat t)=0$, \eqref{eq:helix} gives
$Z(s)-\lambda_{1}=(1-\langle s,\hat t\rangle^{k})(Z^{|\hat t}(s)-\lambda_{1})$, with $1-\langle s,\hat t\rangle^{k}>0$
for $s\ne\hat t$ as $k$ is odd. As $Z(s)\le\lambda_{1}$, this gives $\lambda_{2}\le\lambda_{1}$. Hence
\begin{equation}\label{eq:order}
\lmax(\Omega)\ \le\ \lambda_{2}\ \le\ \lambda_{1} .
\end{equation}
For $t\in M$ let $M_{t}=\Sph^{d-1}\setminus\{\pm t\}$ and let $\zeta^{t}$ be the helix~\eqref{eq:helix} of the noise $X$.
The second maximum at $t$ is
\[
\lambda_{2}^{t}=\sup_{M_{t}}Z^{|t}.
\]
As $Z^{|t}$ is continuous at $-t$, it is also $\sup_{s\ne t}Z^{|t}(s)$, so
$\lambda_{2}^{\hat t}=\lambda_{2}$.
The helix process is made precise in Appendix~\ref*{app:layers} of the Supplementary Material, where statements
and equations are numbered by appendix (Lemma~B.1, equation~(B.1a)).

\subsection{The noise level is not determined by the observation}\label{sub:rank}

The kernel of $X$ has finite rank. Write $f_{\ba}(t)=\binom k\ba^{1/2}t^{\ba}$ for the multi-indices
$|\ba|=k$. The multinomial theorem gives
\begin{subequations}\label{eq:rank}
\begin{align}
\langle s,t\rangle^{k}&=\sum_{|\ba|=k}f_{\ba}(s)f_{\ba}(t),\qquad \sum_{|\ba|=k}f_{\ba}(t)^{2}=1,\label{eq:rank-kernel}\\
X&=\sum_{|\ba|=k}\zeta_{\ba}f_{\ba},\qquad (\zeta_{\ba})\ \text{independent standard Gaussian variables,}\label{eq:rank-KL}
\end{align}
\end{subequations}
and the $m=\binom{d+k-1}{k}$ functions $f_{\ba}$ are linearly independent on $\Sph^{d-1}$, since a homogeneous
polynomial of degree $k$ that vanishes on the sphere vanishes everywhere, $p(x)=|x|^{k}p(x/|x|)$. So the
observation $(\sigma X(t))_{t\in M}$ determines the $m$ Gaussian variables $\sigma\zeta_{\ba}$ and nothing more.
We identify the symmetric tensors of order $k$, a Euclidean space of dimension~$m$ for $\langle\cdot,\cdot\rangle$,
with $\mathbb R^{m}$ through an orthonormal basis; by~\eqref{eq:rank-KL}, $W$ is a standard Gaussian vector of
this space, and the model is the finite-dimensional regression
\begin{equation}\label{eq:regmodel}
Y\sim N_{m}(\lambda\psi_{e},\sigma^{2}\mathrm{Id}_{m}),
\end{equation}
whose mean lies on the cone $\{\lambda\psi_{t}:\lambda\ge0,\ t\in\Sph^{d-1}\}$, of dimension $d$. The parameter
$\sigma$ is identifiable, the laws $N_{m}(0,\sigma^{2}\mathrm{Id}_{m})$ being distinct for distinct $\sigma$, but
it is not a function of the observation: these laws are mutually absolutely continuous, so no measurable function
of~$Y$ equals $\sigma$ almost surely, and no estimator of $\sigma$ is exact.

\medskip

The estimator of \citet{azais2024second} is built on the regression of $Y$ on~$(Z(t),\nabla Z(t))$ at a fixed
$t$. The regression spends $d$ of the $m$ coordinates, and the regressed field $Z^{|t}$ keeps the other $\kappa$
coordinates $(\bar\zeta_{i})_{i\le\kappa}$ of its Karhunen--Lo\`eve expansion; the estimator is their mean square.
It is independent in this sense: under the null, it is independent of the pair $(Z(t),\nabla Z(t))$ on which the
regression is made. It also has a closed form in terms of~$Y$.

\begin{proposition}[The residual estimator of the noise level]\label{prop:sigmahat}
Fix $t\in M$. Then
\begin{subequations}\label{eq:sigmahat}
\begin{align}
\kappa&=m-d,\label{eq:kappa}\\
\hat\sigma_{t}^{2}&=\frac1\kappa\sum_{i\le\kappa}\bar\zeta_{i}^{2}\ \sim\ \sigma^{2}\,\chi^{2}_{\kappa}/\kappa
\quad\text{under the null, independent of $(Z(t),\nabla Z(t))$,}\label{eq:sigmahat-law}\\
\kappa\,\hat\sigma_{t}^{2}&=\|Y\|^{2}-Z(t)^{2}-\frac1k\,\|\nabla Z(t)\|^{2}
\quad\text{under the null and under every alternative.}\label{eq:sigmahat-norm}
\end{align}
\end{subequations}
\end{proposition}

The proof is in Appendix~\ref*{app:proofs} of the Supplementary Material: \eqref{eq:kappa}
and~\eqref{eq:sigmahat-law} are taken from \citet{azais2024second}, and~\eqref{eq:sigmahat-norm} follows from the
orthonormality of $\psi_{t}$ and of the derivatives of $t\mapsto\psi_{t}$.

\medskip

Since $\kappa$ is an integer fixed by $(k,d)$, the known-variance test is not a limiting case of the~$t$-test
\emph{within} the model: the Student law tends to the Gaussian only as $m\to\infty$, that is across models. The
two tests are two different procedures, and only the second is exact without external information.

\subsection{The \texorpdfstring{$t$}{t}-spacing test in closed form}\label{sub:closed}

Let $T_{i}=\lambda_{i}/\hat\sigma$ with $\hat\sigma=\hat\sigma_{\hat t}$, which is positive almost surely since
$\kappa\hat\sigma^{2}=\|Y\|^{2}-\lambda_{1}^{2}>0$ by~\eqref{eq:sigmahat-norm}, as $d<m$. For a symmetric matrix $r$,
\begin{equation}\label{eq:Tdef}
H_{r}(\ell)=\int_{\ell}^{\infty}\det\big(u\,\mathrm{Id}-r\big)\,f_{m-1}\Big(u\sqrt{\tfrac{m-1}{\kappa}}\Big)\,\mathrm{d}u,
\qquad
\mathbf T=\frac{H_{\Omega/\hat\sigma}(T_{1})}{H_{\Omega/\hat\sigma}(T_{2})},
\end{equation}
$f_{\nu}$ being the Student density with $\nu$ degrees of freedom. By \citet[Theorem~3]{azais2024second},
$\mathbf T$ is uniform on $(0,1)$ under the null, whatever $\sigma$; the test of level $\alpha$ rejects when
$\mathbf T\le\alpha$. With $\varphi$ in place of the Student density, the known-variance statistic is
\begin{equation}\label{eq:Sdef}
G_{r}(\ell)=\int_{\ell}^{\infty}\det\big(u\,\mathrm{Id}-r\big)\,\varphi(u)\,\mathrm{d}u,
\qquad
\mathbf S=\frac{G_{\Omega/\sigma}(\lambda_{1}/\sigma)}{G_{\Omega/\sigma}(\lambda_{2}/\sigma)};
\end{equation}
it is likewise uniform on $(0,1)$ under the null when $\sigma$ is known \citep[Theorem~2]{azais2024second}.
Both $G_{r}$ and $H_{r}$ are finite sums of
incomplete moments of the Gaussian and of the Student densities, computed in closed form in
Appendix~\ref*{app:closed} of the Supplementary Material. Both are positive and strictly decreasing on
$(\lmax(r),\infty)$, with limit $0$ at infinity, the integrand being positive there.

\medskip

Recall that $k$ is odd. For $k$ even, $Z(-t)=Z(t)$ and $\langle\psi_{s},\psi_{t}\rangle=1$ at $s=-t$, so the assumptions of
\citet{azais2024second} fail on $\Sph^{d-1}$. We do not treat this case.

The statistic $\mathbf T$ is invariant under $Y\mapsto aY$, $a>0$, and $\mathbf S$ under
$(Y,\sigma)\mapsto(aY,a\sigma)$, so the laws of both depend on $(\lambda,\sigma)$ only through $\lambda/\sigma$; from
Section~\ref{sec:reverse} on we take $\sigma=1$ and call $\lambda$ the signal-to-noise ratio. Both statistics are also invariant under the orthogonal group. For $Q\in O(d)$, the field of $Q^{\otimes k}Y$ is
$Z(Q^{\top}\cdot)$ and $\|Q^{\otimes k}Y\|=\|Y\|$. The two maxima, the spectrum of the normalised Hessian at the
maximiser and the noise estimate do not change, so both statistics take the same value at $Y$ and at
$Q^{\otimes k}Y$. Also $Q^{\otimes k}Y=\lambda\psi_{Qe}+\sigma Q^{\otimes k}W$, and $Q^{\otimes k}W$ has the law of
$W$, since $X(Q^{\top}\cdot)$ has the covariance $\langle s,t\rangle^{k}$ of $X$. So the law of each statistic is
the same at $e$ and at $Qe$. The group acts transitively on the sphere, so these laws do not depend on $e$. We
write $\PP_{\lambda}$ for the law of the observation at signal $\lambda$ and any fixed direction $e$.

\section{The rejection region and the power}\label{sec:reverse}

The level $\alpha\in(0,1)$ is fixed. The rejection region $\{\mathbf T\le\alpha\}$ is sandwiched between two ratio
tests on $\lambda_{1}/\lambda_{2}$ (Proposition~\ref{prop:sandwich}), whose exponent $\kappa$ is the number of degrees
of freedom of $\hat\sigma$. On the event where the test fails, $\lambda_{1}$ is at most a constant times the noise
estimate $\hat\sigma$ (Corollary~\ref{cor:def}). The first fact bounds the power from above
(Theorem~\ref{thm:upper}). The second bounds the deficiency $1-\PP_{\lambda}\{\mathbf T\le\alpha\}$ from above at every
$\lambda$ (Corollary~\ref{cor:closed}), hence its exponent from below. One configuration of the normalised Hessian
and of the noise estimate gives the matching bound, so the exponent is identified (Theorem~\ref{thm:exponent}). The
last subsection lets the dimension grow (Theorem~\ref{thm:sqrtd}). All these statements are
scale-free, as the statistic itself must be. For the known-variance statistic $\mathbf S$ the region lies between a
shift of the squared level by $c_{\alpha}$ and a shift by $2\log(1/\alpha)-2\log(\lambda_{1}/\lambda_{2})$. Almost
surely it contains $\{\lambda_{1}^{2}\ge\lambda_{2}^{2}+c_{\alpha}\}$, with
$c_{\alpha}$ the $(1-\alpha)$-quantile of the $\chi^{2}_{2d}$ law (Corollary~\ref{cor:Sclosed}), and it is contained
in $\{\lambda_{1}^{2}\ge\lambda_{2}^{2}+2\log(1/\alpha)-2\log(\lambda_{1}/\lambda_{2})\}$ (Theorem~\ref{thm:upper}).
The difference between a dilation and a shift is the price of not knowing $\sigma$.

\medskip

The reason is the tail of the reference density. In $H_{r}$ the weight $\det(u\,\mathrm{Id}-r)$ has
degree~$N$ and the Student density with $\nu=m-1$ degrees of freedom decays like $u^{-(\nu+1)}$, so the integrand
decays like $u^{-\kappa-1}$ and $H_{r}$ like $u^{-\kappa}$, since $\nu-N=(m-1)-(d-1)=\kappa$: multiplying the level by $\vartheta$ divides
$H_{r}$ by $\vartheta^{\kappa}$. In $G_{r}$ the Gaussian factor turns a dilation of the level into a shift of half
its square. The next lemma makes both statements non-asymptotic. Write $\lmax(r)^{+}=\max(0,\lmax(r))$ and
$\lmin(r)^{-}=\max(0,-\lmin(r))$, and recall from~\eqref{eq:Tdef} that the reference density in
$H_{r}$ is $f_{\nu}(u/\gamma)$ with the scale $\gamma=\sqrt{\kappa/\nu}$, so it is proportional
to~$(\kappa+u^{2})^{-(\nu+1)/2}$, and $\gamma^{2}\nu=\kappa$ is exactly the number of degrees of freedom of $\hat\sigma$.

\begin{lemma}[Dilation bounds]\label{lem:ratio}
Let $r$ be a symmetric matrix of size $N$ and let $\lmax(r)^{+}<\ell_{2}\le\ell_{1}$; put $\vartheta=\ell_{1}/\ell_{2}$
and $g=\tfrac12(\ell_{1}^{2}-\ell_{2}^{2})$. Then
\begin{subequations}
\begin{equation}\label{eq:dilH}
\vartheta^{-\kappa}\Big(1+\frac{\lmin(r)^{-}}{\ell_{2}}\Big)^{-N}\ \le\
\frac{H_{r}(\ell_{1})}{H_{r}(\ell_{2})}\ \le\
\vartheta^{N+1}\Big(1-\frac{\lmax(r)^{+}}{\ell_{2}}\Big)^{-N}
\Big(\frac{\kappa+\ell_{2}^{2}}{\kappa+\ell_{1}^{2}}\Big)^{\frac{\nu+1}{2}},
\end{equation}
\begin{equation}\label{eq:dilG}
\frac{\ell_{2}}{\ell_{1}}\,e^{-g}\ \le\ \frac{G_{r}(\ell_{1})}{G_{r}(\ell_{2})}\ \le\
e^{-g}\Big(1+\frac{g}{\ell_{2}\big(\ell_{2}-\lmax(r)^{+}\big)}\Big)^{N}.
\end{equation}
\end{subequations}
The two bounds in~\eqref{eq:dilH} differ by a factor tending to $1$ as $\ell_{2}\to\infty$ at fixed $\vartheta$ and
fixed $r$, and likewise in~\eqref{eq:dilG} at fixed $g$.
\end{lemma}

The proof is in Appendix~\ref*{app:proofs} of the Supplementary Material.

Lemma~\ref{lem:ratio} is applied below at $r=\Omega/\hat\sigma$ and $\ell_{2}=T_{2}$, which requires
$T_{2}>\lmax(\Omega/\hat\sigma)^{+}$. This is a property of the field at the maximiser $\hat t$, a random point.
The next lemma reduces such a property to the same property at each fixed point $t$, conditionally on
$\nabla Z(t)=0$, because $\hat t$ is a critical point of $Z$ and the Kac--Rice formula of Appendix~\ref*{app:kr}
counts them.

\begin{lemma}[Transfer to the argmax]\label{lem:transfer}
Let $\lambda\ge0$ and let $f:M\times\mathbb R^{m}\to\mathbb R$ be continuous; write $f_{t}=f(t,Y)$. For $t\in M$, let
$\Pi_{t}$ be the orthogonal projection of $\mathbb R^{m}$ onto the kernel of the linear map $Y\mapsto\nabla Z(t)$. By
Gaussian regression, the law of $Y$ given $\nabla Z(t)=0$ is the law of $\Pi_{t}Y$, that is
\[
\PP_{\lambda}\{f_{t}=0\mid\nabla Z(t)=0\}\ =\ \PP_{\lambda}\{f(t,\Pi_{t}Y)=0\}.
\]
If this probability is $0$ for every $t\in M$, then $f_{\hat t}\ne0$ almost surely under $\PP_{\lambda}$.
\end{lemma}

\begin{proof}
The maximiser $\hat t$ is a critical point of $Z$. So $\PP_{\lambda}\{f_{\hat t}=0\}$ is at most the expected
number of critical points $t$ of $Z$ with $f_{t}=0$. By~(\ref*{eq:markedkr}) for the closed condition
$\{f_{t}=0\}$, with~$B=M$, and by monotone convergence as $u\to-\infty$, this expected number is
\[
\int_{M}\E\Big[\big|\det\nabla^{2}Z(t)\big|\,\ind\{f_{t}=0\}\ \Big|\ \nabla Z(t)=0\Big]\,p_{\nabla Z(t)}(0)\,\mathrm{d}t\ =\ 0,
\]
since $\ind\{f_{t}=0\}$ vanishes almost surely given $\nabla Z(t)=0$.
\end{proof}

\begin{lemma}[The second maximum is positive and exceeds the Hessian]\label{lem:pos}
Let $k$ be odd. For every $Y$ and every $t\in M$, $Z(-t)=-Z(t)$ and $Z^{|t}(-t)=0$. For every $\lambda\ge0$,
$\lambda_{2}>\max\big(\lmax(\Omega),0\big)$ almost surely under $\PP_{\lambda}$.
\end{lemma}

\begin{proof}
As $k$ is odd, $\psi_{-t}=-\psi_{t}$, so $Z(-t)=\langle Y,\psi_{-t}\rangle=-Z(t)$. At $s=-t$ in~\eqref{eq:helix},
$\langle s,t\rangle^{k}=-1$ and $\langle s,\nabla Z(t)\rangle=0$, as $\nabla Z(t)$ is tangent to $M$ at $t$. So
$Z^{|t}(-t)=\big(Z(-t)+Z(t)\big)/2=0$, and $\lambda_{2}^{t}\ge0$ for every $t$.

Let $f_{t}=f(t,Y)=\lambda_{2}^{t}-\lmax(\Omega_{t})^{+}$. The function $f$ is continuous on
$M\times\mathbb R^{m}$ by Lemma~\ref*{lem:lam2cont}. By~\eqref{eq:order} and
$\lambda_{2}\ge0$, $f_{\hat t}=\lambda_{2}-\lmax(\Omega)^{+}\ge0$. So the claim is that $f_{\hat t}\ne0$ almost surely
under $\PP_{\lambda}$. The event $\{f_{\hat t}=0\}$ is determined by $Y$, and the laws of $Y$ under $\PP_{\lambda}$
and under $\PP_{0}$ are equivalent (Appendix~\ref*{app:kr}). So it suffices to take $\lambda=0$.

Fix $t\in M$. Under $\PP_{0}$, $Z=\sigma X$, so $Z^{|t}=\sigma\zeta^{t}$ and $\Omega_{t}=\sigma G/\rho$, with $G$ the
matrix~(\ref*{eq:G}). By Lemma~\ref*{lem:poles}(ii), $\zeta^{t}$ extends to an almost surely continuous field
$\tilde\zeta$ on the blow-up $[0,\pi]\times\Sph^{d-2}$ of the sphere at $\pm t$, with
$\tilde\zeta(0,v)=v^{\top}Gv/\rho$ and $\tilde\zeta(\pi,\cdot)\equiv0$. So the maxima of $\sigma\tilde\zeta$ on the two
pole fibres $\theta=0$ and $\theta=\pi$ are $\lmax(\Omega_{t})$ and $0$. The blow-up is compact, so the supremum of
$\tilde\zeta$ is attained. By \citet[Proposition~4.1 of arXiv:2611.09459v1]{kss}, almost surely it is not attained
on a pole fibre. It is then attained at a point of $M_{t}$, where $\sigma\tilde\zeta=Z^{|t}$. So $\lambda_{2}^{t}$ is
the supremum of $\sigma\tilde\zeta$ and it exceeds the maxima on the two pole fibres. Hence $\PP_{0}\{f_{t}=0\}=0$.
By Lemma~2 of \citet{azais2024second}, under $\PP_{0}$ the helix $\zeta^{t}$ and the matrix $G$ are independent of
$(X(t),\nabla X(t))$, as recalled in Appendix~\ref*{app:layers} at~(\ref*{eq:G}) and after the layer
decomposition. All these variables are jointly Gaussian, so the pair $(\zeta^{t},G)$ is independent of
$(X(t),\nabla X(t))$, and the regression on $\nabla Z(t)$ does not change its law. As $f_{t}$ is a function of $(\zeta^{t},G)$ and $\nabla Z(t)=\sigma\nabla X(t)$,
$\PP_{0}\{f_{t}=0\mid\nabla Z(t)=0\}=\PP_{0}\{f_{t}=0\}=0$ for every $t$. Lemma~\ref{lem:transfer} above, with
$\lambda=0$, gives $f_{\hat t}\ne0$ almost surely under $\PP_{0}$.
\end{proof}

\begin{proposition}[The rejection region is sandwiched between two ratio tests]\label{prop:sandwich}
Almost surely, with $\vartheta=\lambda_{1}/\lambda_{2}\ge1$,
\begin{equation}\label{eq:sandwich}
\vartheta^{-\kappa}\Big(1+\frac{\lmin(\Omega)^{-}}{\lambda_{2}}\Big)^{-N}\ \le\ \mathbf T\ \le\
\vartheta^{-\kappa}\Big(1-\frac{\lmax(\Omega)^{+}}{\lambda_{2}}\Big)^{-N}
\Big(1+\frac{\kappa\hat\sigma^{2}}{\lambda_{2}^{2}}\Big)^{\frac{\nu+1}{2}} .
\end{equation}
The left-hand side involves neither $\hat\sigma$ nor the scale of the observation; the right-hand side involves
$\hat\sigma$ only through $\kappa\hat\sigma^{2}/\lambda_{2}^{2}$.
\end{proposition}

\begin{proof}
Apply Lemma~\ref{lem:ratio} with $r=\Omega/\hat\sigma$ and $\ell_{i}=T_{i}=\lambda_{i}/\hat\sigma$. Its hypothesis
$\lmax(r)^{+}<\ell_{2}$ holds almost surely: by Lemma~\ref{lem:pos}, $T_{2}>0$ and
$\lmax(\Omega/\hat\sigma)<T_{2}$, as $\hat\sigma>0$. All three quantities~$\ell_{1}/\ell_{2}$, $\lmin(r)^{-}/\ell_{2}$ and
$\lmax(r)^{+}/\ell_{2}$ are unchanged when $r$ and the $\ell_{i}$ are divided by $\hat\sigma$, which gives the
$\hat\sigma$-free form of $\vartheta$ and of the two Hessian corrections. For the last factor,
\[\Big(\frac{\kappa+\ell_{2}^{2}}{\kappa+\ell_{1}^{2}}\Big)^{\frac{\nu+1}{2}}
\le\Big(\frac{\kappa+\ell_{2}^{2}}{\ell_{1}^{2}}\Big)^{\frac{\nu+1}{2}}
=\vartheta^{-(\nu+1)}\Big(1+\frac{\kappa}{\ell_{2}^{2}}\Big)^{\frac{\nu+1}{2}},\]
and $\vartheta^{N+1}\vartheta^{-(\nu+1)}=\vartheta^{-\kappa}$, while $\kappa/\ell_{2}^{2}=\kappa\hat\sigma^{2}/\lambda_{2}^{2}$.
\end{proof}

\subsection{An upper bound on the power}

The left-hand inequality of~\eqref{eq:sandwich} costs nothing and yields an upper bound on the power.

\begin{theorem}[Upper bound on the power]\label{thm:upper}
Let $\lambda\ge0$, $\sigma>0$ and $\alpha\in(0,1)$. Then
\begin{subequations}
\begin{equation}\label{eq:upper}
\PP_{\lambda}\big\{\mathbf T\le\alpha\big\}\ \le\
\PP_{\lambda}\bigg\{\lambda_{2}\ \le\ \alpha^{1/\kappa}\,\lambda_{1}
\Big(1+\frac{\lmin(\Omega)^{-}}{\lambda_{2}}\Big)^{N/\kappa}\bigg\},
\end{equation}
equivalently
\begin{equation}
1-\PP_{\lambda}\big\{\mathbf T\le\alpha\big\}\ \ge\
\PP_{\lambda}\bigg\{\lambda_{2}\ >\ \alpha^{1/\kappa}\,\lambda_{1}
\Big(1+\frac{\lmin(\Omega)^{-}}{\lambda_{2}}\Big)^{N/\kappa}\bigg\} .
\end{equation}
\end{subequations}
For the known-variance test, $\PP_{\lambda}\{\mathbf S\le\alpha\}\le
\PP_{\lambda}\{(\lambda_{1}^{2}-\lambda_{2}^{2})/\sigma^{2}\ge2\log(1/\alpha)-2\log(\lambda_{1}/\lambda_{2})\}$. Both bounds
hold under the null as well, where the left-hand side is exactly $\alpha$.
\end{theorem}

\begin{proof}
On $\{\mathbf T\le\alpha\}$ the left-hand inequality of~\eqref{eq:sandwich} gives
$\vartheta^{-\kappa}(1+\lmin(\Omega)^{-}/\lambda_{2})^{-N}\le\alpha$; raising to the power $1/\kappa$ and
multiplying by $\lambda_{1}(1+\lmin(\Omega)^{-}/\lambda_{2})^{N/\kappa}$, which turns $\vartheta^{-1}\lambda_{1}$ into
$\lambda_{2}$, gives the event on the right of~\eqref{eq:upper}, and the complement is the second display. For~$\mathbf S$, use the left-hand inequality of~\eqref{eq:dilG} in the same way.
\end{proof}

Under the null the bound is not vacuous and can be read as a statement about the null law of the triple: since
$\PP_{0}\{\mathbf T\le\alpha\}=\alpha$ exactly, \eqref{eq:upper} says that
$\lambda_{1}/\lambda_{2}\ge\alpha^{-1/\kappa}\big(1+\lmin(\Omega)^{-}/\lambda_{2}\big)^{-N/\kappa}$ with
probability at least $\alpha$. The correction lowers the threshold, so it does not by itself bound
$\PP_{0}\{\lambda_{1}/\lambda_{2}\ge\alpha^{-1/\kappa}\}$. At
$(k,d)=(3,3)$ and $\alpha=0.05$ one has $\kappa=7$ and $\alpha^{-1/\kappa}=1.534$; at $(k,d)=(5,5)$,
$\kappa=121$ and $\alpha^{-1/\kappa}=1.025$. The exponent $1/\kappa$ is what makes the $t$-version
demanding when $m$ is small: with few degrees of freedom the test must see a large \emph{ratio} of the two
maxima, not merely a large gap.

\subsection{A lower bound on the power}

A lower bound on the power needs an upper bound on~$\mathbf T$. It comes from Lemmas~\ref{lem:T2sharp}
and~\ref{lem:mono}, through Corollary~\ref{cor:def} below. They rest on a decomposition of $P_{t}^{\perp}Y$ into layers,
one for each power of the tangent direction.

\begin{lemma}[Layers at a point]\label{lem:layers}
Fix $t\in M$ and identify the tangent space at $t$ with~$t^{\perp}$. For a unit vector $v\perp t$ and $0\le j\le k$
let $\Psi_{j}(v)=\binom kj^{1/2}\operatorname{Sym}\big(t^{\otimes(k-j)}\otimes v^{\otimes j}\big)$, where
$\operatorname{Sym}$ averages over the permutations of the $k$ factors, let $\mathcal V_{j}$ be the span of
the~$\Psi_{j}(v)$, and let $w_{j}(\theta)=\binom kj^{1/2}\cos^{k-j}\theta\,\sin^{j}\theta/(1-\cos^{k}\theta)$, the weights of
Appendix~\ref*{app:layers}. Then~$\langle\Psi_{j}(v),\Psi_{j}(v')\rangle=\langle v,v'\rangle^{j}$,
$\mathcal V_{j}\perp\mathcal V_{j'}$ for $j\ne j'$, the range of $P_{t}^{\perp}$ is
$\mathcal V_{2}\oplus\dots\oplus\mathcal V_{k}$, and the map $vv^{\top}\mapsto\Psi_{2}(v)$ extends to a linear isometry
$\iota$ from the symmetric matrices on $t^{\perp}$, with the Frobenius norm, onto $\mathcal V_{2}$. Let $P_{\mathcal V_{j}}$ be the orthogonal projection onto $\mathcal V_{j}$. For every $Y$ and
every~$s=\cos\theta\,t+\sin\theta\,v\in M_{t}$,
\[
\begin{gathered}
Z^{|t}(s)=\sum_{j=2}^{k}w_{j}(\theta)\,\langle Y,\Psi_{j}(v)\rangle,\qquad
v^{\top}\Omega_{t}v=\rho^{-1}\langle Y,\Psi_{2}(v)\rangle,\\
\kappa\hat\sigma_{t}^{2}=\rho^{2}\|\Omega_{t}\|_{F}^{2}+\textstyle\sum_{j=3}^{k}\|P_{\mathcal V_{j}}Y\|^{2}.
\end{gathered}
\]
\end{lemma}

The proof is in Appendix~\ref*{app:proofs} of the Supplementary Material.

\begin{lemma}[The second maximum and the spectrum of the Hessian]\label{lem:T2sharp}
For every $Y$ and every $t$,
$\rho^{2}(\lambda_{2}^{t})^{2}+\rho^{2}\big(\|\Omega_{t}\|_{F}^{2}-\|\Omega_{t}\|_{\rm op}^{2}\big)\le\kappa\hat\sigma_{t}^{2}$.
In particular, with $T_{\max}=\sqrt\kappa/\rho$ and $r=\Omega/\hat\sigma$, almost surely
$T_{2}^{2}+\|r\|_{F}^{2}-\|r\|_{\rm op}^{2}\le T_{\max}^{2}$ and $T_{2}\le T_{\max}$.
\end{lemma}

\begin{proof}
Let $\Delta_{t}^{2}=\kappa\hat\sigma_{t}^{2}-\rho^{2}\|\Omega_{t}\|_{F}^{2}=\sum_{j\ge3}\|P_{\mathcal V_{j}}Y\|^{2}$. By
Lemma~\ref{lem:layers}, $|\langle Y,\Psi_{2}(v)\rangle|\le\rho\|\Omega_{t}\|_{\rm op}$, and
$|\langle Y,\Psi_{j}(v)\rangle|\le\|P_{\mathcal V_{j}}Y\|$ for $j\ge3$, as $\|\Psi_{j}(v)\|=1$. The Cauchy--Schwarz
inequality, used twice, gives for $s\in M_{t}$
\[
|Z^{|t}(s)|\ \le\ |w_{2}(\theta)|\,\rho\|\Omega_{t}\|_{\rm op}+\Big(\sum_{j\ge3}w_{j}(\theta)^{2}\Big)^{1/2}\Delta_{t}
\ \le\ V(\theta)^{1/2}\big(\rho^{2}\|\Omega_{t}\|_{\rm op}^{2}+\Delta_{t}^{2}\big)^{1/2},
\]
with $V=\sum_{j\ge2}w_{j}^{2}\le\rho^{-2}$ by Lemma~\ref*{lem:variance}. Take the supremum over $M_{t}$ and square. At
$t=\hat t$ divide by $\rho^{2}\hat\sigma^{2}$, which is positive almost surely, and use $\|r\|_{F}\ge\|r\|_{\rm op}$.
\end{proof}

The bound $T_{\max}$ on $T_{2}$ is $3.06$ at $(k,d)=(3,3)$, $6.32$ at $(3,5)$ and $13.9$ at $(5,5)$. On the
simulations of Section~\ref{sec:num} it is never exceeded, and it is approached on the replicates where the test
fails. With the next lemma, Lemma~\ref{lem:T2sharp} describes the deficiency event without reference
to~$\lambda_{2}$ or to the Hessian.

\begin{lemma}[Monotonicity in the spectrum]\label{lem:mono}
Let $r$ be a symmetric matrix of size $N$ and~$\lmax(r)\le\ell_{2}\le\ell_{1}$. The ratio $H_{r}(\ell_{1})/H_{r}(\ell_{2})$
is nondecreasing in $\ell_{2}$, and in each eigenvalue of $r$ on $(-\infty,\ell_{2}]$, the other eigenvalues being
fixed.
\end{lemma}

The proof is in Appendix~\ref*{app:proofs} of the Supplementary Material.

For a symmetric matrix $r$ of size $N$ and $v\ge\lmax(r)$, the function $H_{r}$ is positive and continuous on
$[v,\infty)$, strictly decreasing, with limit $0$. So there is a unique $u_{\alpha}(r,v)>v$ with
\begin{subequations}\label{eq:const}
\begin{align}
H_{r}\big(u_{\alpha}(r,v)\big)&=\alpha\,H_{r}(v),\label{eq:ualpha}\\
\intertext{and, when $T_{2}=v$ and $\Omega/\hat\sigma=r$, the test rejects if and only if $T_{1}\ge u_{\alpha}(r,v)$.
Let}
A&=\frac1{T_{\max}}\max\Big\{u_{\alpha}\big(r,\lmax(r)\big):\ r\succeq0,\ \|r\|_{F}=T_{\max}\Big\}.\label{eq:Adef}
\end{align}
\end{subequations}
The maximum exists: $H_{r}(u)-\alpha H_{r}(\lmax(r))$ is continuous in $(r,u)$ and strictly decreasing in $u$, so
$u_{\alpha}(r,\lmax(r))$ is a continuous function of the eigenvalues of $r$, which range over a compact set. A matrix of rank one in this set has $\lmax(r)=T_{\max}$, so $A>1$. The values of $A$ are in
Remark~\ref{rem:bracket}.

\begin{corollary}[The deficiency event]\label{cor:def}
Let $A$ be the constant~\eqref{eq:Adef}. Almost surely,
\[\{\mathbf T>\alpha\}\ \subseteq\ \Big\{\lambda_{1}<\frac{A\sqrt\kappa}{\rho}\,\hat\sigma\Big\}.\]
\end{corollary}

\begin{proof}
Almost surely the maximiser is unique and $\hat\sigma>0$. Take such an observation with~$\mathbf T>\alpha$, let
$r=\Omega/\hat\sigma$ have the eigenvalues $r_{1}\ge\dots\ge r_{N}$, and let $r_{i}^{+}=\max(r_{i},0)$. The
sum~$\|r\|_{F}^{2}-\|r\|_{\rm op}^{2}$ contains all the $r_{i}^{2}$ but one of largest modulus, so it is at least
$\sum_{i\ge2}(r_{i}^{+})^{2}$. Lemma~\ref{lem:T2sharp} gives
\[
T_{2}\ \le\ v':=\Big(T_{\max}^{2}-\sum_{i\ge2}(r_{i}^{+})^{2}\Big)^{1/2},
\]
and $r_{i}^{+}\le T_{2}\le T_{1}$ by~\eqref{eq:order} and Lemma~\ref{lem:pos}. Let $r''$ have the eigenvalues $v',r_{2}^{+},\dots,r_{N}^{+}$: it is
positive semidefinite, $\|r''\|_{F}=T_{\max}$ and $\lmax(r'')=v'$. If $T_{1}\le v'$, then $T_{1}<u_{\alpha}(r'',v')$.
Otherwise Lemma~\ref{lem:mono}, applied to the negative eigenvalues raised to $0$, then to $\ell_{2}$ raised
from~$T_{2}$ to $v'$, then to $r_{1}^{+}$ raised to $v'$, gives $\alpha<\mathbf T\le H_{r''}(T_{1})/H_{r''}(v')$, hence
$T_{1}<u_{\alpha}(r'',v')$, as~$H_{r''}$ decreases. In both cases $T_{1}<u_{\alpha}(r'',\lmax(r''))\le AT_{\max}$, that
is $\lambda_{1}<AT_{\max}\hat\sigma$.
\end{proof}

\begin{corollary}[A power bound in closed form]\label{cor:closed}
Let $\lambda\ge0$, let $A$ be the constant of Corollary~\ref{cor:def}, and let $\xi\sim N(0,1)$ and $\chi_{m-1}$ be
independent, where $\chi_{m-1}$ is the norm of a standard Gaussian vector of $\mathbb R^{m-1}$. Then
\begin{subequations}\label{eq:closed}
\begin{align}
1-\PP_{\lambda}\{\mathbf T\le\alpha\}&\ \le\ \PP\Big\{\lambda+\xi<\frac{A}{\rho}\,\chi_{m-1}\Big\},\label{eq:closed-nct}\\
\intertext{and, if $\rho\lambda\ge A\sqrt{m-1}$,}
1-\PP_{\lambda}\{\mathbf T\le\alpha\}&\ \le\
\exp\Big(-\frac{\big(\rho\lambda-A\sqrt{m-1}\big)^{2}}{2(A^{2}+\rho^{2})}\Big).\label{eq:closed-chernoff}
\end{align}
\end{subequations}
\end{corollary}

\begin{proof}[Proof of Corollary~\ref{cor:closed}]
In the coordinates of~\eqref{eq:regmodel}, with $\sigma=1$, $Y=\lambda\psi_{e}+W$ and~$W$ is standard Gaussian in
$\mathbb R^{m}$. The vectors collinear with some $\psi_{t}$ form a set of dimension $d<m$, so almost surely $Y$ is not
one of them, and $|Z(t)|<\|Y\|$ for every $t$ by the Cauchy--Schwarz inequality. As $\nabla Z(\hat t)=0$,
\eqref{eq:sigmahat-norm} at $t=\hat t$ gives $\kappa\hat\sigma^{2}=\|Y\|^{2}-\lambda_{1}^{2}$, and Corollary~\ref{cor:def}
gives $\{\mathbf T>\alpha\}\subseteq\{\rho\lambda_{1}<A(\|Y\|^{2}-\lambda_{1}^{2})^{1/2}\}$. Let
$\xi=\langle W,\psi_{e}\rangle$ and~$\chi=\|W-\xi\psi_{e}\|$. Then $Z(e)=\lambda+\xi$ and
$Y-Z(e)\psi_{e}=W-\xi\psi_{e}$, so $\|Y\|^{2}=Z(e)^{2}+\chi^{2}$. The map~$x\mapsto x(\|Y\|^{2}-x^{2})^{-1/2}$ is
increasing on $(-\|Y\|,\|Y\|)$, and $Z(e)\le\lambda_{1}$. Hence, on $\{\mathbf T>\alpha\}$,
$\rho Z(e)<A(\|Y\|^{2}-Z(e)^{2})^{1/2}=A\chi$, that is $\rho(\lambda+\xi)<A\chi$, whatever the sign of $\lambda+\xi$.
The variables $\xi$ and $W-\xi\psi_{e}$ are the projections of $W$ onto $\psi_{e}$ and onto its orthogonal
complement. So $\xi\sim N(0,1)$ is independent of $\chi$, and $\chi$ has the law of~$\chi_{m-1}$. This
proves~\eqref{eq:closed-nct}. For~\eqref{eq:closed-chernoff}, $A\chi-\rho\xi$ is the supremum of
$\langle Au-\rho\psi_{e},W\rangle$ over the unit vectors~$u$ orthogonal to $\psi_{e}$. This is a centred Gaussian
process with variance~$A^{2}+\rho^{2}$, and the mean of its supremum is $A\,\E\chi\le A(\E\chi^{2})^{1/2}=A\sqrt{m-1}$.
The concentration inequality for the supremum of a Gaussian process \citep{azais2009level} gives
$\PP\{A\chi-\rho\xi>\rho\lambda\}\le\exp\big(-(\rho\lambda-A\sqrt{m-1})^{2}/(2(A^{2}+\rho^{2}))\big)$ when
$\rho\lambda\ge A\sqrt{m-1}$.
\end{proof}

The right-hand side of~\eqref{eq:closed-nct} is the distribution function of the noncentral Student law with~$m-1$
degrees of freedom and noncentrality $\lambda$, at the point $A\sqrt{m-1}/\rho$. At $\alpha=0.05$ it is below
$0.05$ from about $\lambda=11.0$ at $(k,d)=(3,5)$, and from $\lambda=11.4$ at $(3,3)$.

The same route gives a power bound for the known-variance statistic $\mathbf S$, once the upper bound
of~\eqref{eq:dilG} is replaced by a bound free of the spectrum of $r$.

\begin{lemma}[A shift bound]\label{lem:Gratio}
Let $r$ be a symmetric matrix of size $N$, let $\lmax(r)\le\ell_{2}$ and~$0\le\ell_{2}\le\ell_{1}$, and put
$g=\tfrac12(\ell_{1}^{2}-\ell_{2}^{2})$. Then
\[\frac{G_{r}(\ell_{1})}{G_{r}(\ell_{2})}\ \le\ e^{-g}\sum_{j=0}^{N}\frac{g^{j}}{j!}.\]
\end{lemma}

\begin{proof}
For $g=0$ both sides are $1$; let $g>0$. Write $D(u)=\det(u\,\mathrm{Id}-r)=\prod_{i}(u-r_{i})$, the $r_{i}$ being
the eigenvalues of $r$. It is positive for $u>\ell_{2}$, since $r_{i}\le\ell_{2}$. As $\ell_{2}\ge0$, the
substitution $u=(\ell_{2}^{2}+2s)^{1/2}$ maps $s\in(0,\infty)$ onto $u\in(\ell_{2},\infty)$ and $s=g$ to
$u=\ell_{1}$, with~$\mathrm{d}u=\mathrm{d}s/u$ and $\varphi(u)=e^{-s}\varphi(\ell_{2})$. Let
$g_{N}(s)=s^{N}e^{-s}/N!$. Then, for $h\in\{0,g\}$,
\begin{subequations}
\begin{equation}
G_{r}\big((\ell_{2}^{2}+2h)^{1/2}\big)=N!\,\varphi(\ell_{2})\int_{h}^{\infty}g_{N}(s)\,q(s)\,\mathrm{d}s,
\qquad q(s)=\frac1u\prod_{i}\frac{u-r_{i}}{s}.
\end{equation}
The function $q$ is positive and nonincreasing: $1/u$ decreases, and $s$ increases with $u$, with
\begin{equation}
\frac{\mathrm{d}}{\mathrm{d}u}\log\frac{u-r_{i}}{s}=\frac1{u-r_{i}}-\frac{2u}{(u-\ell_{2})(u+\ell_{2})}
\ \le\ \frac1{u-\ell_{2}}-\frac1{u-\ell_{2}}=0,
\end{equation}
because $r_{i}\le\ell_{2}$ and $2u\ge u+\ell_{2}>0$. Put $a=\int_{g}^{\infty}g_{N}q$, $b=\int_{0}^{g}g_{N}q$ and
$\bar\Gamma=\int_{g}^{\infty}g_{N}$. As $q$ is nonincreasing, $a\le q(g)\bar\Gamma$ and $b\ge q(g)(1-\bar\Gamma)$.
The ratio $a/(a+b)$ increases in $a$ and decreases in $b$, so
\begin{equation}
\frac{G_{r}(\ell_{1})}{G_{r}(\ell_{2})}=\frac{a}{a+b}\ \le\ \frac{q(g)\bar\Gamma}{q(g)\bar\Gamma+q(g)(1-\bar\Gamma)}=\bar\Gamma,
\end{equation}
\end{subequations}
as the factor $q(g)>0$ cancels. Integration by parts gives
$\bar\Gamma=e^{-g}\sum_{j\le N}g^{j}/j!$.
\end{proof}

\begin{corollary}[A power bound for the known-variance test]\label{cor:Sclosed}
Let $\sigma=1$ and $\lambda\ge0$, let~$c_{\alpha}$ be the $(1-\alpha)$-quantile of the $\chi^{2}_{2d}$ law, and let
$\xi\sim N(0,1)$ and $\chi_{m-1}$ be independent, as in Corollary~\ref{cor:closed}. Then, almost surely,
$\{\mathbf S>\alpha\}\subseteq\{\lambda_{1}^{2}<\lambda_{2}^{2}+c_{\alpha}\}$, and
\begin{subequations}\label{eq:Sclosed}
\begin{align}
1-\PP_{\lambda}\{\mathbf S\le\alpha\}&\ \le\ \PP\big\{\rho^{2}(\lambda+\xi)^{2}<\chi_{m-1}^{2}+\rho^{2}c_{\alpha}\big\},\label{eq:Sclosed-nct}\\
\intertext{and, if $\rho\lambda\ge\sqrt{m-1}+\rho\sqrt{c_{\alpha}}$,}
1-\PP_{\lambda}\{\mathbf S\le\alpha\}&\ \le\
\exp\Big(-\frac{\big(\rho\lambda-\sqrt{m-1}-\rho\sqrt{c_{\alpha}}\big)^{2}}{2(1+\rho^{2})}\Big).\label{eq:Sclosed-gauss}
\end{align}
\end{subequations}
\end{corollary}

\begin{proof}
\emph{The inclusion.} With $\sigma=1$, $\mathbf S=G_{\Omega}(\lambda_{1})/G_{\Omega}(\lambda_{2})$.
By~\eqref{eq:order}, $\lmax(\Omega)\le\lambda_{2}\le\lambda_{1}$, and $\lambda_{2}\ge0$ almost surely by
Lemma~\ref{lem:pos}. Lemma~\ref{lem:Gratio} with $r=\Omega$ and $\ell_{i}=\lambda_{i}$ gives
$\mathbf S\le e^{-g}\sum_{j\le N}g^{j}/j!$ with $g=\tfrac12(\lambda_{1}^{2}-\lambda_{2}^{2})$. As $N+1=d$,
the density $g_{N}$ of its proof is the density of half a $\chi^{2}_{2d}$ variable, so the right-hand side is
$\PP\{\chi^{2}_{2d}>2g\}$. It is strictly decreasing in $g$ and equals $\alpha$ at $2g=c_{\alpha}$. So
$\mathbf S>\alpha$ forces $\lambda_{1}^{2}-\lambda_{2}^{2}<c_{\alpha}$.

\emph{The bound.} Lemma~\ref{lem:pos} and Lemma~\ref{lem:T2sharp} at $t=\hat t$ give $0\le\rho\lambda_{2}\le\sqrt\kappa\,\hat\sigma$, as $\|\Omega\|_{F}\ge\|\Omega\|_{\rm op}$,
and~\eqref{eq:sigmahat-norm} at $\hat t$, where $\nabla Z(\hat t)=0$, gives $\kappa\hat\sigma^{2}=\|Y\|^{2}-\lambda_{1}^{2}$.
With the inclusion, on $\{\mathbf S>\alpha\}$,
\[\rho^{2}\lambda_{1}^{2}\ <\ \rho^{2}\lambda_{2}^{2}+\rho^{2}c_{\alpha}\ \le\ \|Y\|^{2}-\lambda_{1}^{2}+\rho^{2}c_{\alpha}.\]
Let $\xi$ and $\chi$ be as in the proof of Corollary~\ref{cor:closed}, so that $Z(e)=\lambda+\xi$,
$\|Y\|^{2}=Z(e)^{2}+\chi^{2}$, $\xi\sim N(0,1)$ is independent of $\chi$, and $\chi$ has the law of $\chi_{m-1}$.
Also $\lambda_{1}\ge Z(e)$ and $\lambda_{1}\ge Z(-e)=-Z(e)$ by Lemma~\ref{lem:pos}, so $Z(e)^{2}\le\lambda_{1}^{2}$.
Hence, on $\{\mathbf S>\alpha\}$, $(1+\rho^{2})Z(e)^{2}<Z(e)^{2}+\chi^{2}+\rho^{2}c_{\alpha}$, that is
$\rho^{2}(\lambda+\xi)^{2}<\chi^{2}+\rho^{2}c_{\alpha}$. This proves~\eqref{eq:Sclosed-nct}. On this event
$\rho(\lambda+\xi)\le\rho|\lambda+\xi|<\chi+\rho\sqrt{c_{\alpha}}$, that is $\chi-\rho\xi>\rho\lambda-\rho\sqrt{c_{\alpha}}$.
The variable $\chi-\rho\xi$ is the supremum of~$\langle u-\rho\psi_{e},W\rangle$ over the unit vectors $u$ orthogonal
to $\psi_{e}$. This is a centred Gaussian process with variance $1+\rho^{2}$, and the mean of its supremum is
$\E\chi\le\sqrt{m-1}$. The concentration inequality used for Corollary~\ref{cor:closed}
gives~\eqref{eq:Sclosed-gauss} when $\rho\lambda-\rho\sqrt{c_{\alpha}}\ge\sqrt{m-1}$.
\end{proof}

At $\alpha=0.05$, $c_{\alpha}=12.6$ for $d=3$ and $18.3$ for $d=5$. The bound~\eqref{eq:Sclosed-nct} is below
$0.05$ from about $\lambda=6.8$ at $(k,d)=(3,3)$ and from $\lambda=10.0$ at $(3,5)$, against $11.4$ and $11.0$
for~$\mathbf T$ in Corollary~\ref{cor:closed}. By~\eqref{eq:Sclosed-gauss},
$\liminf_{\lambda\to\infty}-\lambda^{-2}\log\big(1-\PP_{\lambda}\{\mathbf S\le\alpha\}\big)\ge\rho^{2}/(2(1+\rho^{2}))=k/(2(3k-2))$,
the value in the column $\mathbf S$ of Remark~\ref{rem:bracket}; that the exponent of $\mathbf S$ is not larger is
not proved.

\subsection{The deficiency exponent}

Corollary~\ref{cor:closed} bounds the exponential rate of the deficiency from below. A configuration of the
normalised Hessian and of the noise estimate at which the test fails bounds it from above.

\begin{theorem}[Deficiency exponent]\label{thm:exponent}
Let $\alpha\in(0,1)$ and let $A$ be the constant~\eqref{eq:Adef}. Then
\[\lim_{\lambda\to\infty}\ -\frac1{\lambda^{2}}\log\big(1-\PP_{\lambda}\{\mathbf T\le\alpha\}\big)
\ =\ E:=\frac{\rho^{2}}{2(A^{2}+\rho^{2})}.\]
\end{theorem}

\begin{proof}
For $\rho\lambda\ge A\sqrt{m-1}$, \eqref{eq:closed-chernoff} gives
\[-\frac1{\lambda^{2}}\log\big(1-\PP_{\lambda}\{\mathbf T\le\alpha\}\big)\ \ge\
\frac{\big(\rho-A\sqrt{m-1}/\lambda\big)^{2}}{2(A^{2}+\rho^{2})},\]
and the right-hand side tends to $E$ as $\lambda\to\infty$. So the $\liminf$ is at least $E$. The matching bound,
$\limsup\le E$, is proved in Appendix~\ref*{app:proofs} of the Supplementary Material, through one configuration of
the normalised Hessian and of the noise estimate at which the test fails, and an open cone around it.
\end{proof}

\begin{remark}[The constant $A$ and the exponent]\label{rem:bracket}
The decimal values of $A$ below are computed by a maximisation over the spectra of $r$. For $d=3$ the spectrum is $(a,b)$ with $a^{2}+b^{2}=T_{\max}^{2}$,
and the maximum is in one variable; at $(3,3)$ it is attained at $(a,b)=(0.995,0.103)\,T_{\max}$. For $d\ge4$ the
column $A$ gives the maximum over the spectra $(a,b,\dots,b)$, which is at most $A$. When the top eigenvalue is $a$,
the others are at most $\min\{a,(T_{\max}^{2}-a^{2})^{1/2}\}$, so Lemma~\ref{lem:mono} bounds $A$ from above by a
maximum in one variable. This upper bound exceeds the column by $1.3\cdot10^{-3}$ at $(3,4)$, $4.6\cdot10^{-4}$ at
$(3,5)$, $1.9\cdot10^{-6}$ at $(5,5)$ and $3.6\cdot10^{-8}$ at $(7,5)$, and the column $E$ keeps the digits that both
values give. A search over all spectra from random starting points returns the value of the column.
At~$\alpha=0.05$, with the lower bound $k/(2(3k-2))$ on the exponent of the known-variance test~$\mathbf S$ from
Corollary~\ref{cor:Sclosed}:

\begin{center}\small
\begin{tabular}{@{}lcccc@{}}
\toprule
$(k,d)$ & $\kappa$ & $A$ & $E$ & $\mathbf S$, lower bound\\
\midrule
$(3,3)$ & 7   & 2.312 & 0.0615 & 0.2143\\
$(3,4)$ & 16  & 1.543 & 0.120  & 0.2143\\
$(3,5)$ & 30  & 1.287 & 0.156  & 0.2143\\
$(5,3)$ & 18  & 1.440 & 0.1158 & 0.1923\\
$(5,5)$ & 121 & 1.064 & 0.1778 & 0.1923\\
$(7,5)$ & 325 & 1.023 & 0.1789 & 0.1842\\
\bottomrule
\end{tabular}
\end{center}

The slopes measured in Section~\ref{sec:num} lie between $10\%$ below and $7\%$ above $E$.
\end{remark}

\subsection{Benchmarks at high signal}

The benchmarks at high signal are the maximum test and the Bayes-optimal test. Under the uniform prior on $e$, which
is least favourable by $O(d)$-invariance, the Bayes-optimal test statistic is the mixture likelihood ratio
\[\bar L_{\lambda}(Y)=e^{-\lambda^{2}/2}\int_{\Sph^{d-1}}e^{\lambda Z(e')}\,\mathrm{d}e',\]
with $\mathrm{d}e'$ the uniform probability measure on the sphere ($|B|$ denotes the measure of $B\subset\Sph^{d-1}$),
a Laplace transform of the profile likelihood at inverse temperature $\lambda$. Since
$\tfrac1\lambda\log\int e^{\lambda Z}\mathrm{d}e'\to\sup_{M}Z$ as $\lambda\to\infty$, the maximum test is the
zero-temperature limit of the Bayes-optimal test. It is not available here: its threshold needs the noise level,
which the observation does not determine (Section~\ref{sub:rank}), and the null law of $\sup Z$, which has no closed
form. The following estimate shows that the two tests differ by a level shift of order $d\log\lambda/\lambda$.

\begin{proposition}[Zero-temperature estimate]\label{prop:laplace}
Let $\mathcal L=\tfrac1\lambda\log\int_{\Sph^{d-1}}e^{\lambda Z}\mathrm{d}e'$ and $K=\sup_{M}\|\nabla^{2}Z\|_{\rm op}$. Then, on the event
$\{\lambda K\ge4(d-1)/\pi^{2}\}$,
\[0\ \le\ \sup_{M}Z-\mathcal L\ \le\ \frac{d-1}{2\lambda}\Big(1+\log\frac{K\lambda}{d-1}\Big)
+\frac{(d-2)\log(\pi/2)+\log\big(\pi(d-1)\big)}{\lambda}.\]
\end{proposition}

\begin{proof}
The upper bound $\mathcal L\le\sup Z$ is immediate. For the lower one, let $\mathcal C_{r}$ be the geodesic cap of
radius $r\le\pi/2$ around $\hat t$. Since $\nabla Z(\hat t)=0$ and $\|\nabla^{2}Z\|_{\rm op}\le K$ on $M$, Taylor's
inequality along the geodesic gives $Z\ge\sup Z-\tfrac12Kr^{2}$ on $\mathcal C_{r}$, so
$\int e^{\lambda Z}\mathrm{d}e'\ge|\mathcal C_{r}|\,e^{\lambda(\sup Z-Kr^{2}/2)}$ and
\[\sup Z-\mathcal L\ \le\ \frac{Kr^{2}}{2}+\frac{\log\big(1/|\mathcal C_{r}|\big)}{\lambda}.\]
From $\sin\theta\ge2\theta/\pi$ on $[0,\pi/2]$ and $\int_{0}^{\pi}\sin^{d-2}\le\pi$ one gets
$|\mathcal C_{r}|\ge(2/\pi)^{d-2}r^{d-1}/(\pi(d-1))$, so the second term is at most
$\big[(d-1)\log(1/r)+(d-2)\log(\pi/2)+\log(\pi(d-1))\big]/\lambda$. Minimising over $r$ gives
$r^{2}=(d-1)/(K\lambda)$, at most $\pi^{2}/4$ on this event, and the stated bound.
\end{proof}

\begin{remark}[What exactness costs]\label{rem:high}
Take $\sigma=1$ known. An oracle that knew $e$ would test~$N(0,1)$ against $N(\lambda,1)$ in one dimension, because
$\|\psi_{e}\|=1$. Its deficiency is $\Phi(z_{\alpha}-\lambda)=e^{-\lambda^{2}(1+o(1))/2}$, with $z_{\alpha}$ the
$(1-\alpha)$-quantile of $N(0,1)$, and by the Neyman--Pearson lemma no test of level $\alpha$ has a smaller
deficiency at $e$. The maximum test has deficiency at most~$\PP\{Z(e)\le u\}=\Phi(u-\lambda)$ at its null threshold
$u$, so its exponent is $1/2$. The Bayes-optimal test maximises the power averaged over $e$ among the tests of
level $\alpha$. By invariance, its power and that of the maximum test do not depend on $e$. So its deficiency is
also at most $\Phi(u-\lambda)$, and its exponent is $1/2$ as well. For the $t$-spacing test the exponent is
$E=\rho^{2}/(2(A^{2}+\rho^{2}))$ (Theorem~\ref{thm:exponent}), below $1/2$.

With $\sigma$ known, the spacing test rejects on a shift of the squared level. The liminf of
$-\lambda^{-2}\log(1-\PP_{\lambda}\{\mathbf S\le\alpha\})$ is at least $\rho^{2}/(2(1+\rho^{2}))=k/(2(3k-2))$, that is
$k/(3k-2)$ of $1/2$ (Corollary~\ref{cor:Sclosed}); that the exponent of $\mathbf S$ is below $1/2$ is not proved.
This lower bound exceeds the exponent $E$ of $\mathbf T$, since $A>1$. So the $t$-test pays more than the
known-variance test on the exponential scale: the exponent of $\mathbf T$ is below the liminf for $\mathbf S$, at every
$(k,d)$ and every level.
\end{remark}

\subsection{Growing dimension}\label{sub:dim}

\begin{theorem}[Detection at the scale $\sqrt d$]\label{thm:sqrtd}
Let $\alpha\in(0,1)$ and $\tau,\delta>0$, and put
\begin{subequations}
\begin{equation}
\begin{gathered}
\Upsilon_{k}=\sup_{0<\theta<\pi}\frac{(|\cos\theta|+\sin\theta)^{k}-|\cos\theta|^{k}-k|\cos\theta|^{k-1}\sin\theta}{1-\cos^{k}\theta},\\
\bar{\mathcal X}=\sqrt{kd}+\tau,\qquad \bar v=(\Upsilon_{k}+k-1)\bar{\mathcal X}+(k-1)\tau,\\
\bar\ell=(k-1)(2\bar{\mathcal X}+\tau),\qquad \bar n^{2}=(\lambda+\tau)^{2}+m-1+\delta .
\end{gathered}
\end{equation}
If $\lambda\ge\bar v+\tau$, $m(\lambda-\tau)^{2}\ge d\,\bar n^{2}$ and
\begin{equation}\label{eq:sqrtd-H}
d\log\frac{\lambda-\tau}{\bar v}+N\log\frac{\bar v}{\bar v-\bar\ell}
+\frac m2\log\Big(1-\frac{(\lambda-\tau)^{2}-\bar v^{2}}{\bar n^{2}}\Big)\ \le\ \log\alpha,
\end{equation}
then
\begin{equation}
1-\PP_{\lambda}\{\mathbf T\le\alpha\}\ \le\ e^{-\tau^{2}/2}+2\big(1-\Phi(\tau)\big)+\PP\{\chi^{2}_{m-1}>m-1+\delta\}.
\end{equation}
In particular, let $v_{k}=(\Upsilon_{k}+k-1)\sqrt k$ and $\ell_{k}=2(k-1)\sqrt k$, and let $C_{k}>v_{k}$ solve
\begin{equation}\label{eq:sqrtd-C}
\tfrac12\big(C_{k}^{2}-v_{k}^{2}\big)-\log\frac{C_{k}}{v_{k}}\ =\ \log\frac{v_{k}}{v_{k}-\ell_{k}}.
\end{equation}
If $k$ is fixed and
$\lambda_{d}/\sqrt d\to\bar\lambda\in(C_{k},\infty)$ as $d\to\infty$, then $\PP_{\lambda_{d}}\{\mathbf T\le\alpha\}\to1$.
Moreover
\begin{equation}\label{eq:sqrtd-Cbound}
C_{k}\ \le\ 1.72\cdot2^{k/2}\sqrt k .
\end{equation}
\end{subequations}
\end{theorem}

In the proof, $\bar{\mathcal X}$ bounds $\sup X$, $\bar v$ bounds $\lambda_{2}$, $\bar\ell$ bounds $\lmax(\Omega)$ and
$\bar n$ bounds $\|Y\|$, outside an event whose probability is at most the bound of the theorem.
The proof is in Appendix~\ref*{app:proofs} of the Supplementary Material; the bound on $C_{k}$ is in
Lemma~\ref*{lem:Ck} there. Numerically \mbox{$\Upsilon_{3}=2.214$}, and
\mbox{$C_{3}=7.70$}, \mbox{$C_{5}=21.74$} and \mbox{$C_{7}=47.06$}.
As the order $k$ grows, the constant~$C_{k}$ grows like $2^{k/2}\sqrt k$: $C_{k}/(2^{k/2}\sqrt k)\to1$
as $k\to\infty$ along the odd integers (Lemma~\ref*{lem:Ck} of the Supplementary Material). Figure~\ref{fig:Ck} shows
this ratio for the odd orders~$k\le999$. The ratio is largest at $k=5$, where it is $1.72$. It decreases from
$k=5$ to~$k=999$, and it is below $1.02$ for~$k\ge21$.

\begin{figure}[ht]
\centering
\includegraphics[width=0.55\textwidth]{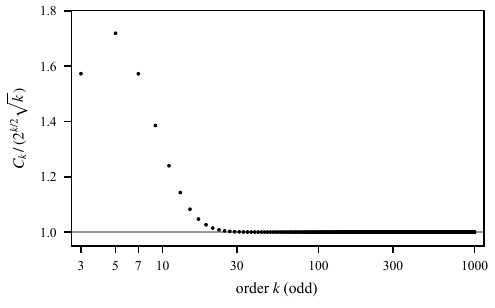}
\caption{The ratio $C_{k}/(2^{k/2}\sqrt k)$ for the odd orders $k=3,5,\dots,999$, where $C_{k}$ is defined
by~\eqref{eq:sqrtd-C} in Theorem~\ref{thm:sqrtd}. The axis of the orders is logarithmic. The ratio tends to one
(Lemma~\ref*{lem:Ck} of the Supplementary Material).}
\label{fig:Ck}
\end{figure}

\begin{remark}[Growing dimension]\label{rem:dim}
Let $k$ and $\alpha$ be fixed and $d\to\infty$. The constant $A$ tends to $1$ (Lemma~\ref*{lem:Alim} of the
Supplementary Material). So the exponent of Theorem~\ref{thm:exponent}, which is a limit as $\lambda\to\infty$ at
fixed $d$, tends to $k/(2(3k-2))$. This is the lower bound for $\mathbf S$ of Corollary~\ref{cor:Sclosed}. As $A>1$, the
bound~\eqref{eq:closed-nct} is at least $\PP\{\lambda+\xi<\chi_{m-1}/\rho\}$, which tends to one if~$\lambda/\sqrt m$ tends
to a constant below $1/\rho$, since $\chi_{m-1}/\sqrt m\to1$ in probability. So the bound is informative only when
$\lambda$ is at least of order $\sqrt m$, that is $d^{k/2}$. This is a property of the bound, not of the test: by
Theorem~\ref{thm:sqrtd} the power tends to one on the scale $\sqrt d$. The power of a test at the worst $e$ is at
most its power averaged over a uniform $e$. For $\lambda=\bar\lambda\sqrt d$ with a fixed $\bar\lambda$ below a
constant that depends on $k$, $1.12$ for $k=3$, the total variation distance between the null law and this average tends to zero
\citep[Theorem~3.2]{perry2020statistical}. So no test of level $\alpha$ has power larger than $\alpha+o(1)$ at every
$e$, and the power of $\mathbf T$, which does not depend on $e$, is at most $\alpha+o(1)$. We do not know the
constant of the threshold of the~$t$-spacing test.
\end{remark}

\section{Minimax estimation of the direction}\label{sec:main}

\subsection{Notation}

For $\lambda>0$ and $b\in(0,\lambda^{2}]$, the cap of size $b$ around $e$ is
\[\mathcal B_{b}=\big\{t\in M:\langle e,t\rangle\ge\sqrt{1-b/\lambda^{2}}\big\}.\]
A point $t$ is in $\mathcal B_{b}$ if and only if $\langle e,t\rangle\ge0$ and $\sin^{2}\angle(t,e)\le b/\lambda^{2}$.

Let $\delta^{\lambda}(u;B)$ be the expected number
of critical points of $Z$ in $B=\{t:\langle e,t\rangle\in I\}$, for an interval $I\subset[-1,1]$, with value at least $u$. By the Kac--Rice formula,
\[\delta^{\lambda}(u;B)=\omega_{d-2}\Big(\frac{k-1}{\pi}\Big)^{\frac N2}\!\int_{I}\!(1-c^{2})^{\frac{d-3}{2}}e^{-\frac{k\lambda^{2}}{2}c^{2k-2}(1-c^{2})}\!\int_{u}^{\infty}\!\!A_{\lambda}(\rho x;c)\varphi(x-\lambda c^{k})\,\mathrm{d}x\,\mathrm{d}c ,\]
where
\begin{align*}
\omega_{d-2}&=2\pi^{(d-1)/2}/\Gamma(\tfrac{d-1}{2})\ \text{is the area of $\Sph^{d-2}$;}\\
A_{\lambda}(z;c)&=\E\big|\det(z\,\mathrm{Id}-G_{\lambda})\big|\ \text{is the mean absolute determinant.}
\end{align*}

\subsection{Bounds on the risk}

The test returns the argmax $\hat t$. At fixed $(k,d)$ the classical theory gives, for the maximum likelihood
estimator, a limit law with covariance $(k\lambda^{2})^{-1}$ on the tangent space; the bounds below hold at
finite $\lambda$, and the upper ones contain the control of the far critical points of $Z$, a global statement
which the local theory does not give. Let the loss be $\ell(\hat e,e)=1-\langle\hat e,e\rangle^{2}=
\sin^{2}\angle(\hat e,e)$. For an estimator $\hat e=\hat e(Y)$ with values in $\Sph^{d-1}$, let
\begin{subequations}\label{eq:risk}
\begin{align}
R_{\lambda}(\hat e)&=\E\big[\sin^{2}\angle(\hat e,e)\big],\label{eq:risk-est}\\
R(\lambda)&=\inf_{\hat e}\,R_{\lambda}(\hat e),\label{eq:risk-bayes}
\end{align}
\end{subequations}
the expectation being under the joint law of $(e,Y)$ with $e$ uniform on $\Sph^{d-1}$, and the infimum over all
estimators. So $R_{\lambda}(\hat e)$ is the Bayes risk of $\hat e$ under the uniform prior, and $R(\lambda)$ is the
Bayes risk. The minimax risk is at least $R(\lambda)$, a supremum over $e$ being at least an average. The risk of
$\hat t$, which is $O(d)$-equivariant, does not depend on $e$, so it equals $R_{\lambda}(\hat t)$. A bracket on
$R(\lambda)$ whose upper side is attained by $\hat t$ is therefore a statement about the minimax risk.

Part~{\rm(iii)} of the theorem holds beyond a threshold $\lambda_{0}(k,d)$, which is defined by explicit
inequalities. Let $m_{N}=\sqrt{2N}+\frac56\sqrt{\max(2N,8)}$, $K_{N}=2\sqrt2\,(\frac54+m_{N}/N)^{N}$ and
$\omega_{d-1}=2\pi^{d/2}/\Gamma(\frac d2)$, the area of $\Sph^{d-1}$. Let $\lambda_{0}(k,d)$ be the least number such that
the five conditions
\begin{subequations}\label{eq:lam0}
\begin{gather}
\lambda^{2}\ge128\log\lambda,\label{eq:lam0a}\\
\lambda>\frac{4N}{3\rho},\qquad \lambda\ge\frac{16\sqrt{2N}}{9\rho},\label{eq:lam0b}\\
g_{1}:=\tfrac{87}{128}\rho\lambda-\sqrt{2N}\ge0,\qquad K_{N}\lambda^{2}e^{-g_{1}^{2}/4}\le1,\label{eq:lam0c}\\
g_{3}:=\tfrac{9}{16}\rho\lambda-\sqrt{2N}\ge0,\qquad
\Big(\frac54+K_{N}e^{-g_{3}^{2}/4}\Big)e^{\frac{128N(N+1)}{169\lambda^{2}}}\,\frac{4k-5}{3k-4}\le2,\label{eq:lam0d}\\
\lambda\ge4\sqrt{d+3},\qquad
\tfrac{10(4k-1)}{3}\,\omega_{d-1}\big(\tfrac{k-1}{\pi}\big)^{N/2}\big(\rho\lambda+\rho\sqrt N+m_{N}\big)^{N}
\lambda^{4}e^{-\lambda^{2}/32}\le1\label{eq:lam0e}
\end{gather}
\end{subequations}
hold at every $\lambda>\lambda_{0}(k,d)$. This number is finite, and each condition, once true at
some~\mbox{$\lambda\ge e$}, stays true for larger $\lambda$ (Appendix~\ref*{app:goe} of the Supplementary Material).
Conditions~\eqref{eq:lam0a} and~\eqref{eq:lam0b} make the choice of $\theta_{1}$ and $u$ of part~{\rm(iii)} admissible in
part~{\rm(ii)}, and~\eqref{eq:lam0c}, \eqref{eq:lam0d} and~\eqref{eq:lam0e} bound its three error terms.

\begin{theorem}[Estimation of the direction: lower and upper bounds]\label{thm:estim}
\begin{subequations}
\begin{itemize}
\item[\rm(i)] For every estimator $\hat e$,
\begin{equation}
R_{\lambda}(\hat e)\ \ge\ \frac{\big((d-1)-d\,R_{\lambda}(\hat e)\big)^{2}}{(d-1)\,k\,\lambda^{2}} .
\end{equation}
\item[\rm(ii)] For the argmax $\hat t$, let $\theta_{1}\in(0,\pi/2)$ and let $u$ satisfy
\[\begin{gathered}\frac{d-1}{\rho}<u<\lambda\cos^{k}\theta_{1},\\
u\ \ge\ \max\Big\{4(k-1)\lambda\sin^{2}\theta_{1},\ (k-1)\lambda\sin^{2}\theta_{1}+\frac{\sqrt{2(d-1)}}{\rho}\Big\}.\end{gathered}\]
Then
\begin{equation}
R_{\lambda}(\hat t)\ \le\
\frac{d-1}{k\lambda^{2}}\cdot\frac{1+\epsilon_{A}}{(1-\varrho)^{\frac{d+1}{2}}}\ +\ \mathrm{Rem}_{\rm far}\ +\ \Phi(u-\lambda),
\end{equation}
with $\varrho=1-\inf_{0<\theta\le\theta_{1}}\cos^{2k-2}\theta\,(\sin\theta/\theta)^{2}$, $\epsilon_{A}$ the
relative error of Lemma~\ref*{lem:Aupper} and $\mathrm{Rem}_{\rm far}=\delta^{\lambda}(u;\{\theta>\theta_{1}\})$, $\theta$ the angle between $t$ and $e$.
\item[\rm(iii)] For the argmax $\hat t$ and $\lambda>\lambda_{0}(k,d)$, the threshold defined by~\eqref{eq:lam0},
\begin{equation}
R_{\lambda}(\hat t)\ \le\ \frac{d-1}{k\lambda^{2}}\exp\Big(\frac{c_{1}\log\lambda+c_{1}'}{\lambda^{2}}\Big)
+2\,\PP\big\{\chi^{2}_{d-1}>8\log\lambda\big\}+\frac{1}{4\lambda^{4}},
\end{equation}
with $c_{1}=\frac{27(3k-2)(d+1)}{16k}+\frac{12(k-1)}{k}$ and $c_{1}'=\frac89\big((d-1)^{2}+2\big)$.
\end{itemize}
\end{subequations}
\end{theorem}

Since $R(\lambda)\le R_{\lambda}(\hat t)$, {\rm(i)} at a Bayes rule and {\rm(ii)} or {\rm(iii)} bracket $R(\lambda)$.

The three error terms of {\rm(ii)} have simple bounds, proved in Lemmas~\ref*{lem:simple} and~\ref*{lem:far} of the
Supplementary Material. For every $\theta_{1}\in(0,\pi/2)$, $\varrho\le(k-\frac23)\theta_{1}^{2}$; $\varrho$ is the
relative loss in the exponent of the Gaussian weight inside the cap $\theta\le\theta_{1}$. Assume now the hypotheses
of {\rm(ii)}. If also $u\ge4\sqrt{2(d-1)}/\rho$, then
\begin{subequations}
\begin{equation}
\epsilon_{A}\ \le\ \frac{(k-1)\lambda\sin^{2}\theta_{1}}{u}+\frac{3(d-1)^{2}}{2\lambda^{2}\cos^{2k}\theta_{1}}
+3e^{\frac54(d-1)-\frac{\rho^{2}u^{2}}{16}},
\end{equation}
the three terms coming from the spike of the Hessian, from the Gaussian moments of the value and from the negative
part of the determinant. If $3\rho u\ge4\sqrt{2(d-1)}$, then
\begin{equation}
\mathrm{Rem}_{\rm far}\ \le\ C_{M}\,\PP\big\{\chi^{2}_{d-1}>k\lambda^{2}\cos^{2k-2}\theta_{1}\sin^{2}\theta_{1}\big\}+F_{\rm far},
\end{equation}
\end{subequations}
the chi-square tail counting the critical points near $e$ beyond the angle $\theta_{1}$, and $F_{\rm far}$ those far
from $e$. The factor $C_{M}$ and the remainder $F_{\rm far}$ are explicit functions of $(k,d,\lambda,u)$.
For $\lambda>\lambda_{0}(k,d)$ and the choice of $\theta_{1}$ and $u$ below, $C_{M}\le2$ and $F_{\rm far}\le\lambda^{-4}/10$.

In {\rm(iii)} the additive term is $O(\lambda^{-4}(\log\lambda)^{(d-3)/2})$: for $\log\lambda\ge(d-3)/2$,
$2\,\PP\{\chi^{2}_{d-1}>8\log\lambda\}\le8(4\log\lambda)^{(d-3)/2}\lambda^{-4}/(3\,\Gamma(\frac{d-1}{2}))$
(Appendix~\ref*{app:goe}). Rounded up, $\lambda_{0}(k,d)=28.83$, $33.11$, $29.50$, $38.52$ and $34.52$ at
$(k,d)=(3,3)$, $(3,5)$, $(5,3)$, $(3,8)$ and $(7,5)$. Part {\rm(iii)} is {\rm(ii)} at
$\theta_{1}=3\lambda^{-1}\sqrt{\log\lambda/k}$ and $u=\lambda-2\sqrt{2\log\lambda}$, with the bounds above.

\begin{proof}
(i) Write $R=R_{\lambda}(\hat e)$. Let $v_{1},\dots,v_{d-1}$ be an orthonormal frame of $T_{e}\Sph^{d-1}$ and
\[S_{i}=\partial_{v_{i}}\log\big(\pi(e)f(Y\mid e)\big)\]
the score of the joint law of $(e,Y)$. The prior being uniform,
$\pi$ is constant and $S_{i}=\lambda\langle Y,\partial_{v_{i}}\psi_{e}\rangle$; since
$\langle\psi_{e},\partial_{v}\psi_{e}\rangle=\tfrac12\partial_{v}\|\psi_{e}\|^{2}=0$ and $Y=\lambda\psi_{e}+W$,
$S_{i}=\lambda\langle W,\partial_{v_{i}}\psi_{e}\rangle$. Differentiating
$\langle\psi_{e},\psi_{e'}\rangle=\langle e,e'\rangle^{k}$ twice gives
$\langle\partial_{v_{i}}\psi_{e},\partial_{v_{j}}\psi_{e}\rangle=k\delta_{ij}$, so $\E\|S\|^{2}=(d-1)k\lambda^{2}$.

Write $a=\langle\hat e,e\rangle$ and $\ell=1-a^{2}$, and take the tangent field
\begin{subequations}
\begin{equation}
g(e)=a\,\big(\hat e-a\,e\big),\qquad \|g\|^{2}=a^{2}\ell .
\end{equation}
The factor $a$ is what makes the argument insensitive to the sign of $a$, which the loss does not see and which no
estimator can control: the field $\hat e-a\,e$ alone would produce the factor~$(\E a)^{2}$, and $\E a$ vanishes
when the data are uninformative. Moving $e$ along a unit tangent~$v$, with~$\partial_{v}a=\langle\hat e,v\rangle$
and $\langle e,v\rangle=0$,
\begin{equation}
\partial_{v}g=\langle\hat e,v\rangle\hat e-2a\langle\hat e,v\rangle e-a^{2}v,\qquad
\langle\partial_{v}g,v\rangle=\langle\hat e,v\rangle^{2}-a^{2},
\end{equation}
and summing over the frame, with $\sum_{i}\langle\hat e,v_{i}\rangle^{2}=\|\hat e-a\,e\|^{2}=\ell$,
\begin{equation}
\operatorname{div}_{e}g=\ell-(d-1)a^{2}=\ell-(d-1)(1-\ell).
\end{equation}
\end{subequations}
For fixed $Y$ the field $e\mapsto g(e)$ is a polynomial in $e$, $\hat e=\hat e(Y)$ being a fixed unit vector,
so the divergence theorem on the closed sphere gives
$\int\langle g,\nabla_{e}(\pi f(Y\mid e))\rangle\,\mathrm{d}e=-\int\operatorname{div}_{e}g\,\pi f(Y\mid e)\,\mathrm{d}e$
for each $Y$; integrating in $Y$ (Fubini, the integrand being integrable) and using
$S\,\pi f=\nabla_{e}(\pi f)$, $\E\langle g,S\rangle=-\E[\operatorname{div}_{e}g]=(d-1)-d\,R$. No regularity of
$Y\mapsto\hat e(Y)$ beyond measurability is used. Cauchy--Schwarz with
$\E\|g\|^{2}=\E[a^{2}\ell]\le R$ gives $\big((d-1)-dR\big)^{2}\le R\,(d-1)k\lambda^{2}$.

(ii) The hypotheses require $\lambda>(d-1)/(\rho\cos^{k}\theta_{1})$; the hypotheses other than $u>(d-1)/\rho$ are used only inside
Lemma~\ref*{lem:Aupper}. On $\{\sin^{2}\angle(\hat t,e)>b/\lambda^{2}\}$ the argmax lies outside $\mathcal B_{b}$.
Either $\sup Z\le u$, of probability at most $\PP\{Z(e)\le u\}=\Phi(u-\lambda)$ since $Z(e)\sim N(\lambda,1)$, or
$\hat t$ is a critical point outside $\mathcal B_{b}$ with value above $u$, of probability at most
$\delta^{\lambda}(u;\mathcal B^{c}_{b})$ by Markov. Lemma~\ref*{lem:Aupper} bounds that count by
\begin{subequations}
\begin{equation}
\frac{1+\epsilon_{A}}{(1-\varrho)^{\frac{d-1}{2}}}\ \PP\big\{\chi^{2}_{d-1}>kb(1-\varrho)\big\}+\mathrm{Rem}_{\rm far}.
\end{equation}
Integrating the tail of the loss,
\begin{equation}
\begin{aligned}R_{\lambda}(\hat t)&=\frac{1}{\lambda^{2}}\int_{0}^{\lambda^{2}}\PP\Big\{\sin^{2}\angle>\frac{b}{\lambda^{2}}\Big\}\,\mathrm{d}b\\
&\le\frac{1}{\lambda^{2}}\int_{0}^{\infty}\frac{(1+\epsilon_{A})\,\PP\{\chi^{2}_{d-1}>kb(1-\varrho)\}}{(1-\varrho)^{\frac{d-1}{2}}}\,\mathrm{d}b
+\mathrm{Rem}_{\rm far}+\Phi(u-\lambda),\end{aligned}
\end{equation}
\end{subequations}
and $\int_{0}^{\infty}\PP\{\chi^{2}_{d-1}>kb(1-\varrho)\}\,\mathrm{d}b=\E\chi^{2}_{d-1}/(k(1-\varrho))=(d-1)/(k(1-\varrho))$.

Part {\rm(iii)} is proved in Appendix~\ref*{app:goe} of the Supplementary Material.
\end{proof}

\begin{corollary}[Asymptotic minimaxity with the sharp constant]\label{cor:minimax}
As $\lambda\to\infty$,
\[\begin{aligned}\frac{d-1}{k\lambda^{2}}\Big(1-\frac{d\,R(\lambda)}{d-1}\Big)^{2}\ &\le\ R(\lambda)\ \le\ R_{\lambda}(\hat t)\\
&\le\ \frac{d-1}{k\lambda^{2}}\Big(1+O\Big(\frac{d^{2}+d\log\lambda}{\lambda^{2}}\Big)\Big)+O\big(\lambda^{-4}(\log\lambda)^{\frac{d-3}{2}}\big),\end{aligned}\]
so that $\lambda^{2}R_{\lambda}(\hat t)\to(d-1)/k$ and $\lambda^{2}R(\lambda)\to(d-1)/k$: the argmax is
asymptotically minimax for the direction, with the sharp constant.
\end{corollary}

\begin{proof}
The upper bound is Theorem~\ref{thm:estim}{\rm(iii)}, with $e^{x}-1\le xe^{x}$, $c_{1}\le\frac{81}{16}(d+1)+12$,
$c_{1}'\le\frac89(d^{2}+2)$ and the bound on the chi-square tail given after the theorem. It forces $R(\lambda)\to0$,
and the lower bound is Theorem~\ref{thm:estim}{\rm(i)} applied to the Bayes rule, a top eigenvector of
$\E[ee^{\top}\mid Y]$, whose risk is $R(\lambda)$.
\end{proof}

\begin{remark}
The localisation radius is never chosen: it is integrated out, and the constant~$(d-1)/k$ is the whole
chi-square tail. This is why the two sides of Theorem~\ref{thm:estim} carry the same constant. Solving the
inequality of part {\rm(i)} gives, at $(k,d)=(3,5)$, the values $0.0472$, $0.01291$, $0.003306$ and $0.000832$ for
$\lambda=5,10,20,40$, which are $88.54\%$, $96.80\%$, $99.18\%$ and $99.79\%$ of $(d-1)/(k\lambda^{2})$.
\end{remark}
\section{Scope: what happens at vanishing signal}\label{sec:low}

At high signal the maximum test is the zero-temperature limit of the Bayes-optimal test
(Proposition~\ref{prop:laplace}). At vanishing signal this picture reverses. The
mixture likelihood ratio has no first-order term, because $k$ is odd and~$\int e'^{\otimes k}\,\mathrm{d}e'=0$ by
antipodal symmetry, so
\[\bar L_{\lambda}=1+\frac{\lambda^{2}}{2}\,(T_{\rm q}-1)+O(\lambda^{3}),\qquad
T_{\rm q}=\int_{\Sph^{d-1}}Z(e')^{2}\,\mathrm{d}e'=\langle Y,\Sigma Y\rangle,\quad \Sigma=\int\psi_{e'}\psi_{e'}^{\top}\,\mathrm{d}e',\]
a quadratic form in the observed tensor. Expanding $\E_{\lambda}\phi=\E_{0}[\phi\bar L_{\lambda}]$ to second order, every test $\phi$ of size $\alpha$ satisfies
$\E_{\lambda}[\phi]=\alpha+\tfrac{\lambda^{2}}{2}\Cov_{0}(\phi,T_{\rm q})+o(\lambda^{2})$, a linear functional of $\phi$
under the constraint $\E_{0}\phi=\alpha$, hence maximised by the upper-tail test on $T_{\rm q}$. \emph{At vanishing signal
the optimal procedure is the quadratic test, not this one.}

\medskip

The null law of $T_{\rm q}$ is explicit: by Funk--Hecke the kernel $\langle e,e'\rangle^{k}$ has eigenvalue
\[\mu_{j}=\int_{-1}^{1}t^{k}\,\frac{C_{j}^{((d-2)/2)}(t)}{C_{j}^{((d-2)/2)}(1)}\;
\frac{\Gamma(\frac d2)}{\sqrt\pi\,\Gamma(\frac{d-1}{2})}(1-t^{2})^{\frac{d-3}{2}}\,\mathrm{d}t,\]
with $C_{j}^{((d-2)/2)}$ the Gegenbauer polynomial, on the spherical harmonics of degree $j\le k$ with $j\equiv k$ modulo $2$, of dimension
$n_{j}=\binom{d+j-1}{j}-\binom{d+j-3}{j-2}$, and all other eigenvalues vanish,
so $T_{\rm q}\overset{d}{=}\sum_{j}\mu_{j}\chi^{2}_{n_{j}}$ with independent factors. One checks
$\sum_{j}\mu_{j}n_{j}=\E X^{2}=1$ and $\sum_{j}\mu_{j}^{2}n_{j}=\E\,T_{\rm cor}^{2k}$, where
$T_{\rm cor}=\langle e,e'\rangle$ for independent uniform directions, so the mean shift under the alternative is
$\lambda^{2}\E T_{\rm cor}^{2k}$ and the null variance is $2\E T_{\rm cor}^{2k}$, with
$\E T_{\rm cor}^{2k}=(2k-1)!!/\prod_{i<k}(d+2i)$.

\medskip

How much does a test built on the maximum lose? With $q_{\alpha}$ and $t_{\alpha}$ the $(1-\alpha)$-quantiles of
$\sup X$ and of $T_{\rm q}$ under the null, its local coefficient is
\begin{subequations}
\begin{equation}
\tfrac12\Cov_{0}(\ind\{\sup X>q_{\alpha}\},T_{\rm q}),
\end{equation}
and the optimal one is
\begin{equation}
\tfrac12\Cov_{0}(\ind\{T_{\rm q}>t_{\alpha}\},T_{\rm q})=\tfrac12(\E[T_{\rm q};T_{\rm q}>t_{\alpha}]-\alpha);
\end{equation}
\end{subequations}
since $T_{\rm q}$ is a chi-square mixture and not Gaussian, the optimal coefficient has to be computed from its own upper tail. At $(k,d)=(3,3)$, where
$T_{\rm q}=\tfrac15\chi^{2}_{3}+\tfrac2{35}\chi^{2}_{7}$, simulation gives a local efficiency of $0.810$ at $\alpha=0.05$
and $0.777$ at $\alpha=0.01$.

\medskip

When $\sigma$ is unknown the quadratic test is not available. Among tests invariant under $Y\mapsto aY$, $a>0$,
and under the orthogonal group, the locally most powerful one rejects for large
$\tilde T=\sum_{j}\mu_{j}Q_{j}/\sum_{j}Q_{j}$, $Q_{j}=\|P_{j}Y\|^{2}$ the squared norm of the projection on the
harmonics of degree $j$: the likelihood ratio of $Y/\|Y\|$, averaged over the orbit of $e$, is
$1+\frac{\lambda^{2}}{2}(m\tilde T-1)+O(\lambda^{3})$. At $(k,d)=(3,3)$, $\tilde T$ is increasing in
$Q_{1}/Q_{3}$ and $\tfrac73Q_{1}/Q_{3}$ is $F_{3,7}$ under the null; the local coefficient of this $F$-test is
$1.54\cdot10^{-2}$ at $\alpha=0.05$, $42\%$ of the quadratic test's and $52\%$ of the maximum
test's, and its power is $0.100$ at $\lambda=2$ and $0.191$ at $\lambda=5$.

\medskip

The window in which this matters is narrow. Comparing the two powers at fixed $(k,d)$ rather than their local
coefficients is a finite computation: under the alternative the quadratic statistic is the noncentral mixture
\[T_{\rm q}\ \overset{d}{=}\ \sum_{j}\mu_{j}\,\chi^{2}_{n_{j}}\big(\lambda^{2}\mu_{j}n_{j}\big),\]
because the projection of $\lambda\psi_{e}$ on the harmonics of degree $j$ has squared norm
$\lambda^{2}\mu_{j}n_{j}$, while the power of the maximum test is obtained by simulating the field. At
$(k,d)=(3,3)$:

\begin{center}\small
\begin{tabular}{@{}lcccccc@{}}
\toprule
& $\lambda=1.4$ & $1.6$ & $1.8$ & $2.0$ & $3.0$ & $5.0$\\
\midrule
quadratic test, $\alpha=0.05$ & 0.136 & 0.166 & 0.203 & 0.246 & 0.534 & 0.963\\
maximum test, $\alpha=0.05$ & 0.131 & 0.165 & 0.205 & 0.254 & 0.584 & 0.982\\
$\mathbf S$, $\alpha=0.05$ & 0.112 & 0.138 & 0.173 & 0.213 & 0.510 & 0.962\\
$\mathbf T$, $\alpha=0.05$ & 0.069 & 0.078 & 0.090 & 0.105 & 0.233 & 0.650\\
\bottomrule
\end{tabular}
\end{center}

The first two cross at $\lambda\approx1.7$ for $\alpha=0.05$ and at $\lambda\approx1.46$ for $\alpha=0.01$,
where the common power is $0.19$ and $0.045$ respectively. \emph{The quadratic test is preferable only where no
test has useful power}, and the maximum test leads everywhere else. The spacing tests, built on the same
maximum, pay in addition the price of exactness of Remark~\ref{rem:high}: $\mathbf S$ is below the maximum
test at every $\lambda$ and within $10^{-3}$ of the quadratic test at $\lambda=5$, and $\mathbf T$, which also
estimates $\sigma$ with $\kappa=7$ degrees of freedom, has power $0.65$ at $\lambda=5$. The rows $\mathbf S$
and~$\mathbf T$ come from $5\cdot10^{5}$ replicates each of the Karhunen--Lo\`eve simulation of the field used
for the measured exponent of Section~\ref{sec:num}, with standard error below $10^{-3}$; both statistics pass a
Kolmogorov--Smirnov test of uniformity on $10^{6}$ null replicates, $p=0.28$ and $0.99$.

\medskip

No single procedure is therefore optimal at both ends, and the end at which the maximum test is not optimal is the
end at which nothing works. What the $t$-spacing test has and the quadratic test has not is an exact level at every
$\lambda$, no need for the noise level, and the maximiser as a by-product, which Theorem~\ref{thm:estim} and
Corollary~\ref{cor:minimax} show to be an asymptotically minimax estimator of the direction.

\section{Numerical illustration}\label{sec:num}

The level maps and the power bound are computed from the closed forms of Appendix~\ref*{app:closed} of the
Supplementary Material, which were checked against numerical quadrature to ten significant digits. The estimation
sandwich also uses a Kac--Rice integral under the alternative, evaluated numerically. The simulations are described
in each paragraph.

\paragraph{The model constants} The Karhunen--Lo\`eve order, the degrees of freedom of $\hat\sigma^{2}$, and the
degrees of freedom of the reference Student density.

\begin{center}\small
\begin{tabular}{@{}lcccc@{}}
\toprule
$(k,d)$ & $(3,3)$ & $(3,5)$ & $(3,8)$ & $(5,6)$\\
\midrule
$m=\binom{d+k-1}{k}$ & 10 & 35 & 120 & 252\\
$\kappa=m-d$ & 7 & 30 & 112 & 246\\
$\nu=m-1$ & 9 & 34 & 119 & 251\\
\bottomrule
\end{tabular}
\end{center}

\paragraph{The price of knowing $\sigma$} The level map $u_{\alpha}(r,v)$ of~\eqref{eq:ualpha} gives the threshold of the
$t$-test on~$T_{1}$ at $T_{2}=v$ and $\Omega/\hat\sigma=r$. The same definition with $G_{r}$ gives the threshold of
$\mathbf S$.
Compare these two level maps at the same $\alpha$ and the same second-maximum threshold $v$, for one draw of
$\Omega=G/\rho$, $G\sim\GOE(N)$, at each $(k,d)$. The draw satisfies $v>\lmax(\Omega)$, and we take
$\hat\sigma=\sigma=1$, so that both maps are evaluated at $r=\Omega$. The $t$-test asks for a higher
first maximum, and that gap is the price of not having $\sigma$.

\begin{center}\small
\begin{tabular}{@{}lccccc@{}}
\toprule
$(k,d)$ & $\alpha$ & $v$ & $u_{\alpha}$, known variance & $u_{\alpha}$, $t$-test & extra level\\
\midrule
$(3,5)$ & 0.05 & 4 & 4.897 & 5.486 & $+0.589$\\
$(3,5)$ & 0.01 & 4 & 5.276 & 6.171 & $+0.896$\\
$(3,8)$ & 0.05 & 4 & 5.036 & 5.202 & $+0.166$\\
$(5,6)$ & 0.05 & 4 & 4.862 & 4.923 & $+0.061$\\
\bottomrule
\end{tabular}
\end{center}

The draw at $(3,5)$ is on the high side: over $3000$ draws the median extra level is $0.505$ and this one
exceeds $90\%$ of them; at $(3,8)$ and $(5,6)$ the draws are close to the median. The gap is not the inflation of
a single quantile: with the scale $\gamma$ of $H_{r}$, the $0.95$-quantile of the reference Student law is $3.4\%$
below the Gaussian one at $(3,5)$, and its tail is heavier only beyond $u=2.7$. The level map compares two
tail integrals of a polynomial of degree~$N$ against the two densities, so the heavier Student tail is penalised
$N$ times over. It falls below~$0.1$ in level as soon as $m$ is a few hundred.

\paragraph{The power bound} At $(k,d)=(3,5)$ and $\alpha=0.05$, the column of Remark~\ref{rem:bracket} and its upper
bound give $1.2873\le A\le1.2878$, and $m-1=34$. The bound~\eqref{eq:closed-nct} of Corollary~\ref{cor:closed} is $0.05$ at
$\lambda=11.0$ and $0.01$ at $\lambda=12.0$. Its Gaussian form~\eqref{eq:closed-chernoff} reaches these values at
$\lambda=13.1$ and $14.1$. In the table of the next paragraph, the deficiency is $0.05$ near $\lambda=6.2$ and $0.01$
near $\lambda=7.2$, by log-linear interpolation. The bound therefore needs a signal $1.7$ to $1.8$ times larger.

\paragraph{The true exponent, measured} The deficiency $1-\PP_{\lambda}\{\mathbf T\le\alpha\}$ can be simulated
directly, the field being generated from its Karhunen--Lo\`eve expansion and the two maxima, the Hessian and
$\hat\sigma$ computed from the realised path. The maximum $\lambda_{1}$ is found by gradient ascent on the sphere
from $41$ starting points, $40$ random ones and $e$. The supremum~\eqref{eq:lam2} is found in the same way, from
$40$ random points, the directions of the eigenvectors of $\Omega$ at $\hat t$ and the other local maxima of $Z$,
and $\lambda_{2}$ is the larger of this value and its limit $\lmax(\Omega)$ at $\hat t$. On~$2\%$ of the
replicates, run again with four times as many starting points, $\lambda_{1}$ never changed, $\lambda_{2}$ rose
$3$ times in $279\,000$, and no decision changed. Before measuring anything one checks the calibration:
with $10^{6}$ null replicates at each pair the statistic passes a Kolmogorov--Smirnov test of uniformity, with
$p=0.18$, $0.31$, $0.56$ and $0.88$ at $(3,3)$, $(3,4)$, $(3,5)$ and $(5,3)$, and the order
relation~\eqref{eq:order} holds in every replicate under the alternative. Under the alternatives the number of replicates
grows with $\lambda$, from $2\cdot10^{4}$ to $2\cdot10^{6}$ per point and $1.0\cdot10^{7}$ in all, so that at each
point the test fails to reject at least $78$ times. At $\alpha=0.05$:

\begin{center}\footnotesize\setlength{\tabcolsep}{3pt}
\begin{tabular}{@{}lccccccc@{}}
\toprule
$\lambda$ & 5 & 6 & 7 & 8 & 9 & 10 & 12\\
\midrule
$(3,3)$, $\kappa=7$  &  & $0.183$ & $0.0791$ & $2.90\cdot10^{-2}$ & $1.00\cdot10^{-2}$ & $2.84\cdot10^{-3}$ & $1.57\cdot10^{-4}$\\
$(3,4)$, $\kappa=16$ & $0.251$ & $8.65\cdot10^{-2}$ & $2.21\cdot10^{-2}$ & $4.01\cdot10^{-3}$ & $5.75\cdot10^{-4}$ & $6.10\cdot10^{-5}$ & \\
$(3,5)$, $\kappa=30$ & $0.238$ & $7.40\cdot10^{-2}$ & $1.40\cdot10^{-2}$ & $1.88\cdot10^{-3}$ & $1.20\cdot10^{-4}$ &  & \\
$(5,3)$, $\kappa=18$ & $0.175$ & $5.05\cdot10^{-2}$ & $1.02\cdot10^{-2}$ & $1.53\cdot10^{-3}$ & $1.95\cdot10^{-4}$ &  & \\
\bottomrule
\end{tabular}
\end{center}

Regressing $-\log$ of these on $\lambda^{2}$ gives a slope $0.065\pm0.002$ at $(3,3)$, flat to within $4\%$
over $\lambda\in[6,12]$, and $0.115\pm0.012$ at $(3,4)$, where
the slope is still climbing and the quoted uncertainty is dominated by that drift. At $(5,3)$ the slope is
$0.124\pm0.003$, flat over $\lambda\in[6,9]$. At $(3,5)$ the slope is not stable over the window. Weighted fits
give $0.136$, $0.143$ and~$0.153$ over $\lambda\in[6,9.5]$, $[7,9.5]$ and $[8,9.5]$, and we quote $0.14\pm0.01$.
The local slopes over successive intervals are $0.106$, $0.128$, $0.134$, $0.162$ and, on $[9,9.5]$,
$0.12\pm0.02$; the deficiency is $3.9\cdot10^{-5}$ at $\lambda=9.5$, from $78$ failures in $2\cdot10^{6}$
replicates. Two readings follow. The
exponent depends on the dimension, as $E$ does through $\kappa$ and $N$: a
$d$-free exponent cannot describe the three rows with $k=3$. And the measured values are close to the exponent of
Theorem~\ref{thm:exponent}:

\begin{center}\small
\begin{tabular}{@{}lcc@{}}
\toprule
$(k,d)$ & measured & Theorem~\ref{thm:exponent}, $E$\\
\midrule
$(3,3)$ & $0.065\pm0.002$ & 0.0615\\
$(3,4)$ & $0.115\pm0.012$ & 0.120\\
$(3,5)$ & $0.14\pm0.01$ & 0.156\\
$(5,3)$ & $0.124\pm0.003$ & 0.1158\\
\bottomrule
\end{tabular}
\end{center}

A slope over a finite window is not the limit, and the theorem gives the limit only: the measured slope is $6\%$ above
$E$ at $(3,3)$, $7\%$ above at $(5,3)$, where the difference is $2.7$ times the quoted uncertainty, within its error at $(3,4)$,
and $10\%$ below at $(3,5)$. At every $\lambda$ of the table above, $-\lambda^{-2}\log$ of the deficiency is
below $E$: the deficiency exceeds $e^{-E\lambda^{2}}$, by a factor~$1.7$ at $\lambda=6$ and $1.1$ at $\lambda=12$ at
$(3,3)$.

\paragraph{The estimation sandwich} Parts {\rm(i)} and {\rm(ii)} of Theorem~\ref{thm:estim} bracketed numerically at
$(k,d)=(3,5)$, the middle column being the common target $(d-1)/(k\lambda^{2})$.

\begin{center}\small
\begin{tabular}{@{}lcccccc@{}}
\toprule
$\lambda$ & $u$ & lower bound, {\rm(i)} & $(d-1)/(k\lambda^{2})$ & upper bound, Kac--Rice & ratio\\
\midrule
20 & 16 & $3.306\cdot10^{-3}$ & $3.333\cdot10^{-3}$ & $3.390\cdot10^{-3}$ & 1.017\\
20 & 18 & $3.306\cdot10^{-3}$ & $3.333\cdot10^{-3}$ & $2.604\cdot10^{-2}$ & 7.81\\
40 & 34 & $8.316\cdot10^{-4}$ & $8.333\cdot10^{-4}$ & $8.349\cdot10^{-4}$ & 1.0019\\
40 & 36 & $8.316\cdot10^{-4}$ & $8.333\cdot10^{-4}$ & $8.666\cdot10^{-4}$ & 1.0399\\
\bottomrule
\end{tabular}
\end{center}

The upper column is the bound $\int\sin^{2}\theta\,\delta^{\lambda}(u;\mathrm{d}\theta)+\Phi(u-\lambda)$ reached in the
proof of {\rm(ii)} before Lemma~\ref*{lem:Aupper} is applied, the Kac--Rice integral being evaluated numerically.
The right-hand side of {\rm(ii)}, with $\epsilon_{A}$ bounded as in Lemma~\ref*{lem:Aupper}, $\mathrm{Rem}_{\rm far}$
evaluated numerically, and optimised over $\theta_{1}$, is looser: $4.261\cdot10^{-3}$,
$2.695\cdot10^{-2}$, $8.934\cdot10^{-4}$ and $9.244\cdot10^{-4}$ in the four rows, ratios $1.278$, $8.09$,
$1.072$ and $1.109$. Optimised over $(\theta_{1},u)$, the explicit bound of {\rm(ii)} is within a factor $2$ of the
target from $\lambda\approx11$, $13$, $13$ and $16$ at $(3,3)$, $(5,3)$, $(3,5)$ and~$(3,8)$, and within $10\%$ from
$\lambda\approx26$, $28$, $34$ and $45$ at the same pairs. At $\lambda=10$ it is $4$ to~$46$ times the target: the
remainder $\varepsilon_{\det}$ of Lemma~\ref*{lem:Aupper}, the gap term $\Phi(u-\lambda)$ and~$\mathrm{Rem}_{\rm far}$
are then all of the order of the target. The Kac--Rice bound from which it is derived is within~$10\%$ of the lower
bound from $\lambda\approx9$ ($10.4$ at $(3,8)$). With the Kac--Rice column the two sides bracket the Bayes risk
within $2.6\%$ at $\lambda=20$ and within $0.40\%$ at $\lambda=40$. The rows
with the larger $u$ show the role of the gap $\lambda-u$: the term $\Phi(u-\lambda)$ enters the risk additively,
so it must be pushed below the target, which asks for $\lambda-u\gtrsim\sqrt{2\log\lambda^{2}}$;
at~$\lambda=40$, $\lambda-u=4$ costs $3.8\%$ and $\lambda-u=6$ costs nothing.

\paragraph{The figure of the introduction} Figure~\ref{fig:intro} is computed on the Karhunen--Lo\`eve
simulation above at $(k,d)=(3,8)$, $\kappa=112$, with $10^{5}$ replicates per point and the same replicates
for every test. Besides $\mathbf T$, $\mathbf S$ and the maximum test, it shows four competitors. The
studentised maximum $\sup Z/\hat\sigma$ is scale-free; its null law has no closed form and is simulated once
per~$(k,d)$, here on $10^{6}$ replicates. Compared instead with the quantile of $\sup Z/\sigma$, as if $\hat\sigma$
were exact, it rejects with probability $0.097$ at nominal level $0.05$ and $0.026$ at $0.01$. The locally
optimal invariant test of Section~\ref{sec:low} rejects for large $Q_{1}/Q_{3}$, with
$(n_{3}/n_{1})Q_{1}/Q_{3}\sim F_{8,112}$ under the null; its power is $0.20$ at $\lambda=4$ and $0.94$ at
$\lambda=12$. The unfolding test rejects for large~$s_{1}^{2}/\|Y\|^{2}$, $s_{1}$ the largest singular
value of the $d\times d^{k-1}$ unfolding of $Y$ \citep{montanari2014statistical}, calibrated by simulation; its
power is above that of $\mathbf T$ by at most $0.011$. The most powerful test invariant under rotations and
scalings, at a given $\lambda/\sigma$, rejects for large
\[\tilde L_{\lambda}(Y)=e^{-\lambda^{2}/2}\int_{\Sph^{d-1}}\frac{I(\lambda Z(e')/\|Y\|)}{I(0)}\,\mathrm{d}e',\qquad
I(c)=\int_{0}^{\infty}r^{m-1}e^{-r^{2}/2+cr}\,\mathrm{d}r,\]
the likelihood ratio of the maximal invariant \citep{wijsman1967cross}; the group being amenable, this test is
maximin over $(e,\sigma)$ at fixed $\lambda/\sigma$ among all tests of level $\alpha$, by the Hunt--Stein
theorem \citep[Section~8.5]{lehmann2005testing}. The integral over the sphere is estimated by importance
sampling with $2000$ points per replicate and the threshold is calibrated on the same statistic, so that the
level is exact and the power is a lower bound for that of the ideal test; on $3000$ replicates a cloud four
times larger changes the power by at most $0.002$. Its power exceeds that of the studentised maximum by at
most $0.011$ at every $\lambda$, the finite-$\lambda$ form of Proposition~\ref{prop:laplace}: the maximum
carries almost all the information of the likelihood ratio. The~$t$-spacing test is below both by at most
$0.08$, at $\lambda=5$, and the tests with $\sigma$ known are above it by at most $0.16$. Panel (c) shows $(3,8)$ and
$(5,3)$; the same simulation was also run at $(3,3)$ and $(3,5)$. The risk does not depend on $e$, by invariance, so it is
the Bayes risk of Theorem~\ref{thm:estim}. At $\lambda=10$ the risk of the argmax is within $10\%$ of the lower
bound at the four pairs and the unfolding estimator is above it by $7$ to $24\%$; at $\lambda=20$ and $(3,8)$ the
two figures are $2.5\%$ and $5.5\%$.

\paragraph{Local efficiency at vanishing signal} At $(k,d)=(3,3)$, the local coefficients of
Section~\ref{sec:low}, the field being simulated on $40\,000$ points of $\Sph^{2}$ with $2\cdot10^{5}$ replicates
and the optimal coefficient cross-checked against the exact mixture $\tfrac15\chi^{2}_{3}+\tfrac2{35}\chi^{2}_{7}$.

\begin{center}\small
\begin{tabular}{@{}lccccc@{}}
\toprule
$\alpha$ & $q_{\alpha}$ & $t_{\alpha}$ & $C_{\max}$ & $C_{\rm opt}$ & efficiency\\
\midrule
0.05 & 3.181 & 2.022 & $2.959\cdot10^{-2}$ & $3.655\cdot10^{-2}$ & 0.810\\
0.01 & 3.678 & 2.732 & $8.398\cdot10^{-3}$ & $1.081\cdot10^{-2}$ & 0.777\\
\bottomrule
\end{tabular}
\end{center}

\section*{Acknowledgments}
Yohann De Castro acknowledges the support of the Institut Universitaire de France (IUF).

\section*{Statements and Declarations}
\paragraph{Competing interests} The authors have no competing interests to declare that are relevant to the content of this article.
\paragraph{Data availability} The scripts that produce every number and table of this article are available at \url{https://github.com/ydecastro/TENSOR-KSS}, folder \texttt{paper2}. No other data were generated.

\section*{Disclosure of AI use}
In line with the Leiden Declaration on Artificial Intelligence and Mathematics
\citeyearpar{leiden2026declaration}, we disclose the use of AI agents (Claude Code, Anthropic) in the
preparation of this paper. These agents were used to audit the statements, the proofs and the
numerical tables, through repeated adversarial verification passes recorded in a change-log, and
to help with editorial tasks: consistency of notation, bibliography formatting, and typesetting.
The authors produced and checked every result and every proof themselves and take sole
responsibility for the correctness of the paper and for its citations.

\appendix
\makeatletter
\@addtoreset{theorem}{section}\@addtoreset{proposition}{section}\@addtoreset{lemma}{section}
\@addtoreset{remark}{section}
\makeatother
\renewcommand{\thetheorem}{\thesection.\arabic{theorem}}
\renewcommand{\theproposition}{\thesection.\arabic{proposition}}
\renewcommand{\thelemma}{\thesection.\arabic{lemma}}
\renewcommand{\theremark}{\thesection.\arabic{remark}}
\renewcommand{\theHtheorem}{\thesection.\arabic{theorem}}
\renewcommand{\theHproposition}{\thesection.\arabic{proposition}}
\renewcommand{\theHlemma}{\thesection.\arabic{lemma}}
\renewcommand{\theHremark}{\thesection.\arabic{remark}}
\numberwithin{equation}{section}

\section*{Supplementary Material}
The Supplementary Material contains the Kac--Rice formula with marks
(Appendix~\ref{app:kr}), the layer decomposition of the helix process and its drift, restated
from \citet{kss}, and the continuity of the second maximum at a point (Appendix~\ref{app:layers}), the lemmas on the Gaussian orthogonal ensemble used in the
proofs and the proof of Theorem~\ref*{thm:estim}{\rm(iii)} of the main paper (Appendix~\ref{app:goe}), the proofs of
Proposition~\ref*{prop:sigmahat}, Lemmas~\ref*{lem:ratio}, \ref*{lem:layers} and~\ref*{lem:mono} and
Theorem~\ref*{thm:sqrtd} of the main paper, the matching bound of Theorem~\ref*{thm:exponent}, the limit of the constant $A$ as the dimension grows and the
growth of the constant $C_{k}$ of Theorem~\ref*{thm:sqrtd} with the order $k$
(Appendix~\ref{app:proofs}), the normalisations of the literature translated into our variables
(Appendix~\ref{app:norm}), and the closed forms of the functions $G_{r}$ and $H_{r}$ of the two statistics (Appendix~\ref{app:closed}). Statements and equations are numbered by appendix. References to the main
paper are marked as such.

\section{The Kac--Rice formula with marks}\label{app:kr}

Let $\lambda\ge0$. Laws and expectations are under $\PP_{\lambda}$.
Let $f:M\times\mathbb R^{m}\to\mathbb R$ be continuous and write $f_{t}=f(t,Y)$. A mark attached to a critical point $t$ is an event
$\mathcal M_{t}$ of the form $\{f_{t}\ne0\}$, an open condition, or $\{f_{t}=0\}$, a closed condition. With $f_{t}=1$
the first form is the mark that always holds. Let $B\subseteq M$ be Borel and $u\in\mathbb R$. We use the formula
\begin{equation}\label{eq:markedkr}
\begin{split}
&\E\,\#\big\{t\in B:\nabla Z(t)=0,\ Z(t)\ge u,\ \mathcal M_{t}\big\}\\
&\qquad=\int_{B}\E\Big[\big|\det\nabla^{2}Z(t)\big|\ind\{Z(t)\ge u\}\ind_{\mathcal M_{t}}\ \Big|\ \nabla Z(t)=0\Big]\,p_{\nabla Z(t)}(0)\,\mathrm{d}t .
\end{split}
\end{equation}
With $\{Z(t)>u\}$ in place of $\{Z(t)\ge u\}$ and a mark of the first form, this is the Kac--Rice formula with a
lower semi-continuous weight of \citet[Theorem~6.1]{armentano2023general}, an
extension by monotone convergence of the weighted formula \citep[Theorem~6.4]{azais2009level}. For a Gaussian
field it requires only $C^{2}$ paths of $Z$ and a density of $\nabla Z(t)$ continuous and bounded on compacts, both
immediate here, and it counts all zeros, degenerate ones included. We apply it on the charts of a finite atlas,
to the gradient of $Z$ read in the chart. Given $\nabla Z(t)=0$, the integrand of the theorem in a chart
is~$|\det\nabla^{2}Z(t)|\,p_{\nabla Z(t)}(0)\,\mathrm{d}t$, whatever the chart: the Jacobian of the chart cancels
between the determinant, the density and the volume. The pair $(\nabla Z(t),Z)$ is
degenerate, so we check the condition of the theorem on the conditional law directly. Let $\Pi_{t}$ be the
orthogonal projection onto the kernel of the linear map $Y\mapsto\nabla Z(t)$. The covariance of $Y$ is
$\sigma^{2}\mathrm{Id}_{m}$. So, by Gaussian regression, the law of $Y$ given $\nabla Z(t)=x$ is the law of
$\Pi_{t}Y+y_{x}$, where $y_{x}$ is linear in~$x$. This law is weakly continuous in $x$. The path of $Z$ is linear in $Y$, so its
conditional law is weakly continuous in $x$ too, as the theorem asks. At $x=0$ it is the
conditional law of~\eqref{eq:markedkr} and of Lemma~\ref*{lem:transfer} of the main paper. The weight
$\ind\{Z(t)>u\}\ind\{f_{t}\ne0\}$ is the indicator of an open subset of $M\times\mathbb R^{m}$. A homogeneous
polynomial that vanishes on an open subset of the sphere vanishes on an open cone, hence everywhere. So, as in
Section~\ref*{sub:rank} of the main paper, $Y$ is a linear function of the values of $Z$ at $m$ points of the
chart. Hence this weight is lower semi-continuous in $t$ and in the path for the weak topology of the continuous
functions on the chart domain, for which point evaluations are continuous.

The closed conditions are reached by dominated convergence. Let $N_{Z}$ be the number of critical points of $Z$ on
$M$. The theorem with the weight $1$ gives $\E N_{Z}=\int_{M}D(t)\,\mathrm{d}t$, where
$D(t)=\E\big[|\det\nabla^{2}Z(t)|\mid\nabla Z(t)=0\big]\,p_{\nabla Z(t)}(0)$. The function $D$ is bounded on~$M$.
Indeed, for a unit vector $v\perp t$, two derivatives of $\theta\mapsto Z(\cos\theta\,t+\sin\theta\,v)$ at $0$ give
$v^{\top}\nabla^{2}Z(t)v=-kZ(t)+k(k-1)\langle Y,t^{\otimes(k-2)}\otimes v\otimes v\rangle$. So
$|v^{\top}\nabla^{2}Z(t)v|\le k^{2}\|Y\|$ and $|\det\nabla^{2}Z(t)|\le k^{2(d-1)}\|Y\|^{d-1}$. Given $\nabla Z(t)=0$, $Y$ has
the law of $\Pi_{t}Y$, and $\|\Pi_{t}Y\|\le\|Y\|$. So
$\E\big[\|Y\|^{d-1}\mid\nabla Z(t)=0\big]\le\E\|Y\|^{d-1}$. Also $p_{\nabla Z(t)}(0)\le(2\pi k\sigma^{2})^{-(d-1)/2}$,
as $\nabla Z(t)$ has the covariance $k\sigma^{2}\mathrm{Id}_{d-1}$. Hence, for every $t\in M$,
\[
D(t)\ \le\ k^{2(d-1)}\,\E\|Y\|^{d-1}\,(2\pi k\sigma^{2})^{-(d-1)/2},
\]
and $\E N_{Z}<\infty$. Let $\epsilon>0$. Put $\mathcal M_{t}^{\epsilon}=\{|f_{t}|<\epsilon\}$ if
$\mathcal M_{t}=\{f_{t}=0\}$, and $\mathcal M_{t}^{\epsilon}=\mathcal M_{t}$ if $\mathcal M_{t}=\{f_{t}\ne0\}$. In both
cases $\mathcal M_{t}^{\epsilon}$ is a mark of the first form, as $\{|f_{t}|<\epsilon\}$ is the event that
$(\epsilon-|f_{t}|)^{+}\ne0$. So the theorem gives~\eqref{eq:markedkr} with $\{Z(t)>u-\epsilon\}$ and
$\mathcal M_{t}^{\epsilon}$ in place of $\{Z(t)\ge u\}$ and $\mathcal M_{t}$. As $\epsilon$ decreases to $0$, these
two conditions decrease to $\{Z(t)\ge u\}$ and $\mathcal M_{t}$. On the left-hand side, the set of critical points
that is counted decreases to the set counted in~\eqref{eq:markedkr}. It is finite when $N_{Z}$ is finite, so the
count converges almost surely, and it is at most $N_{Z}$, which is integrable. On the right-hand side, the variable
under the conditional expectation decreases to that of~\eqref{eq:markedkr}. It is at most $|\det\nabla^{2}Z(t)|$,
which is integrable given $\nabla Z(t)=0$, and the integrand in $t$ is at most $D(t)$. Both sides converge by
dominated convergence. So~\eqref{eq:markedkr} holds for the two forms of marks. Lemma~\ref*{lem:transfer} of the
main paper uses it with a closed condition.

In our model the maximum of $Z$ is almost surely attained at a single point. Indeed the
drift~$\lambda\langle e,\cdot\rangle^{k}=\lambda\Cov(X(\cdot),X(e))$ is
$\lambda$ times the reproducing kernel of $X$ at $e$, hence an element of its Cameron--Martin space, so the law of
$Z$ is equivalent to that of $\sigma X$; and for~$k$ odd $\Var(X(s)-X(t))=2(1-\langle s,t\rangle^{k})$ vanishes
only at $s=t$, so the maximum of $X$ is attained at a single point \citep[Theorem~3]{lifshits1983absolute}. Every
almost sure property of $X$ transfers to $Z$, in particular the non-degeneracy of the critical points and
$\PP\{\sup Z=u\}=0$.

\section{The layer decomposition and the drift}\label{app:layers}

Fix $t$ and write $s=\cos\theta\,t+\sin\theta\,v$ with $v\in\Sph^{d-2}\subset t^{\perp}$ and $r=\cos\theta$, and
set $M_{t}=\Sph^{d-1}\setminus\{\pm t\}$. Standard computations give $\Cov(X(s),X(t))=r^{k}$,
$\Cov(X(s),\partial_{a}X(t))=k\,r^{k-1}\sin\theta\,v_{a}$, $\Cov(X(t),\nabla X(t))=0$ and
$\Cov(\nabla X(t))=k\,\mathrm{Id}_{d-1}$, so the Gaussian regression of $X(s)$ on $(X(t),\nabla X(t))$ reads
\begin{subequations}
\begin{equation}\label{eq:regression}
X(s)=r^{k}X(t)+r^{k-1}\sin\theta\,\langle v,\nabla X(t)\rangle+R^{t}(s),
\end{equation}
and the helix process is $\zeta^{t}(s)=R^{t}(s)/(1-r^{k})$ on $M_{t}$. The matrix
\begin{equation}\label{eq:G}
G:=\frac{\nabla^{2}X(t)+k\,X(t)\,\mathrm{Id}_{d-1}}{\sqrt{2k(k-1)}}\;\sim\;\GOE(d-1)
\end{equation}
is independent of $(X(t),\nabla X(t))$ and of the layers of index $j\ge3$ below. The
binomial expansion of $\langle s,\cdot\rangle^{k}$ and the regression on $(Z(t),\nabla Z(t))$ give the
\emph{layer decomposition} of the helix: with $w_{j}(\theta)=\binom kj^{1/2}r^{k-j}\sin^{j}\theta/(1-r^{k})$,
\begin{equation}
\zeta^{t}(s)=\sum_{j=2}^{k}w_{j}(\theta)\,Z_{j}(v),
\end{equation}
\end{subequations}
where $Z_{2}(v)=v^{\top}Gv$ with $G\sim\GOE(d-1)$ and, for $j\ge3$, $Z_{j}$ is an independent unit-variance
Kostlan--Shub--Smale field of degree $j$ on $\Sph^{d-2}$; all of them are independent of~$(Z(t),\nabla Z(t))$.

Under the alternative the same computation, applied to $\lambda\langle e,\cdot\rangle^{k}$ with
$c=\langle e,t\rangle$, replaces each layer by a shifted one: the terms
$j=0$ and $j=1$ of the expansion are exactly what the regression removes, and
\begin{subequations}
\begin{equation}\label{eq:drift}
\mu_{t}(s)=\sum_{j=2}^{k}w_{j}(\theta)\,\beta_{j}(c)\,\langle e_{\perp},v\rangle^{j},\qquad
\beta_{j}(c)=\lambda\binom kj^{1/2}c^{k-j}(1-c^{2})^{j/2},
\end{equation}
an identity for every $k$. Cauchy--Schwarz, $\sum_{j}w_{j}^{2}\le\rho^{-2}=2(k-1)/k$,
$\binom kj\le\binom k2\binom{k-2}{j-2}$ for~$j\ge2$ and $\binom k2^{1/2}/\rho=k-1$ give the sharp drift bound
\begin{equation}\label{eq:driftbound}
\sup_{M_{t}}|\mu_{t}|\ \le\ (k-1)\,\lambda\,(1-c^{2}) .
\end{equation}
\end{subequations}
At the pole the drift is $\beta_{2}(c)\langle e_{\perp},v\rangle^{2}/\rho$, whence $G_{\lambda}=G+\beta_{2}(c)e_{\perp}e_{\perp}^{\top}$
in~(\ref*{eq:omega}) of the main paper. The layer decomposition is proved in \citet{kss}; Lemmas~\ref{lem:variance} and~\ref{lem:poles}
below are quoted from there without proof.

Let $\Mb=[0,\pi]\times\Sph^{d-2}$ be the blow-up of the sphere at $\pm t$, whose boundary components~$\{0\}\times\Sph^{d-2}$ and $\{\pi\}\times\Sph^{d-2}$ are the pole fibres.

\begin{lemma}[Variance; see \citet{kss}]\label{lem:variance}
$V(\theta):=\Var\zeta^{t}(s)=\|w(\theta)\|^{2}=\dfrac{1-r^{2k}-k\,r^{2k-2}(1-r^{2})}{(1-r^{k})^{2}}\le\rho^{-2}=\dfrac{2(k-1)}{k}$ on $(0,\pi)$, with equality only in the limit $\theta\to0$ when $k$ is odd.
\end{lemma}

\begin{lemma}[Poles; see \citet{kss}]\label{lem:poles}
(i) Each $w_{j}$ is $C^{\infty}$ on~$[0,\pi)$, and as $\theta\to0$
\[
w_{j}(\theta)=\tfrac2k\tbinom kj^{1/2}\theta^{\,j-2}\big(1+O(\theta^{2})\big),\qquad
w_{2}(0)=\rho^{-1}.
\]
(ii) $\zeta^{t}$ extends to an almost surely continuous field $\tilde\zeta(\theta,v)=\sum_{j}w_{j}(\theta)Z_{j}(v)$ on $\Mb$, with $\tilde\zeta(0,v)=v^{\top}Gv/\rho$, and $\tilde\zeta(\pi,\cdot)\equiv0$ for $k$ odd.
\end{lemma}

The last lemma of this appendix is not taken from \citet{kss}. The proof of Lemma~\ref*{lem:pos} of the main paper
uses it. Its own proof uses Lemma~\ref*{lem:layers} of the main paper, which is stated after
Lemma~\ref*{lem:pos}. The proof of Lemma~\ref*{lem:layers}, in Appendix~\ref{app:proofs}, uses neither
Lemma~\ref*{lem:transfer} nor Lemma~\ref*{lem:pos} of the main paper. For a general field,
\citet[Appendix~C]{azais2024second} write that the regularity of $\lambda_{2}^{t}$ as a
function of $t$ is difficult to control. Here the model has finite rank. Lemma~\ref{lem:lam2cont} obtains
the continuity of $\lambda_{2}^{t}$ in $(t,Y)$ from the finite sum of Lemma~\ref*{lem:layers} of the main paper.

\begin{lemma}[Continuity of the second maximum at a point]\label{lem:lam2cont}
Let $k$ be odd. The maps $(t,Y)\mapsto\lambda_{2}^{t}$ and $(t,Y)\mapsto\lmax(\Omega_{t})$ are continuous on
$M\times\mathbb R^{m}$.
\end{lemma}

\begin{proof}
For $t\in M$, $v\in\mathbb R^{d}$ and $2\le j\le k$ let
$\Psi_{j}^{t}(v)=\binom kj^{1/2}\operatorname{Sym}\big(t^{\otimes(k-j)}\otimes v^{\otimes j}\big)$. For a unit vector
$v\perp t$ this is the tensor $\Psi_{j}(v)$ of Lemma~\ref*{lem:layers} of the main paper at the point $t$. Its
coordinates are polynomials in $(t,v)$. Let $S_{t}$ be the set of the unit vectors $v\perp t$. For $t\in M$,
$\theta\in[0,\pi]$, $v\in\mathbb R^{d}$ and $Y\in\mathbb R^{m}$ let
\begin{subequations}
\begin{equation}\label{eq:Fext}
F(t,\theta,v,Y)=\sum_{j=2}^{k}w_{j}(\theta)\,\big\langle Y,\Psi_{j}^{t}(v)\big\rangle .
\end{equation}
The weights $w_{j}$ are continuous on $[0,\pi]$, with $w_{2}(0)=\rho^{-1}$, $w_{j}(0)=0$ for $j\ge3$ and
$w_{j}(\pi)=0$. On $[0,\pi)$ this is Lemma~\ref{lem:poles}(i). At $\pi$ the denominator $1-\cos^{k}\theta$ tends to
$2$, as $k$ is odd, and~$\sin^{j}\theta$ tends to $0$. So $F$ is continuous in $(t,\theta,v,Y)$. By
Lemma~\ref*{lem:layers} of the main paper, for every~$Y$, $F(t,\theta,v,Y)=Z^{|t}(\cos\theta\,t+\sin\theta\,v)$ for
$(\theta,v)\in(0,\pi)\times S_{t}$, and $F(t,0,v,Y)=v^{\top}\Omega_{t}v$ for $v\in S_{t}$. The points
$\cos\theta\,t+\sin\theta\,v$ with $(\theta,v)\in(0,\pi)\times S_{t}$ are the points of $M_{t}$,
and~$(0,\pi)\times S_{t}$ is dense in the compact set $[0,\pi]\times S_{t}$. So
\begin{equation}\label{eq:lam2max}
\lambda_{2}^{t}=\max_{(\theta,v)\in[0,\pi]\times S_{t}}F(t,\theta,v,Y),\qquad
\lmax(\Omega_{t})=\max_{v\in S_{t}}F(t,0,v,Y).
\end{equation}
The sets in~\eqref{eq:lam2max} depend on $t$. Fix $t_{0}\in M$ and let $U=\{t\in M:\langle t_{0},t\rangle>-1\}$, an
open neighbourhood of $t_{0}$. For $t\in U$ let
\begin{equation}\label{eq:Qt}
Q_{t}=\Big(\mathrm{Id}-\frac{(t_{0}+t)(t_{0}+t)^{\top}}{1+\langle t_{0},t\rangle}\Big)\big(\mathrm{Id}-2\,t_{0}t_{0}^{\top}\big).
\end{equation}
Both factors are reflections, the first one because $\|t_{0}+t\|^{2}=2(1+\langle t_{0},t\rangle)$. The second
maps~$t_{0}$ to $-t_{0}$, and the first maps $-t_{0}$ to $t$. So $Q_{t}$ is orthogonal, $Q_{t}t_{0}=t$, and $Q_{t}$ maps
$S_{t_{0}}$ onto~$S_{t}$. It is continuous in $t\in U$. By~\eqref{eq:lam2max}, for $t\in U$,
\begin{equation}\label{eq:lam2maxfixed}
\lambda_{2}^{t}=\max_{(\theta,v)\in[0,\pi]\times S_{t_{0}}}F(t,\theta,Q_{t}v,Y),\qquad
\lmax(\Omega_{t})=\max_{v\in S_{t_{0}}}F(t,0,Q_{t}v,Y).
\end{equation}
\end{subequations}
The function $g(t,Y,\theta,v)=F(t,\theta,Q_{t}v,Y)$ is continuous on $U\times\mathbb R^{m}\times K$, where
$K=[0,\pi]\times S_{t_{0}}$ is compact and does not depend on $(t,Y)$. Let $(t_{n},Y_{n})\to(t,Y)$ in
$U\times\mathbb R^{m}$. The set $\big(\{(t_{n},Y_{n}):n\ge1\}\cup\{(t,Y)\}\big)\times K$ is compact, so $g$ is
uniformly continuous on it, and~$g(t_{n},Y_{n},\cdot)\to g(t,Y,\cdot)$ uniformly on $K$. The maxima of two functions
over a set differ by at most the supremum of their difference. Hence the two maxima in~\eqref{eq:lam2maxfixed} at
$(t_{n},Y_{n})$ converge to the maxima at $(t,Y)$. So both maps are continuous on $U\times\mathbb R^{m}$, and
$t_{0}$ is arbitrary.
\end{proof}

\section{The Gaussian Orthogonal Ensemble, and the count outside a cap}\label{app:goe}

\begin{lemma}[Determinant of a rank-one spiked GOE]\label{lem:spikedet}
Let $G\sim\GOE(N)$, let $q$ be a unit vector of $\mathbb R^{N}$ and $\beta\in\mathbb R$, and write
$D_{n}(z)=\E\det(z\,\mathrm{Id}_{n}-G_{n})$ for $G_{n}\sim\GOE(n)$, $D_{0}=1$. Then
\begin{subequations}
\begin{equation}
\E\det\big(z\,\mathrm{Id}_{N}-G-\beta\,qq^{\top}\big)=D_{N}(z)-\beta\,D_{N-1}(z)\qquad\text{for every }z\in\mathbb R,
\end{equation}
and $\tfrac34z^{n}\le D_{n}(z)\le z^{n}$ for $z\ge n$, so that
\begin{equation}
\E\big|\det(z\,\mathrm{Id}_{N}-G-\beta qq^{\top})\big|\ \ge\ \tfrac12z^{N}\qquad\text{for }z\ge N,\ |\beta|\le z/4 .
\end{equation}
\end{subequations}
\end{lemma}

\begin{proof}
The law of $G$ is invariant under $G\mapsto R^{\top}GR$, $R\in O(N)$. Picking $R$ with $Re_{1}=q$ gives
$R^{\top}(qq^{\top})R=e_{1}e_{1}^{\top}$, so one may take $q=e_{1}$. The matrix determinant lemma
$\det(M-\beta e_{1}e_{1}^{\top})=\det M-\beta\,e_{1}^{\top}\mathrm{adj}(M)e_{1}$ with $M=z\,\mathrm{Id}-G$
gives $\det(z\,\mathrm{Id}-G-\beta e_{1}e_{1}^{\top})=\det(z\,\mathrm{Id}-G)-\beta\det(z\,\mathrm{Id}_{N-1}-G')$,
the $(1,1)$ entry of the adjugate being the corresponding principal minor; $G'=(G_{ij})_{2\le i,j\le N}$ is a
$\GOE(N-1)$ with the same normalisation, whence the identity. For the bracket, $\E\det G'=\sum_{\pi}\mathrm{sgn}(\pi)\E\prod_{i}G'_{i\pi(i)}$
survives only for $\pi$ a product of disjoint transpositions, which gives
$D_{n}(z)=\sum_{i\le n/2}(-1)^{i}t_{i}$ with $t_{i}=\frac{n!}{i!(n-2i)!4^{i}}z^{n-2i}\ge0$ and
$t_{i+1}/t_{i}=\frac{(n-2i)(n-2i-1)}{4(i+1)z^{2}}\le\frac{n^{2}}{4z^{2}}\le\frac14$ for $z\ge n$; an
alternating series with decreasing terms is bracketed by its first two partial sums. Hence
$D_{N}-\beta D_{N-1}\ge\frac34z^{N}-\frac z4z^{N-1}=\frac12z^{N}$ for $0\le\beta\le z/4$, while for
$\beta\le0$ the subtracted term is non-positive and the bound is $\frac34z^{N}$; and $\E|\det|\ge\E\det$.
\end{proof}

We also use, from \citet{kss}, the Gaussian concentration of the GOE norm, $\PP\{\|G\|>\sqrt{2N}+s\}\le2e^{-s^{2}/2}$
for $s\ge0$. It gives the moment bound $\|\,\|G\|\,\|_{q}\le\sqrt{2N}+\tfrac56\sqrt q$ for $q\ge8$: integrating the
tail, $\E(\|G\|-\sqrt{2N})_{+}^{q}\le2q\int_{0}^{\infty}s^{q-1}e^{-s^{2}/2}\,\mathrm{d}s=q\,2^{q/2}\Gamma(q/2)$, whose
$q$-th root is at most $\tfrac56\sqrt q$ for $q\ge8$ ($2.29$ against $2.36$ at $q=8$, and the ratio of the two sides
decreases in $q$), and Minkowski's inequality in $L^{q}$ gives the bound.

\begin{lemma}[Upper bound on the count outside a cap]\label{lem:Aupper}
Let $\theta_{1}\in(0,\pi/2)$, let $\beta_{\star}:=\lambda\binom k2^{1/2}\sin^{2}\theta_{1}=\rho(k-1)\lambda\sin^{2}\theta_{1}$,
let $u$ satisfy $(d-1)/\rho<u<\lambda\cos^{k}\theta_{1}$ and $\rho u\ge\max\{4\beta_{\star},\ \beta_{\star}+\sqrt{2N}\}$, and set
\begin{subequations}
\begin{equation}
\varrho=1-\inf_{0<\theta\le\theta_{1}}\cos^{2k-2}\theta\Big(\frac{\sin\theta}{\theta}\Big)^{2} .
\end{equation}
Then, for every $b\in(0,\lambda^{2}]$,
\begin{equation}
\delta^{\lambda}\big(u;\mathcal B_{b}^{\,c}\big)\ \le\ \frac{1+\epsilon_{A}}{(1-\varrho)^{\frac{d-1}{2}}}\;
\PP\big\{\chi^{2}_{d-1}>k\,b\,(1-\varrho)\big\}\ +\ \mathrm{Rem}_{\rm far},
\end{equation}
where $\mathrm{Rem}_{\rm far}=\delta^{\lambda}(u;\{\theta>\theta_{1}\})$ and
$1+\epsilon_{A}=\big(1+\beta_{\star}/(\rho u)+\varepsilon_{\det}\big)(1+r_{N})$, the three remainders being, with
$q_{\star}=\max(2N,8)$ and $\ell_{1}=\lambda\cos^{k}\theta_{1}$,
\begin{equation}
r_{N}\le e^{N^{2}/(2\ell_{1}^{2})}-1,\qquad
\varepsilon_{\det}\le2\sqrt2\Big(\tfrac54+\frac{\sqrt{2N}+\frac56\sqrt{q_{\star}}}{N}\Big)^{N}
e^{-\frac14(\rho u-\beta_{\star}-\sqrt{2N})^{2}} .
\end{equation}
\end{subequations}
\end{lemma}

\begin{proof}
Three steps. \emph{The determinant.} On the event that all eigenvalues of $G_{\lambda}$ lie below~$z$ the determinant
$\det(z\,\mathrm{Id}-G_{\lambda})$ is positive, so
$A_{\lambda}(z;c)=\E\det(z\,\mathrm{Id}-G_{\lambda})+2\E[(\det)^{-}]$, and
$\E\det\le D_{N}(z)+|\beta_{2}(c)|D_{N-1}(z)\le z^{N}(1+|\beta_{2}(c)|/z)$ by
Lemma~\ref{lem:spikedet}. For the negative part, Cauchy--Schwarz and
$\|\det(z\,\mathrm{Id}-G_{\lambda})\|_{2}\le(z+|\beta_{2}(c)|+\|\,\|G\|\,\|_{q_{\star}})^{N}$
with $\|\,\|G\|\,\|_{q_{\star}}\le\sqrt{2N}+\tfrac56\sqrt{q_{\star}}$ give
\[2\E[(\det)^{-}]\le2\big(z+\tfrac z4+\sqrt{2N}+\tfrac56\sqrt{q_{\star}}\big)^{N}\sqrt{2}\,
e^{-\frac14(z-|\beta_{2}(c)|-\sqrt{2N})^{2}}\le\varepsilon_{\det}\,z^{N},\]
using $|\beta_{2}(c)|\le z/4$ and $z\ge N$. \emph{The Gaussian integral.} Extending the range to $\mathbb R$ and expanding,
$\int_{\mathbb R}x^{N}\varphi(x-\ell_{0})\,\mathrm{d}x=\E(\ell_{0}+\xi)^{N}=\ell_{0}^{N}(1+r_{N})$ with $\ell_{0}=\lambda c^{k}$, the odd moments
vanishing and
\begin{subequations}
\begin{equation}
r_{N}=\sum_{i\ge1}\binom N{2i}(2i-1)!!\,\ell_{0}^{-2i}\le\sum_{i\ge1}\frac{1}{i!}\Big(\frac{N^{2}}{2\ell_{0}^{2}}\Big)^{i}
=e^{N^{2}/(2\ell_{0}^{2})}-1,
\end{equation}
since $\binom N{2i}(2i-1)!!\le N^{2i}/(2i)!!=N^{2i}/(2^{i}i!)$; hence
\begin{equation}
\int_{u}^{\infty}A_{\lambda}(\rho x;c)\varphi(x-\ell_{0})\,\mathrm{d}x\le(\rho\ell_{0})^{N}(1+\epsilon_{A}).
\end{equation}
\end{subequations}
The extension does not decrease the integral. For $N$ even the integrand is non-negative. For~$N$ odd the contribution of
$x<0$ is at least $-\varphi(\ell_{0})N!/\ell_{0}^{N+1}$, since $\varphi(x+\ell_{0})\le\varphi(\ell_{0})e^{-x\ell_{0}}$
for $x\ge0$, while that of $(0,u)$ is at least $\varphi(\ell_{0})u^{N+1}/(N+1)$, since $\ell_{0}x-x^{2}/2\ge0$
for~$0\le x\le u\le\ell_{0}$; the latter dominates because $u\ell_{0}\ge u^{2}>N^{2}$ by $\rho u>N$ and $\rho<1$,
whence~$(u\ell_{0})^{N+1}\ge N^{2N+2}\ge(N+1)!$.
\emph{The Laplace step.} For $b\ge\lambda^{2}\sin^{2}\theta_{1}$, $\mathcal B_{b}^{c}\subseteq\{\theta>\theta_{1}\}$ and the bound
holds with $\mathrm{Rem}_{\rm far}$ alone. For $b<\lambda^{2}\sin^{2}\theta_{1}$, in the variable $\theta$, on $(\theta_{b},\theta_{1}]$
with~$\sin\theta_{b}=\sqrt b/\lambda$, use $\sin^{d-2}\theta\le\theta^{d-2}$, $\cos^{kN}\theta\le1$ and
$\tfrac{k\lambda^{2}}{2}\cos^{2k-2}\theta\sin^{2}\theta=\tfrac{w^{2}}{2}q(\theta)$
with~$q\ge1-\varrho$ and $w=\lambda\sqrt k\,\theta$, so that
\begin{subequations}
\begin{equation}
\int_{\theta_{b}}^{\theta_{1}}\le\frac{1}{(\lambda\sqrt k)^{d-1}}\int_{\sqrt{kb}}^{\infty}w^{d-2}e^{-\frac{w^{2}}{2}(1-\varrho)}\,\mathrm{d}w .
\end{equation}
The prefactor normalises exactly:
\begin{equation}
\omega_{d-2}\big(\tfrac{k-1}{\pi}\big)^{\frac{d-1}{2}}\rho^{d-1}k^{-\frac{d-1}{2}}\int_{0}^{\infty}w^{d-2}e^{-w^{2}/2}\,\mathrm{d}w
=2^{\frac{d-1}{2}}\rho^{d-1}\big(\tfrac{k-1}{k}\big)^{\frac{d-1}{2}}=1,
\end{equation}
\end{subequations}
because $\rho^{2}=k/(2(k-1))$. Rescaling
$w\mapsto w\sqrt{1-\varrho}$ turns the truncated integral into
$(1-\varrho)^{-(d-1)/2}\PP\{\chi^{2}_{d-1}>kb(1-\varrho)\}$; the range $\theta>\theta_{1}$ is
$\mathrm{Rem}_{\rm far}$.
\end{proof}

\begin{lemma}[Simple bounds]\label{lem:simple}
For $\theta_{1}\in(0,\pi/2)$, the quantity $\varrho$ of Lemma~\ref{lem:Aupper} satisfies
\begin{subequations}
\begin{equation}
\varrho=1-\cos^{2k-2}\theta_{1}\Big(\frac{\sin\theta_{1}}{\theta_{1}}\Big)^{2}\ \le\ \Big(k-\frac23\Big)\theta_{1}^{2}.
\end{equation}
Under the hypotheses of Lemma~\ref{lem:Aupper},
\begin{equation}
1+\epsilon_{A}\le\Big(1+\frac{(k-1)\lambda\sin^{2}\theta_{1}}{u}+\varepsilon_{\det}\Big)e^{\frac{N^{2}}{2\ell_{1}^{2}}},
\qquad \varepsilon_{\det}\le3\,e^{\frac54N-\frac14(\rho u-\beta_{\star}-\sqrt{2N})^{2}}.
\end{equation}
If also $\rho u\ge4\sqrt{2N}$, then $\varepsilon_{\det}\le3e^{\frac54N-\frac{\rho^{2}u^{2}}{16}}\le3e^{-\frac34N}$ and
\begin{equation}
\epsilon_{A}\le\frac{(k-1)\lambda\sin^{2}\theta_{1}}{u}+\frac{3N^{2}}{2\ell_{1}^{2}}+3e^{\frac54N-\frac{\rho^{2}u^{2}}{16}} .
\end{equation}
\end{subequations}
\end{lemma}

\begin{proof}
On $(0,\pi/2)$ the functions $\cos^{2k-2}\theta$ and $(\sin\theta/\theta)^{2}$ are positive and decreasing, so the
infimum that defines $\varrho$ is attained at $\theta_{1}$. Next, $\sin x\ge x\cos^{1/3}x$ on $[0,\pi/2)$. Indeed,
with~$c=\cos x$ and $f(x)=c^{-1/3}\sin x-x$, one has $f(0)=0$ and $f'(x)=\frac23c^{2/3}+\frac13c^{-4/3}-1\ge0$, by the
weighted arithmetic-geometric inequality and $(c^{2/3})^{2/3}(c^{-4/3})^{1/3}=1$. Hence~$1-\varrho\ge\cos^{p}\theta_{1}$
with $p=2k-\frac43\ge1$. With $1-c^{p}\le p(1-c)$ for $c\in[0,1]$ and~$1-\cos\theta_{1}\le\theta_{1}^{2}/2$, this gives
$\varrho\le(k-\frac23)\theta_{1}^{2}$.

By Lemma~\ref{lem:Aupper}, $1+\epsilon_{A}=(1+a+\varepsilon_{\det})(1+r_{N})$ with
$a=\beta_{\star}/(\rho u)=(k-1)\lambda\sin^{2}\theta_{1}/u$ and $1+r_{N}\le e^{N^{2}/(2\ell_{1}^{2})}$. The constant in
the bound of Lemma~\ref{lem:Aupper} on $\varepsilon_{\det}$ is $2\sqrt2B_{N}^{N}$ with~$B_{N}=\frac54+m_{N}/N$ and
$m_{N}=\sqrt{2N}+\frac56\sqrt{q_{\star}}$. One has $B_{2}=\frac94+\frac56\sqrt2<3.43$, $B_{3}<2.86$ and, for $N\ge4$,
where $q_{\star}=2N$, $B_{N}=\frac54+\frac{11}{6}\sqrt{2/N}<2.55$. Since $e^{5/4}>3.49$, $2\sqrt2B_{N}^{N}\le3e^{5N/4}$
for every $N\ge2$. This proves the second display.

Let now $\rho u\ge4\sqrt{2N}$. With $\rho u\ge4\beta_{\star}$, the gap $g=\rho u-\beta_{\star}-\sqrt{2N}$ is at least
$\rho u/2$, so~$g^{2}/4\ge\rho^{2}u^{2}/16\ge2N$. This gives the two bounds on $\varepsilon_{\det}$, and
$\varepsilon_{\det}\le3e^{-3/2}<0.67$. Put~$x=N^{2}/(2\ell_{1}^{2})$. Since $\ell_{1}>u>N/\rho$ and $\rho^{2}\le\frac34$,
$x<\rho^{2}/2\le\frac38$, so $r_{N}\le e^{x}-1\le xe^{x}\le e^{3/8}x<0.73\,N^{2}/\ell_{1}^{2}$. Also $a\le\frac14$, so
$1+a+\varepsilon_{\det}<1.92$, and $\epsilon_{A}=a+\varepsilon_{\det}+r_{N}(1+a+\varepsilon_{\det})\le
a+\varepsilon_{\det}+1.41\,N^{2}/\ell_{1}^{2}$, which gives the last display.
\end{proof}

\begin{lemma}[The far zone]\label{lem:far}
Let $\theta_{1}$ and $u$ satisfy the hypotheses of Lemma~\ref{lem:Aupper}, and let $3\rho u\ge4\sqrt{2N}$. Define
$\theta_{3}\in(0,\pi/2)$ by $\sin^{2}\theta_{3}=u/(4(k-1)\lambda)$, and let $\ell_{3}=\lambda\cos^{k}\theta_{3}$,
$h(\theta)=\cos^{2k-2}\theta\sin^{2}\theta$, $m_{N}=\sqrt{2N}+\frac56\sqrt{q_{\star}}$,
$K_{N}=2\sqrt2\,(\frac54+m_{N}/N)^{N}$, and $\omega_{d-1}=2\pi^{d/2}/\Gamma(\frac d2)$, the area of $\Sph^{d-1}$. Put
\begin{subequations}
\begin{align}
C_{M}&=\Big(\frac54+K_{N}e^{-\frac14(\frac34\rho u-\sqrt{2N})^{2}}\Big)\frac{e^{N(N+1)/(2\ell_{3}^{2})}}{1-(k-1)\tan^{2}\theta_{3}},\\
F_{\rm far}&=\omega_{d-1}\Big(\frac{k-1}{\pi}\Big)^{\frac N2}\Big(1+\frac{(k-1)\lambda}{u}\Big)
\big(\rho\lambda+\rho\sqrt N+m_{N}\big)^{N}e^{-E_{\rm far}},\\
E_{\rm far}&=\sup_{\theta_{3}\le\theta_{s}\le\pi/2}\min\Big\{\frac{k\lambda^{2}}{2}\min\big(h(\theta_{3}),h(\theta_{s})\big),
\ \frac12\big(u-\lambda\cos^{k}\theta_{s}\big)_{+}^{2}\Big\}.
\intertext{Then $\theta_{1}\le\theta_{3}$, $(k-1)\tan^{2}\theta_{3}\le(k-1)/(4k-5)$, and}
\mathrm{Rem}_{\rm far}&=\delta^{\lambda}\big(u;\{\theta>\theta_{1}\}\big)\ \le\ C_{M}\,\PP\big\{\chi^{2}_{N}>k\lambda^{2}h(\theta_{1})\big\}+F_{\rm far}.
\end{align}
\end{subequations}
If also $u\ge3\lambda/4$, then $E_{\rm far}\ge\lambda^{2}/32$ and $\cos^{k}\theta_{3}\ge13/16$.
\end{lemma}

\begin{proof}
The hypothesis $\rho u\ge4\beta_{\star}$ of Lemma~\ref{lem:Aupper} reads $\sin^{2}\theta_{1}\le\sin^{2}\theta_{3}$, so
$\theta_{1}\le\theta_{3}$. Since $u<\lambda$, $s_{3}:=\sin^{2}\theta_{3}<1/(4(k-1))$, and
$(k-1)\tan^{2}\theta_{3}=(k-1)s_{3}/(1-s_{3})\le(k-1)/(4k-5)<1$.

First, for $\ell>0$ and $\xi\sim N(0,1)$, $\E|\ell+\xi|^{N}\le\ell^{N}e^{N(N+1)/(2\ell^{2})}$. For $N$ even,
$\E(\ell+\xi)^{N}=\ell^{N}(1+r_{N})\le\ell^{N}e^{N^{2}/(2\ell^{2})}$ by the computation in the proof of
Lemma~\ref{lem:Aupper}, which holds for every $\ell>0$. For $N$ odd, H\"older's inequality gives
$\E|\ell+\xi|^{N}\le(\E(\ell+\xi)^{N+1})^{N/(N+1)}\le\ell^{N}e^{N(N+1)/(2\ell^{2})}$.

Second, for $z>0$ and every $c$, $A_{\lambda}(z;c)\le(z+m_{N})^{N}(1+|\beta_{2}(c)|/z)$. The matrix determinant
lemma in the proof of Lemma~\ref{lem:spikedet} gives
\begin{subequations}
\begin{equation}
\det(z\,\mathrm{Id}-G_{\lambda})=\det(z\,\mathrm{Id}-G)-\beta_{2}(c)\det(z\,\mathrm{Id}_{N-1}-G'),
\end{equation}
with $G'$ a principal submatrix of $G$, so that $\|G'\|\le\|G\|$. For $n\in\{N-1,N\}$ and $G_{n}$ the matrix of size
$n$ among $G'$ and $G$, Minkowski's inequality gives $\E|\det(z\,\mathrm{Id}_{n}-G_{n})|\le(z+m_{N})^{n}$, since
$n\le q_{\star}$ and $\|\,\|G\|\,\|_{q_{\star}}\le m_{N}$. Add the two terms and use
$(z+m_{N})^{N-1}\le(z+m_{N})^{N}/z$.

Third, let $(k-1)\tan^{2}\theta_{\max}=1$. On $(0,\theta_{\max})$ the map
$w(\theta)=\lambda\sqrt k\cos^{k-1}\theta\sin\theta$ is increasing, with
$w'(\theta)=\lambda\sqrt k\cos^{k}\theta\,(1-(k-1)\tan^{2}\theta)$ and $w^{2}=k\lambda^{2}h(\theta)$. Also
$\sin^{N-1}\theta\cos^{kN}\theta=(w/(\lambda\sqrt k))^{N-1}\cos^{N+k-1}\theta$, and the normalisation in the proof of
Lemma~\ref{lem:Aupper} reads $\omega_{d-2}((k-1)/\pi)^{N/2}\rho^{N}k^{-N/2}=1/(2^{N/2-1}\Gamma(N/2))$. This gives the
exact identity
\begin{equation}
\omega_{d-2}\Big(\frac{k-1}{\pi}\Big)^{\frac N2}\sin^{N-1}\theta\,(\rho\lambda\cos^{k}\theta)^{N}
e^{-\frac{k\lambda^{2}}{2}h(\theta)}\,\mathrm{d}\theta=f_{\chi_{N}}(w)\,\frac{\cos^{N-1}\theta}{1-(k-1)\tan^{2}\theta}\,\mathrm{d}w,
\end{equation}
\end{subequations}
where $f_{\chi_{N}}$ is the density of the norm of a standard Gaussian vector of $\mathbb R^{N}$.

In the variable $\theta$, with $c=\cos\theta$, the density of $\delta^{\lambda}$ is
\[\omega_{d-2}\Big(\frac{k-1}{\pi}\Big)^{\frac N2}\sin^{N-1}\theta\,e^{-\frac{k\lambda^{2}}{2}h(\theta)}
\int_{u}^{\infty}A_{\lambda}(\rho x;c)\varphi(x-\lambda c^{k})\,\mathrm{d}x .\]
Split $\{\theta>\theta_{1}\}$ into $(\theta_{1},\theta_{3}]$ and $(\theta_{3},\pi]$.

On $(\theta_{1},\theta_{3}]$, $c\in(0,1)$ and $|\beta_{2}(c)|\le\rho(k-1)\lambda\sin^{2}\theta_{3}=\rho u/4$. For $x\ge u$,
$z=\rho x\ge\rho u>N$, $|\beta_{2}(c)|\le z/4$ and $z-|\beta_{2}(c)|-\sqrt{2N}\ge\frac34\rho u-\sqrt{2N}\ge0$. The
determinant step of the proof of Lemma~\ref{lem:Aupper}, run with these bounds, gives
$A_{\lambda}(\rho x;c)\le(\rho x)^{N}(\frac54+\varepsilon_{3})$, with
$\varepsilon_{3}=K_{N}e^{-\frac14(\frac34\rho u-\sqrt{2N})^{2}}$. The first fact with $\ell=\lambda c^{k}\ge\ell_{3}$ gives
$\int_{u}^{\infty}A_{\lambda}(\rho x;c)\varphi(x-\lambda c^{k})\,\mathrm{d}x\le(\frac54+\varepsilon_{3})
e^{N(N+1)/(2\ell_{3}^{2})}(\rho\lambda\cos^{k}\theta)^{N}$. The interval lies in $(0,\theta_{\max})$, and there the
factor $\cos^{N-1}\theta/(1-(k-1)\tan^{2}\theta)$ of the third fact is at most $1/(1-(k-1)\tan^{2}\theta_{3})$. By the
third fact, this zone contributes at most
$C_{M}\int_{w(\theta_{1})}^{\infty}f_{\chi_{N}}=C_{M}\,\PP\{\chi^{2}_{N}>k\lambda^{2}h(\theta_{1})\}$.

On $(\theta_{3},\pi]$, the second fact with $|\beta_{2}(c)|\le\rho(k-1)\lambda$ and $z=\rho x\ge\rho u$ gives
$A_{\lambda}(\rho x;c)\le(1+(k-1)\lambda/u)(\rho x+m_{N})^{N}$ for $x\ge u$. Let $\ell=\lambda c^{k}\le\lambda$. Since
$\E\xi^{2j}=(2j-1)!!\le(2j-1)^{j}$, Lyapunov's inequality with $2j\in\{N,N+1\}$ gives $\|\xi\|_{N}\le\sqrt N$. If $\ell\ge u$, Minkowski's
inequality gives
$\int_{u}^{\infty}(\rho x+m_{N})^{N}\varphi(x-\ell)\,\mathrm{d}x\le\E(\rho|\ell+\xi|+m_{N})^{N}\le(\rho\lambda+\rho\sqrt N+m_{N})^{N}$.
If $\ell<u$, put $x=u+y$ and $g=u-\ell$. Then $\varphi(g+y)\le e^{-g^{2}/2}\varphi(y)$ for $y\ge0$, and the integral is
at most $e^{-g^{2}/2}\E(\rho u+m_{N}+\rho|\xi|)^{N}\le e^{-g^{2}/2}(\rho\lambda+\rho\sqrt N+m_{N})^{N}$. So the integrand
of~$\delta^{\lambda}$ is at most $\omega_{d-2}((k-1)/\pi)^{N/2}(1+(k-1)\lambda/u)(\rho\lambda+\rho\sqrt N+m_{N})^{N}
\sin^{N-1}\theta\,e^{-E(\theta)}$, with $E(\theta)=\frac{k\lambda^{2}}{2}h(\theta)+\frac12(u-\lambda\cos^{k}\theta)_{+}^{2}$.
Fix $\theta_{s}\in[\theta_{3},\pi/2]$. On $(0,\pi/2)$ the sign of $h'$ is that of $\cos^{2}\theta-(k-1)\sin^{2}\theta$, so
$h$ increases and then decreases, and $h\ge\min(h(\theta_{3}),h(\theta_{s}))$ on~$[\theta_{3},\theta_{s}]$. For
$\theta>\theta_{s}$, $\cos^{k}\theta<\cos^{k}\theta_{s}$, $k$ being odd, so
$(u-\lambda\cos^{k}\theta)_{+}\ge(u-\lambda\cos^{k}\theta_{s})_{+}$. Hence $E(\theta)\ge E_{\rm far}$ on
$(\theta_{3},\pi]$, and $\omega_{d-2}\int_{0}^{\pi}\sin^{N-1}\theta\,\mathrm{d}\theta=\omega_{d-1}$ gives the bound
$F_{\rm far}$ for this zone.

Let now $u\ge3\lambda/4$. Then $s_{3}\ge3/(16(k-1))$, and $\cos^{k}\theta_{3}=(1-s_{3})^{k/2}\ge1-\frac k2s_{3}\ge
1-\frac{k}{8(k-1)}\ge\frac{13}{16}$. Also $h(\theta_{3})=(1-s_{3})^{k-1}s_{3}\ge(1-(k-1)s_{3})s_{3}\ge\frac34s_{3}\ge
9/(64(k-1))$, so~$\frac{k\lambda^{2}}{2}h(\theta_{3})\ge9\lambda^{2}/128$. Take $\theta_{s}$ with $\cos^{k}\theta_{s}=\frac12$;
it lies in $[\theta_{3},\pi/2]$, since $\cos^{k}\theta_{3}\ge\frac{13}{16}$. Then~$u-\lambda/2\ge\lambda/4$, and the second
term in $E_{\rm far}$ is at least $\lambda^{2}/32$. With $x=2\log2/k\le0.47$, $\sin^{2}\theta_{s}=1-e^{-x}\ge
x(1-x/2)\ge0.76x$ and $\cos^{2k-2}\theta_{s}=2^{-2(k-1)/k}\ge\frac14$, so $\frac k2h(\theta_{s})\ge0.19\log2>\frac1{32}$.
Since also $\frac9{128}>\frac1{32}$, $E_{\rm far}\ge\lambda^{2}/32$.
\end{proof}

The threshold $\lambda_{0}(k,d)$ of Theorem~\ref*{thm:estim}{\rm(iii)} is defined in the main paper: it is the least
number such that the five conditions~(\ref*{eq:lam0a}) to~(\ref*{eq:lam0e}) of the main paper hold at every
$\lambda>\lambda_{0}(k,d)$. The constants $K_{N}$, $m_{N}$ and $\omega_{d-1}$ are those of Lemma~\ref{lem:far}. In the rest of
this appendix, (\ref*{eq:lam0a}) to~(\ref*{eq:lam0e}) are these five conditions, and $L=\log\lambda$. Condition~(\ref*{eq:lam0a}) fails at $\lambda=e$, so every number with this
property is at least $e$; the least one exists, since the set of these numbers contains its infimum. For $\lambda\ge e$
each condition, once true, stays true for larger~$\lambda$. Indeed~$\lambda^{2}/\log\lambda$ increases for
$\lambda>\sqrt e$. In~(\ref*{eq:lam0c}), $g_{1}$ increases, and the logarithm of~$K_{N}\lambda^{2}e^{-g_{1}^{2}/4}$ has
derivative $2/\lambda-\frac{87}{256}\rho g_{1}$, which is negative wherever~(\ref*{eq:lam0c}) holds, because there
$\frac{87}{128}\rho\lambda g_{1}\ge g_{1}^{2}\ge4\log(K_{N}\lambda^{2})>4$. The left side of the second inequality
of~(\ref*{eq:lam0d}) decreases once $g_{3}\ge0$. In~(\ref*{eq:lam0e}) the derivative of the logarithm of the left side is
at most~$(N+4)/\lambda-\lambda/16<0$ for~$\lambda>4\sqrt{d+3}$. Each condition holds for $\lambda$ large, the left side
of the second inequality of~(\ref*{eq:lam0d}) tending to $\frac54\cdot\frac{4k-5}{3k-4}\le\frac74$. So $\lambda_{0}(k,d)$
is finite for every odd~$k\ge3$ and $d\ge3$, and it is the largest of the five points beyond which the conditions
hold. The values of $\lambda_{0}$ given after the theorem in the main paper are computed numerically
from~(\ref*{eq:lam0}) of the main paper; at these pairs~(\ref*{eq:lam0e}) is the last condition to hold.

\begin{proof}[Proof of Theorem~\ref*{thm:estim}{\rm(iii)} of the main paper]
Let $\lambda>\lambda_{0}(k,d)$, so that $L>1$. We apply part {\rm(ii)} with $\theta_{1}=3\lambda^{-1}\sqrt{L/k}$ and
$u=\lambda-2\sqrt{2L}$. By~(\ref*{eq:lam0a}), $\sqrt{8L}\le\lambda/4$, so $u\ge3\lambda/4$, $\lambda/u\le4/3$ and
$L/\lambda\le\lambda/128$. Also $k\lambda^{2}\theta_{1}^{2}=9L$, $\theta_{1}\le3/\sqrt{128k}<1$, and
$\bar\varrho:=(k-\frac23)\theta_{1}^{2}=3(3k-2)L/(k\lambda^{2})\le9L/\lambda^{2}\le9/128$.

The hypotheses of {\rm(ii)} hold. By~(\ref*{eq:lam0b}), $u\ge3\lambda/4>N/\rho$. Since
$1-\cos^{k}\theta\le k(1-\cos\theta)\le k\theta^{2}/2$, $\lambda-\lambda\cos^{k}\theta_{1}\le9L/(2\lambda)<\sqrt{8L}$,
the last step because $\lambda\ge\sqrt{128L}$; so $u<\lambda\cos^{k}\theta_{1}$. Next
$4(k-1)\lambda\sin^{2}\theta_{1}\le4(k-1)\lambda\theta_{1}^{2}\le36L/\lambda\le36\lambda/128<u$, and, by
(\ref*{eq:lam0b}), $(k-1)\lambda\sin^{2}\theta_{1}+\sqrt{2N}/\rho\le9\lambda/128+9\lambda/16<3\lambda/4\le u$. These are
also the hypotheses of Lemma~\ref{lem:Aupper}. Finally $\frac34\rho u\ge\frac9{16}\rho\lambda\ge\sqrt{2N}$
by~(\ref*{eq:lam0b}), so Lemma~\ref{lem:far} applies.

The main term. By Lemma~\ref{lem:simple}, $\varrho\le\bar\varrho$, and with $-\log(1-y)\le y/(1-y)$ and
$1/(1-\bar\varrho)\le128/119<9/8$,
\[(1-\varrho)^{-\frac{d+1}{2}}\le\exp\Big(\frac{9}{16}(d+1)\bar\varrho\Big)=\exp\Big(\frac{27(3k-2)(d+1)L}{16k\lambda^{2}}\Big).\]
Next, $(k-1)\lambda\sin^{2}\theta_{1}/u\le9(k-1)L/(k\lambda u)\le12(k-1)L/(k\lambda^{2})$. Also
$\beta_{\star}\le9\rho L/\lambda\le9\rho\lambda/128$, so $\rho u-\beta_{\star}-\sqrt{2N}\ge g_{1}\ge0$, and
Lemma~\ref{lem:Aupper} and~(\ref*{eq:lam0c}) give $\varepsilon_{\det}\le K_{N}e^{-g_{1}^{2}/4}\le\lambda^{-2}\le u^{-2}$.
With $\ell_{1}>u$, Lemma~\ref{lem:simple} gives
$1+\epsilon_{A}\le\exp(12(k-1)L/(k\lambda^{2})+(N^{2}+2)/(2u^{2}))$, and $(N^{2}+2)/(2u^{2})\le\frac89(N^{2}+2)/\lambda^{2}$.
So the main term of {\rm(ii)} is at most $\frac{d-1}{k\lambda^{2}}\exp((c_{1}L+c_{1}')/\lambda^{2})$.

The far zone. In Lemma~\ref{lem:far}, $\ell_{3}\ge13\lambda/16$ and $\frac34\rho u-\sqrt{2N}\ge g_{3}\ge0$, so
$C_{M}\le(\frac54+K_{N}e^{-g_{3}^{2}/4})e^{128N(N+1)/(169\lambda^{2})}(4k-5)/(3k-4)\le2$ by~(\ref*{eq:lam0d}). Also
$E_{\rm far}\ge\lambda^{2}/32$ and $1+(k-1)\lambda/u\le1+4(k-1)/3=(4k-1)/3$, so $F_{\rm far}\le\lambda^{-4}/10$
by~(\ref*{eq:lam0e}). Since $k\lambda^{2}h(\theta_{1})=k\lambda^{2}\theta_{1}^{2}(1-\varrho)\ge9L(1-9/128)\ge8L$,
$\mathrm{Rem}_{\rm far}\le2\,\PP\{\chi^{2}_{N}>8L\}+\lambda^{-4}/10$.

The gap. $\Phi(u-\lambda)=\Phi(-\sqrt{8L})\le\varphi(\sqrt{8L})/\sqrt{8L}=\lambda^{-4}/(4\sqrt{\pi L})\le0.142\,\lambda^{-4}$,
as $L>1$. Since $\lambda^{-4}/10+0.142\,\lambda^{-4}\le\lambda^{-4}/4$, part {\rm(ii)} gives the bound of {\rm(iii)}.

The chi-square tail. Let $p=N/2-1\ge0$ and $y>N-2$. Integrating by parts once,
\[\int_{y}^{\infty}x^{p}e^{-x/2}\,\mathrm{d}x\ \le\ 2y^{p}e^{-y/2}+\frac{2p}{y}\int_{y}^{\infty}x^{p}e^{-x/2}\,\mathrm{d}x,\]
so $\PP\{\chi^{2}_{N}>y\}\le(y/2)^{p}e^{-y/2}/(\Gamma(N/2)(1-(N-2)/y))$. At $y=8L$ with $L\ge(d-3)/2$,
$(N-2)/y\le\frac14$. This gives
$2\,\PP\{\chi^{2}_{d-1}>8\log\lambda\}\le8(4\log\lambda)^{(d-3)/2}\lambda^{-4}/(3\,\Gamma(\frac{d-1}{2}))$, the bound
stated after the theorem in the main paper.
\end{proof}

Note that $\beta_{\star}$ is the supremum of $|\beta_{2}|$ over the \emph{near} zone $\theta\le\theta_{1}$, not over the
whole sphere, where it would be of order $\lambda$ and the bound would be useless. At the moderate values of $\lambda$
of Section~\ref*{sec:num} of the main paper, $\theta_{1}$ is optimised numerically in part {\rm(ii)}.

\section{Proofs for Sections~\ref*{sec:model} and~\ref*{sec:reverse} of the main paper}\label{app:proofs}

\begin{proof}[Proof of Proposition~\ref*{prop:sigmahat} of the main paper]
The count~(\ref*{eq:kappa}) and the law and the independence in~(\ref*{eq:sigmahat-law}) are
\citet[Assumption~(A5) and Proposition~6]{azais2024second}, applied to $X$ on $M=\Sph^{d-1}$, of dimension $d-1$,
whose Karhunen--Lo\`eve expansion~(\ref*{eq:rank-KL}) of the main paper has order $m$, with $Z=\sigma X$ under the
null. The integer $\kappa$ is also the codimension of the tangent space of the cone at $\lambda\psi_{t}$, $\lambda>0$.
For~(\ref*{eq:sigmahat-norm}), write
$D_{u}\psi_{t}$ for the derivative of $t\mapsto\psi_{t}$ in a tangent direction $u$. Differentiating
$\langle\psi_{s},\psi_{t}\rangle=\langle s,t\rangle^{k}$ once in $s$, then once more in $t$, at $s=t$, shows that
$\psi_{t}$ and the vectors $k^{-1/2}D_{u}\psi_{t}$, for $u$ in an orthonormal basis of the tangent space at $t$, are
orthonormal in $\mathbb R^{m}$. Also $\langle Y,D_{u}\psi_{t}\rangle=\langle u,\nabla Z(t)\rangle$,
as~$Z(t)=\langle Y,\psi_{t}\rangle$. Let $P_{t}^{\perp}$ be the orthogonal projection onto the orthogonal complement of
the span of $\psi_{t}$ and the $k^{-1/2}D_{u}\psi_{t}$. The numerator of~(\ref*{eq:helix}) of the main paper is
$\langle P_{t}^{\perp}Y,\psi_{s}\rangle$, and the~$\psi_{s}$ span $\mathbb R^{m}$. So, up to an orthogonal change of
coordinates, which leaves $\sum_{i}\bar\zeta_{i}^{2}$ unchanged, the $\bar\zeta_{i}$ are the coordinates of
$P_{t}^{\perp}Y$ in an orthonormal basis of its range, and~$\kappa\hat\sigma_{t}^{2}=\|P_{t}^{\perp}Y\|^{2}$, which
is~(\ref*{eq:sigmahat-norm}).
\end{proof}

\begin{proof}[Proof of Lemma~\ref*{lem:ratio} of the main paper]
Write $D(u)=\det(u\,\mathrm{Id}-r)=\prod_{i}(u-r_{i})$, the $r_{i}$ being the eigenvalues of $r$, and $h(u)=(\kappa+u^{2})^{-(\nu+1)/2}$ with $\nu=m-1$, so that
$H_{r}(\ell)=c\int_{\ell}^{\infty}Dh$ for a constant $c>0$ that cancels; $D$ is positive and increasing
on~$(\lmax(r),\infty)$, and the integrals converge because $N<\nu$.

\emph{The Student bounds.} Substituting $u=\vartheta w$, which maps $w=\ell_{2}$ to $u=\ell_{1}$,
\begin{subequations}
\begin{equation}
H_{r}(\ell_{1})=\vartheta\int_{\ell_{2}}^{\infty}D(\vartheta w)\,h(\vartheta w)\,\mathrm{d}w .
\end{equation}
For the density, $\kappa+\vartheta^{2}w^{2}\le\vartheta^{2}(\kappa+w^{2})$ because $\vartheta\ge1$, so
$h(\vartheta w)\ge\vartheta^{-(\nu+1)}h(w)$; and~$w\mapsto(\kappa+w^{2})/(\kappa+\vartheta^{2}w^{2})$ has
derivative of the sign of $\kappa(1-\vartheta^{2})\le0$, so for $w\ge\ell_{2}$ it is at most its value at
$\ell_{2}$ and $h(\vartheta w)\le\big((\kappa+\ell_{2}^{2})/(\kappa+\ell_{1}^{2})\big)^{(\nu+1)/2}h(w)$.
For the determinant, fix $i$ and~$w\ge\ell_{2}>\lmax(r)^{+}$. If $r_{i}\ge0$ then
$\vartheta w-r_{i}\ge\vartheta(w-r_{i})$; if $r_{i}<0$ then, with~$b=-r_{i}\le\lmin(r)^{-}$,
\begin{equation}
\frac{\vartheta w+b}{w+b}=\vartheta\,\frac{w+b/\vartheta}{w+b}\ \ge\ \vartheta\,\frac{w}{w+b}\ \ge\ \vartheta\,\frac{\ell_{2}}{\ell_{2}+\lmin(r)^{-}},
\end{equation}
the last step because $w\mapsto w/(w+b)$ increases. Hence
\begin{equation}
D(\vartheta w)\ge\vartheta^{N}\big(1+\lmin(r)^{-}/\ell_{2}\big)^{-N}D(w).
\end{equation}
In the other direction
\begin{equation}
\frac{\vartheta w-r_{i}}{w-r_{i}}=\vartheta+\frac{(\vartheta-1)r_{i}}{w-r_{i}}
\ \le\ \vartheta+\frac{(\vartheta-1)\lmax(r)^{+}}{\ell_{2}-\lmax(r)^{+}}
\ \le\ \vartheta\,\frac{\ell_{2}}{\ell_{2}-\lmax(r)^{+}},
\end{equation}
\end{subequations}
since $z\mapsto z/(w-z)$ increases and $w\ge\ell_{2}$, so
$D(\vartheta w)\le\vartheta^{N}\big(1-\lmax(r)^{+}/\ell_{2}\big)^{-N}D(w)$. Multiplying the two pairs of bounds and
using $\vartheta\cdot\vartheta^{N}\cdot\vartheta^{-(\nu+1)}=\vartheta^{N-\nu}=\vartheta^{-\kappa}$ gives~(\ref*{eq:dilH}) of the main paper.

\emph{The Gaussian bounds.} Substituting $u=\sqrt{w^{2}+2g}$, which also maps $w=\ell_{2}$ to $u=\ell_{1}$ and has
$\mathrm{d}u=(w/u)\,\mathrm{d}w$ and $\varphi(u)=e^{-g}\varphi(w)$,
\begin{subequations}
\begin{equation}
G_{r}(\ell_{1})=e^{-g}\int_{\ell_{2}}^{\infty}D(u)\,\frac wu\,\varphi(w)\,\mathrm{d}w .
\end{equation}
Here $u\ge w$, so $D(u)\ge D(w)$, and $w\mapsto w/\sqrt{w^{2}+2g}$ increases, so $w/u\ge\ell_{2}/\ell_{1}$; this is
the left inequality. For the right one, $w/u\le1$ and $u-w=2g/(u+w)\le g/\ell_{2}$, so
\begin{equation}
\frac{u-r_{i}}{w-r_{i}}=1+\frac{u-w}{w-r_{i}}\ \le\ 1+\frac{g}{\ell_{2}\big(\ell_{2}-\lmax(r)^{+}\big)} .
\end{equation}
\end{subequations}

\emph{The comparison.} In~(\ref*{eq:dilH}) the ratio of the two bounds is
$\vartheta^{N+1+\kappa}\big((\kappa+\ell_{2}^{2})/(\kappa+\ell_{1}^{2})\big)^{(\nu+1)/2}$ times two factors
that tend to $1$; since $N+1+\kappa=\nu+1$ and
$(\kappa+\ell_{2}^{2})/(\kappa+\ell_{1}^{2})\to\vartheta^{-2}$ as~$\ell_{2}\to\infty$, that ratio tends to
$\vartheta^{\nu+1}\vartheta^{-(\nu+1)}=1$. In~(\ref*{eq:dilG}) it is $(\ell_{1}/\ell_{2})$ times a factor tending to $1$
at fixed $g$.
\end{proof}

\begin{proof}[Proof of Lemma~\ref*{lem:layers} of the main paper]
The tensor $\operatorname{Sym}(t^{\otimes(k-j)}\otimes v^{\otimes j})$ is the mean of the $\binom kj$ tensors with $v$
at $j$ positions and $t$ elsewhere. So the expansion of $s^{\otimes k}$ reads
$\psi_{s}=\sum_{j}\binom kj^{1/2}\cos^{k-j}\theta\,\sin^{j}\theta\,\Psi_{j}(v)$. Two of these tensors, built with $v$ and
$v'$, have the inner product $\langle v,v'\rangle^{j}$ if the positions agree, and $0$ otherwise because $v$ and $v'$
are orthogonal to $t$. This gives the inner products and the orthogonality of the $\mathcal V_{j}$. Each $\psi_{s}$
lies in $\mathcal V_{0}\oplus\dots\oplus\mathcal V_{k}$ and the $\psi_{s}$ span $\mathbb R^{m}$
(Section~\ref*{sub:rank} of the main paper), so this sum is $\mathbb R^{m}$. As $\Psi_{0}=\psi_{t}$ and
$\Psi_{1}(u)=k^{-1/2}D_{u}\psi_{t}$, the range of $P_{t}^{\perp}$ is $\mathcal V_{2}\oplus\dots\oplus\mathcal V_{k}$. The
numerator of~(\ref*{eq:helix}) of the main paper is $\langle P_{t}^{\perp}Y,\psi_{s}\rangle$, which gives the first
formula. Write $Y[x_{1},\dots,x_{k}]=\langle Y,x_{1}\otimes\dots\otimes x_{k}\rangle$. Two derivatives of
$\theta\mapsto Z(\cos\theta\,t+\sin\theta\,v)$ at $0$ give $v^{\top}\nabla^{2}Z(t)v=-kZ(t)+k(k-1)Y[t,\dots,t,v,v]$, and
$Y[t,\dots,t,v,v]=\binom k2^{-1/2}\langle Y,\Psi_{2}(v)\rangle$ as $Y$ is symmetric. So
$v^{\top}\Omega_{t}v=(k-1)\binom k2^{-1/2}\langle Y,\Psi_{2}(v)\rangle$, and $(k-1)\binom k2^{-1/2}=1/\rho$. Since
$\langle\Psi_{2}(v),\Psi_{2}(v')\rangle=\langle v,v'\rangle^{2}=\langle vv^{\top},v'v'^{\top}\rangle_{F}$ and the
$vv^{\top}$ span the symmetric matrices, the map $vv^{\top}\mapsto\Psi_{2}(v)$ extends to the isometry $\iota$. Then
$P_{\mathcal V_{2}}Y=\iota(\rho\,\Omega_{t})$, as both have the inner product $\rho\,v^{\top}\Omega_{t}v$ with
every~$\Psi_{2}(v)$, and the last formula follows from $\kappa\hat\sigma_{t}^{2}=\|P_{t}^{\perp}Y\|^{2}$.
\end{proof}

\begin{proof}[Proof of Lemma~\ref*{lem:mono} of the main paper]
Write $D(u)=\prod_{i}(u-r_{i})$ and $h(u)=(\kappa+u^{2})^{-(\nu+1)/2}$, as in the proof of Lemma~\ref*{lem:ratio} above, so
that $H_{r}(\ell)=c\int_{\ell}^{\infty}Dh$. As $D\ge0$ on~$[\ell_{2},\infty)$, $H_{r}(\ell_{2})$ is nonincreasing in
$\ell_{2}$, and the ratio is nondecreasing in $\ell_{2}$. Fix $i$ with~$r_{i}<\ell_{2}$ and let $U$ have the density
proportional to $Dh$ on $(\ell_{2},\infty)$. Differentiating under the integral, which is licit since
$|\partial_{r_{i}}(Dh)|\le Dh/(\ell_{2}-r_{i})$ there,
\[
\partial_{r_{i}}\log\frac{H_{r}(\ell_{1})}{H_{r}(\ell_{2})}=\E\Big[\frac1{U-r_{i}}\Big]-\E\Big[\frac1{U-r_{i}}\ \Big|\ U>\ell_{1}\Big].
\]
The function $u\mapsto1/(u-r_{i})$ is positive and decreasing on $(\ell_{2},\infty)$, so its conditional mean on
$\{U>\ell_{1}\}$ is at most $1/(\ell_{1}-r_{i})$, which is at most its conditional mean on $\{U\le\ell_{1}\}$. The
mean $\E[1/(U-r_{i})]$ lies between the two, so the derivative is nonnegative. The case $r_{i}=\ell_{2}$ follows by
continuity of $H_{r}$ in $r$.
\end{proof}

\begin{proof}[Proof of the matching bound in Theorem~\ref*{thm:exponent} of the main paper]
With $\sigma=1$, let $E=\rho^{2}/(2(A^{2}+\rho^{2}))$ with $A$ the constant~(\ref*{eq:Adef}) of the main paper.
We show that the~$\limsup$ of $-\lambda^{-2}\log(1-\PP_{\lambda}\{\mathbf T\le\alpha\})$ is at most $E$.

\emph{A configuration.} Fix a positive semidefinite $r$ of size $N$ with $\|r\|_{F}=T_{\max}$, and
$a\in\big(\lmax(r),u_{\alpha}(r,\lmax(r))\big)$. With $t=e$ in Lemma~\ref*{lem:layers} of the main paper, let
$y_{0}=a\,\psi_{e}+\iota(\rho r)$. By the same lemma, its field at $s=\cos\theta\,e+\sin\theta\,v$ is
$z(s)=a\cos^{k}\theta+(1-\cos^{k}\theta)\,\rho\,w_{2}(\theta)\,v^{\top}rv$. By its definition,
$\rho w_{2}=\tfrac k2\cos^{k-2}\theta\sin^{2}\theta/(1-\cos^{k}\theta)$. For $x\in[0,1)$,
$f(x)=1-x^{k}-\tfrac k2x^{k-2}(1-x^{2})$ has $f(1)=0$ and $f'(x)=\tfrac k2(k-2)x^{k-3}(x^{2}-1)\le0$, so $f\ge0$: with
$x=\cos\theta$, $\rho w_{2}\le1$ on $(0,\pi/2]$. On $[\pi/2,\pi)$, $w_{2}\le0$ as $k$ is odd. As $r$ is positive
semidefinite, $\rho\,w_{2}\,v^{\top}rv\le\lmax(r)<a$. So $z(s)<a=z(e)$ for $s\ne e$, the value at $-e$ being $-a$. At
$y_{0}$ the maximiser is $e$ and $\lambda_{1}=a$; $\nabla z(e)=0$, as $y_{0}$ has no component in $\mathcal V_{1}$;
$\Omega=r$, so $\nabla^{2}z(e)=k(r-a\,\mathrm{Id})$ is negative definite;
$\kappa\hat\sigma^{2}=\|\iota(\rho r)\|^{2}=\rho^{2}T_{\max}^{2}=\kappa$; and
$\lambda_{2}=\sup\rho\,w_{2}\,v^{\top}rv=\lmax(r)$, the value on the fibre $\theta=0$. Hence
$\mathbf T(y_{0})=H_{r}(a)/H_{r}(\lmax(r))>\alpha$.

\emph{Continuity at $y_{0}$.} The field $z_{y}(s)=\langle y,\psi_{s}\rangle$ is linear in $y$, so $z_{y}\to z_{y_{0}}$
in $C^{2}$ as $y\to y_{0}$. Take a geodesic ball $B$ around $e$ on which the Hessian of $z_{y_{0}}$ is negative definite. As $e$ is the only
maximiser of $z_{y_{0}}$, there is $\varsigma>0$ with $z_{y_{0}}\le a-\varsigma$ outside $B$. For $y$ close to $y_{0}$,
$z_{y}<a-\varsigma/2<z_{y}(e)$ outside $B$ and the Hessian of $z_{y}$ is negative definite on $B$. So $z_{y}$ has one
maximiser $\hat t(y)$, and it lies in $B$. As $B$ can be taken small, $\hat t(y)\to e$. Then $\lambda_{1}$,
$\hat\sigma^{2}=(\|y\|^{2}-\lambda_{1}^{2})/\kappa$ and the eigenvalues of $\Omega$ are continuous at $y_{0}$. By
Lemma~\ref{lem:lam2cont}, the second maximum $\lambda_{2}^{t}$ of the field $z_{y}$ at a point $t$ is continuous in
$(t,y)$. As $\lambda_{2}(y)$ is its value at $t=\hat t(y)$ and $\hat t(y)\to e$, $\lambda_{2}$ is continuous at
$y_{0}$. The map
$(\ell_{1},\ell_{2},r)\mapsto H_{r}(\ell_{1})/H_{r}(\ell_{2})$ is continuous on $\{\lmax(r)\le\ell_{2}\le\ell_{1}\}$,
where $H_{r}(\ell_{2})>0$, and $(T_{1},T_{2},\Omega/\hat\sigma)$ lies in this set by~(\ref*{eq:order}) of the main
paper. So there is $\eta_{0}>0$ with $\mathbf T>\alpha$ on the ball $\mathbb B(y_{0},\eta_{0})$ of $\mathbb R^{m}$.

\emph{The upper bound.} $\mathbf T$ is invariant under $y\mapsto cy$, $c>0$ (Section~\ref*{sub:closed} of the main
paper), so $\mathbf T>\alpha$ on the open cone $\mathcal U=\{cy:c>0,\ y\in\mathbb B(y_{0},\eta_{0})\}$. Let
$y^{*}=c^{*}y_{0}$ with $c^{*}=\langle y_{0},\psi_{e}\rangle/\|y_{0}\|^{2}=a/(a^{2}+\kappa)$, the point of the ray closest
to $\psi_{e}$, and let $\eta>0$ with $\mathbb B(y^{*},\eta)\subset\mathcal U$. For $\lambda\ge1$,
$\mathbb B(\lambda y^{*},\eta)\subset\lambda\,\mathbb B(y^{*},\eta)\subset\mathcal U$. As $Y=\lambda\psi_{e}+W$ with
$W$ standard Gaussian in $\mathbb R^{m}$,
\[
1-\PP_{\lambda}\{\mathbf T\le\alpha\}\ \ge\ \PP\big\{W\in\mathbb B(\lambda(y^{*}-\psi_{e}),\eta)\big\}
\ \ge\ \frac{|\mathbb B(0,\eta)|}{(2\pi)^{m/2}}\exp\Big(-\tfrac12\big(\lambda\|y^{*}-\psi_{e}\|+\eta\big)^{2}\Big),
\]
with $|\mathbb B(0,\eta)|$ the volume of the ball. Here
$\|y^{*}-\psi_{e}\|^{2}=1-\langle y_{0},\psi_{e}\rangle^{2}/\|y_{0}\|^{2}=\kappa/(\kappa+a^{2})=\rho^{2}/(\rho^{2}+(a/T_{\max})^{2})$.
So the $\limsup$ is at most $\rho^{2}/(2(\rho^{2}+(a/T_{\max})^{2}))$. Let $a$ tend to~$u_{\alpha}(r,\lmax(r))$ and take
$r$ at the maximum in~(\ref*{eq:Adef}) of the main paper: the $\limsup$ is at most~$E$.
\end{proof}

The proof of Theorem~\ref*{thm:sqrtd} of the main paper bounds the supremum of the noise field, then the second
maximum and the Hessian through it. With $\sigma=1$, write $\mathcal X=\sup_{M}X$, $\xi=X(e)$
and~$\hat c=\langle e,\hat t\rangle$. As $X(-t)=-X(t)$ for all $t$ (Lemma~\ref*{lem:pos} of the main paper),
$\mathcal X=\sup_{M}|X|\ge0$.

\begin{lemma}[Supremum of the noise field]\label{lem:supX}
$\E\mathcal X\le\sqrt{kd}$, and $\PP\{\mathcal X>\sqrt{kd}+\tau\}\le e^{-\tau^{2}/2}$ for $\tau\ge0$.
\end{lemma}

\begin{proof}
For $s,t\in M$, $\E(X(s)-X(t))^{2}=2(1-\langle s,t\rangle^{k})$. For $x\in[-1,1)$, $1-x^{k}=(1-x)\sum_{i<k}x^{i}$, and the
sum lies in $[0,k]$, since it equals $(1-x^{k})/(1-x)\ge0$ and each term is at most $1$ in absolute value. So
$\E(X(s)-X(t))^{2}\le2k(1-\langle s,t\rangle)=k\|s-t\|^{2}=\E(X'(s)-X'(t))^{2}$, with $X'(t)=\sqrt k\langle g,t\rangle$ and
$g$ standard Gaussian in $\mathbb R^{d}$. The Sudakov--Fernique inequality \citep[Chapter~2]{azais2009level} gives
$\E\mathcal X\le\E\sup_{M}X'=\sqrt k\,\E\|g\|\le\sqrt{kd}$. The Borell--TIS inequality \citep[Chapter~2]{azais2009level},
with $\sup_{t}\Var X(t)=1$, gives $\PP\{\mathcal X>\E\mathcal X+\tau\}\le e^{-\tau^{2}/2}$.
\end{proof}

\begin{lemma}[Pathwise bounds through the supremum]\label{lem:banach}
For $\theta\in(0,\pi)$ let
\[g_{k}(\theta)=\frac{(|\cos\theta|+\sin\theta)^{k}-|\cos\theta|^{k}-k|\cos\theta|^{k-1}\sin\theta}{1-\cos^{k}\theta},\]
so that $\Upsilon_{k}=\sup g_{k}$ in Theorem~\ref*{thm:sqrtd} of the main paper. Then $k-1<\Upsilon_{k}<\infty$. Almost
surely, for every $t\in M$, every unit $v\perp t$, every $2\le j\le k$ and $s=\cos\theta\,t+\sin\theta\,v\in M_{t}$, with
$\Psi_{j}(v)$ of Lemma~\ref*{lem:layers} of the main paper,
{\rm(a)} $|\langle W,\Psi_{j}(v)\rangle|\le\binom kj^{1/2}\mathcal X$, and
{\rm(b)} $|\zeta^{t}(s)|\le g_{k}(\theta)\,\mathcal X\le\Upsilon_{k}\mathcal X$.
\end{lemma}

\begin{proof}
The form $S_{W}(x_{1},\dots,x_{k})=\langle W,x_{1}\otimes\dots\otimes x_{k}\rangle$ is symmetric and $k$-linear
on~$\mathbb R^{d}$, and $S_{W}(x,\dots,x)=X(x)$. Banach's theorem on symmetric multilinear forms
\citep{banach1938homogene}, stated for tensors in \citet[Section~5]{friedland2018nuclear}, gives
\[\sup\big\{|S_{W}(x_{1},\dots,x_{k})|:\ \|x_{1}\|=\dots=\|x_{k}\|=1\big\}=\sup_{M}|X|=\mathcal X.\]
As $W$ is symmetric, $\langle W,\Psi_{j}(v)\rangle$ is $\binom kj^{1/2}$ times the value of $S_{W}$ at $k-j$ copies
of $t$ and $j$ copies of $v$, which gives (a). Lemma~\ref*{lem:layers} of the main paper with $W$ in place of $Y$ gives
$\zeta^{t}(s)=\sum_{j\ge2}w_{j}(\theta)\langle W,\Psi_{j}(v)\rangle$, so by (a) and the binomial formula,
$|\zeta^{t}(s)|\le\mathcal X\sum_{j\ge2}\binom kj|\cos\theta|^{k-j}\sin^{j}\theta/(1-\cos^{k}\theta)=g_{k}(\theta)\mathcal X$.
Finally $g_{k}(\theta)\to0$ as $\theta\to\pi$, and $g_{k}(\theta)=k-1+\tfrac2k\binom k3\theta+O(\theta^{2})$ as
$\theta\to0$, so $g_{k}$ is bounded and $\Upsilon_{k}>k-1$.
\end{proof}

\begin{lemma}[The argmax under the alternative]\label{lem:argmaxalt}
On the event $\{\lambda>\mathcal X-\xi\}$, $\hat c>0$, $1-\hat c^{2}\le(\mathcal X-\xi)/\lambda$, and
\[\lambda_{2}\ \le\ \Upsilon_{k}\mathcal X+(k-1)(\mathcal X-\xi),\qquad \lmax(\Omega)\ \le\ (k-1)(2\mathcal X-\xi).\]
\end{lemma}

\begin{proof}
$\lambda_{1}=Z(\hat t)\ge Z(e)=\lambda+\xi$ and $Z(\hat t)=\lambda\hat c^{k}+X(\hat t)\le\lambda\hat c^{k}+\mathcal X$, so
$\lambda\hat c^{k}\ge\lambda+\xi-\mathcal X>0$. Hence $\hat c>0$, $\hat c^{2}\ge\hat c^{k}$ and
$1-\hat c^{2}\le1-\hat c^{k}\le(\mathcal X-\xi)/\lambda$. On $M_{\hat t}$, $Z^{|\hat t}=\zeta^{\hat t}+\mu_{\hat t}$,
with~$\zeta^{\hat t}\le\Upsilon_{k}\mathcal X$ by Lemma~\ref{lem:banach}(b) and
$|\mu_{\hat t}|\le(k-1)\lambda(1-\hat c^{2})\le(k-1)(\mathcal X-\xi)$ by~\eqref{eq:driftbound};
and~$Z^{|\hat t}(-\hat t)=0$. This bounds $\lambda_{2}$. For a unit $v\perp\hat t$, Lemma~\ref*{lem:layers} of the main
paper and Lemma~\ref{lem:banach}(a) give
\[\rho\,v^{\top}\Omega v=\langle W,\Psi_{2}(v)\rangle+\lambda\langle\psi_{e},\Psi_{2}(v)\rangle
\le\tbinom k2^{1/2}\big(\mathcal X+\lambda\hat c^{k-2}\langle e,v\rangle^{2}\big)\le\tbinom k2^{1/2}\big(\mathcal X+\lambda(1-\hat c^{2})\big),\]
since $0<\hat c\le1$ and $\langle e,v\rangle^{2}\le1-\hat c^{2}$. As $\binom k2^{1/2}/\rho=k-1$, this gives the bound on
$\lmax(\Omega)$.
\end{proof}

\begin{proof}[Proof of Theorem~\ref*{thm:sqrtd} of the main paper]
Let $\chi^{2}=\|W-\xi\psi_{e}\|^{2}$, independent of $\xi$ with the law $\chi^{2}_{m-1}$, and
$\mathcal G=\{\mathcal X\le\bar{\mathcal X}\}\cap\{|\xi|\le\tau\}\cap\{\chi^{2}\le m-1+\delta\}$. By Lemma~\ref{lem:supX},
$\PP(\mathcal G^{c})$ is at most the bound of the theorem. We show that $\mathbf T\le\alpha$ almost surely on
$\mathcal G$.

On $\mathcal G$, $\lambda\ge\bar v+\tau>\bar{\mathcal X}+\tau\ge\mathcal X-\xi$, so Lemma~\ref{lem:argmaxalt} applies:
$\lambda_{2}\le\Upsilon_{k}\mathcal X+(k-1)(\mathcal X-\xi)\le\bar v$ and
$\lmax(\Omega)^{+}\le(k-1)(2\mathcal X-\xi)\le\bar\ell$, using $2\mathcal X-\xi\ge\mathcal X\ge0$. Also $\bar\ell<\bar v$,
because $\bar v-\bar\ell=(\Upsilon_{k}-k+1)\bar{\mathcal X}$ and $\Upsilon_{k}>k-1$, and
$\lambda_{1}\ge\lambda+\xi\ge\lambda-\tau\ge\bar v$.

Put $r=\Omega/\hat\sigma$ and $v=\bar v/\hat\sigma$. Then $\lmax(r)^{+}<v\le T_{1}$, and $\lmax(r)\le T_{2}\le v$
by~(\ref*{eq:order}) of the main paper. As $H_{r}$ decreases on $[\lmax(r),\infty)$, $\mathbf T\le H_{r}(T_{1})/H_{r}(v)$.
The right-hand inequality of~(\ref*{eq:dilH}) of the main paper, with $\ell_{2}=v$, $\ell_{1}=T_{1}$,
$\lmax(r)^{+}/v\le\bar\ell/\bar v$ and $\nu+1=m$, gives
\[\mathbf T\ \le\ \Big(\frac{\lambda_{1}}{\bar v}\Big)^{N+1}\Big(1-\frac{\bar\ell}{\bar v}\Big)^{-N}
\Big(\frac{\kappa\hat\sigma^{2}+\bar v^{2}}{\kappa\hat\sigma^{2}+\lambda_{1}^{2}}\Big)^{m/2}.\]
By~(\ref*{eq:sigmahat-norm}) of the main paper at $\hat t$, $\kappa\hat\sigma^{2}=\|Y\|^{2}-\lambda_{1}^{2}$. Hence
$\log\mathbf T\le F(\lambda_{1})+N\log(\bar v/(\bar v-\bar\ell))$, with
$F(x)=d\log(x/\bar v)+\tfrac m2\log\big(1-(x^{2}-\bar v^{2})/\|Y\|^{2}\big)$ on $[\bar v,\|Y\|)$.

On $\mathcal G$, $\|Y\|^{2}=(\lambda+\xi)^{2}+\chi^{2}\le\bar n^{2}$. For $x\ge\lambda-\tau$,
$F'(x)=d/x-mx/(\|Y\|^{2}-x^{2}+\bar v^{2})\le d/x-mx/\|Y\|^{2}\le0$, by the second hypothesis. As
$\lambda-\tau\le\lambda_{1}<\|Y\|$ (proof of Corollary~\ref*{cor:closed} of the main paper),
$F(\lambda_{1})\le F(\lambda-\tau)$. Finally $F(\lambda-\tau)$ is nondecreasing in~$\|Y\|^{2}$,
since $(\lambda-\tau)^{2}\ge\bar v^{2}$, so it is at most its value with $\bar n$ in place of $\|Y\|$,
and~(\ref*{eq:sqrtd-H}) of the main paper gives~$\log\mathbf T\le\log\alpha$.
\end{proof}

\begin{proof}[Proof of the last part of Theorem~\ref*{thm:sqrtd} of the main paper]
Put $\tau=d^{1/4}$ and $\delta=m^{3/4}$ in the theorem. By the Chebyshev inequality the last term of its bound is at
most $2(m-1)/\delta^{2}$, so the bound tends to $0$. As $d\to\infty$, $\bar v/\sqrt d\to v_{k}$,
$\bar\ell/\sqrt d\to\ell_{k}$, $(\lambda_{d}-\tau)/\sqrt d\to\bar\lambda$ and~$\bar n^{2}/m\to1$, because
$\lambda_{d}^{2}=O(d)$ and $m\ge d^{k}/k!$. The first hypothesis holds for $d$ large, as~$\bar\lambda>C_{k}>v_{k}$, and
the second because $\bar\lambda>v_{k}>1$. With $\log(1-y)\le-y$, the left-hand side of~(\ref*{eq:sqrtd-H}) of the main
paper is at most
\[d\Big[\log\frac{\bar\lambda}{v_{k}}+\log\frac{v_{k}}{v_{k}-\ell_{k}}-\frac12\big(\bar\lambda^{2}-v_{k}^{2}\big)+o(1)\Big].\]
The map $x\mapsto\tfrac12(x^{2}-v_{k}^{2})-\log(x/v_{k})$ increases on $(v_{k},\infty)$, as $v_{k}>1$, and
equals $\log(v_{k}/(v_{k}-\ell_{k}))>0$ at $x=C_{k}$. So the bracket has a negative limit when $\bar\lambda>C_{k}$,
and~(\ref*{eq:sqrtd-H}) holds for $d$ large.
\end{proof}

\begin{lemma}[The constant in growing dimension]\label{lem:Alim}
Let $k$ and $\alpha$ be fixed. The constant~$A$ of~(\ref*{eq:Adef}) of the main paper tends to $1$ as $d\to\infty$.
\end{lemma}

\begin{proof}
$A>1$ (main paper, after~(\ref*{eq:Adef})). Let $r$ be positive semidefinite with $\|r\|_{F}=T_{\max}$, and
$u=u_{\alpha}(r,\lmax(r))$. Every eigenvalue of $r$ is at most $T_{\max}$. If $u>T_{\max}$, Lemma~\ref*{lem:mono} of the
main paper, applied to $\ell_{2}$ raised from $\lmax(r)$ to $T_{\max}$, then to every eigenvalue raised to $T_{\max}$,
gives $\alpha=H_{r}(u)/H_{r}(\lmax(r))\le P(u)$, with
$P(u)=H_{T_{\max}\mathrm{Id}}(u)/H_{T_{\max}\mathrm{Id}}(T_{\max})$. With $u=T_{\max}(1+x)$ and
$\kappa=\rho^{2}T_{\max}^{2}$, $P(T_{\max}(1+x))=\int_{x}^{\infty}g/\int_{0}^{\infty}g$ with
$g(x)=x^{N}(\rho^{2}+(1+x)^{2})^{-m/2}$, and $P(T_{\max}(1+x))$ decreases in $x$. So it suffices to show, for every
$\epsilon\in(0,1]$, that $P(T_{\max}(1+\epsilon))<\alpha$ for $d$ large: then $u\le T_{\max}(1+\epsilon)$ for every
such $r$, and $A\le1+\epsilon$. As $\rho<1$, $\rho^{2}+(1+x)^{2}\le2(1+x)^{2}$ and
$(\log g)'(x)\le N/x-m/(2(1+x))=\big(N-x(m/2-N)\big)/(x(1+x))$. Let $p=(m/2-N)/2$, which tends to infinity with $d$,
and take $d$ so large that $p>1$ and $\epsilon p\ge4N$. Then $(\log g)'\le-p/(1+x)$ on $[\epsilon/4,\infty)$, so $g$
decreases there and $g(x)\le g(\epsilon/2)\big((1+\epsilon/2)/(1+x)\big)^{p}$ for $x\ge\epsilon/2$. Hence
\[P\big(T_{\max}(1+\epsilon)\big)\ \le\ \frac{\int_{\epsilon}^{\infty}g}{\int_{\epsilon/4}^{\epsilon/2}g}
\ \le\ \frac{4}{\epsilon}\,\frac{1+\epsilon}{p-1}\Big(\frac{1+\epsilon/2}{1+\epsilon}\Big)^{p},\]
which tends to $0$ as $d\to\infty$, since $N/m\to0$.
\end{proof}

The last lemma of this appendix gives the growth, in the order $k$, of the constant $C_{k}$ of
Theorem~\ref*{thm:sqrtd} of the main paper.

\begin{lemma}[Growth of the constant $C_{k}$]\label{lem:Ck}
Let $k\ge3$ be odd, and let $\Upsilon_{k}$, $v_{k}$, $\ell_{k}$ and $C_{k}$ be as in Theorem~\ref*{thm:sqrtd} of the main
paper. Then
\begin{subequations}
\begin{gather}
2^{k/2}\big(1-(k+1)2^{-k}\big)\ \le\ \Upsilon_{k}\ \le\
\max\Big\{\frac{2^{k/2}}{1-\cos^{k}(1/4)},\ k\Big(1+\frac{k}{32}\Big)e^{k/4}\Big\},\label{eq:Ck-ups}\\
v_{k}\ <\ C_{k}\ \le\ v_{k}+\frac{v_{k}}{v_{k}^{2}-1}\,\log\frac{v_{k}}{v_{k}-\ell_{k}} .\label{eq:Ck-bracket}
\end{gather}
In particular $C_{k}/(2^{k/2}\sqrt k)\to1$ as $k\to\infty$ along the odd integers. Moreover
\begin{equation}\label{eq:Ck-bound}
C_{k}\le1.72\cdot2^{k/2}\sqrt k\quad\text{for }k\ge3,\qquad C_{k}\le1.02\cdot2^{k/2}\sqrt k\quad\text{for }k\ge21 .
\end{equation}
\end{subequations}
\end{lemma}

\begin{proof}
Let $g_{k}$ be the function of Lemma~\ref{lem:banach}, so that $\Upsilon_{k}=\sup_{0<\theta<\pi}g_{k}(\theta)$. Write
$a=|\cos\theta|$ and $b=\sin\theta$. The numerator of $g_{k}(\theta)$ is
$n_{k}(\theta)=(a+b)^{k}-a^{k}-ka^{k-1}b=\sum_{j=2}^{k}\binom kj a^{k-j}b^{j}$. It is non-negative, and it is at most
$(a+b)^{k}\le2^{k/2}$. At $\theta=\pi/4$, $a=b=2^{-1/2}$ and
$g_{k}(\pi/4)=\big(2^{k/2}-(k+1)2^{-k/2}\big)/\big(1-2^{-k/2}\big)\ge2^{k/2}-(k+1)2^{-k/2}$, because the
numerator is positive, as $2^{k}>k+1$, and the denominator is in $(0,1)$. This is the lower bound of~\eqref{eq:Ck-ups}.

For the upper bound we split $(0,\pi)$ into three ranges. If $\pi/2\le\theta<\pi$, then $\cos^{k}\theta\le0$, as $k$ is odd, so
the denominator $1-\cos^{k}\theta$ is at least $1$ and $g_{k}(\theta)\le2^{k/2}$. If $1/4\le\theta\le\pi/2$,
then~$0\le\cos^{k}\theta\le\cos^{k}(1/4)$ and $g_{k}(\theta)\le2^{k/2}/(1-\cos^{k}(1/4))$. Let now $0<\theta<1/4$. As~$a\le1$,
$\binom kj\le k^{j}/j!$ and $b\le\theta$,
\[n_{k}(\theta)\ \le\ \sum_{j\ge2}\frac{(k\theta)^{j}}{j!}\ \le\ \frac{(k\theta)^{2}}{2}\,e^{k\theta}.\]
The second derivative of $\log\cos$ is $-1/\cos^{2}\le-1$ on $[0,\pi/2)$, so $\cos\theta\le e^{-\theta^{2}/2}$. With
$y=k\theta^{2}/2$ and $1-e^{-y}\ge y/(1+y)$, the denominator $1-\cos^{k}\theta\ge1-e^{-y}$ is at least $y/(1+y)$. Hence
\[g_{k}(\theta)\ \le\ \frac{(k\theta)^{2}}{2}\,e^{k\theta}\,\frac{1+y}{y}\ =\ k\,e^{k\theta}\Big(1+\frac{k\theta^{2}}{2}\Big)
\ \le\ k\Big(1+\frac{k}{32}\Big)e^{k/4}.\]
The three bounds give the upper bound of~\eqref{eq:Ck-ups}, as $2^{k/2}\le2^{k/2}/(1-\cos^{k}(1/4))$.

By Lemma~\ref{lem:banach}, $\Upsilon_{k}>k-1$, so $v_{k}>\ell_{k}$. Also $\ell_{k}\ge4\sqrt3>1$, as $k\ge3$. So
$L_{k}:=\log\big(v_{k}/(v_{k}-\ell_{k})\big)>0$. Let $F_{k}(x)=\tfrac12(x^{2}-v_{k}^{2})-\log(x/v_{k})$.
By~(\ref*{eq:sqrtd-C}) of the main paper, $C_{k}>v_{k}$ and $F_{k}(C_{k})=L_{k}$. For $x>v_{k}$, $\log(x/v_{k})\le(x-v_{k})/v_{k}$, so
\[F_{k}(x)\ \ge\ (x-v_{k})\Big(\frac{x+v_{k}}{2}-\frac1{v_{k}}\Big)\ \ge\ (x-v_{k})\Big(v_{k}-\frac1{v_{k}}\Big).\]
At $x=C_{k}$ this gives $L_{k}\ge(C_{k}-v_{k})(v_{k}-1/v_{k})$, which is~\eqref{eq:Ck-bracket}, as $v_{k}>1$.

Let $k\to\infty$ along the odd integers. The lower bound of~\eqref{eq:Ck-ups}, divided by $2^{k/2}$, is
$1-(k+1)2^{-k}$, and it tends to $1$. The upper bound, divided by $2^{k/2}$, is the maximum of
$1/(1-\cos^{k}(1/4))$, which tends to $1$ as $\cos(1/4)<1$, and of $k(1+k/32)(e^{1/4}/\sqrt2)^{k}$, which tends to
$0$ as $e^{1/4}<\sqrt2$, that is $e<4$. So $\Upsilon_{k}/2^{k/2}\to1$. As $(k-1)2^{-k/2}\to0$,
$v_{k}/(2^{k/2}\sqrt k)=(\Upsilon_{k}+k-1)/2^{k/2}\to1$, and $v_{k}\to\infty$. Also
$\ell_{k}/v_{k}=2(k-1)/(\Upsilon_{k}+k-1)\to0$, so $L_{k}\to0$. By~\eqref{eq:Ck-bracket}, $C_{k}-v_{k}\to0$, and
the claim follows.

It remains to prove the two bounds on $C_{k}/(2^{k/2}\sqrt k)$. Let first $k\ge125$ be odd. Then
$\cos^{k}(1/4)\le\cos^{125}(1/4)<0.0194$. The function $k\mapsto k(1+k/32)(e^{1/4}/\sqrt2)^{k}$ decreases
on~$[125,\infty)$, as the derivative of its logarithm is $1/k+1/(32+k)+\tfrac14-\tfrac12\log2<0$ there, and it is below
$0.01$ at $k=125$. So~\eqref{eq:Ck-ups} gives $\Upsilon_{k}\le2^{k/2}/0.9806$, and
$v_{k}/(2^{k/2}\sqrt k)=\Upsilon_{k}/2^{k/2}+(k-1)2^{-k/2}<1.0198$, as $(k-1)2^{-k/2}<10^{-16}$. By the lower bound
of~\eqref{eq:Ck-ups}, $\Upsilon_{k}\ge3(k-1)$, so $\ell_{k}/v_{k}\le1/2$ and $L_{k}\le\log2$. As
$v_{k}\ge\ell_{k}\ge248\sqrt{125}$, \eqref{eq:Ck-bracket} gives $C_{k}\le v_{k}\big(1+\log2/(v_{k}^{2}-1)\big)\le(1+10^{-7})\,v_{k}$. Hence
$C_{k}\le1.02\cdot2^{k/2}\sqrt k$.

For the odd $k$ from $3$ to $123$ the bounds are checked by interval arithmetic. On $(0,\pi/2]$,
$1-\cos\theta=\sin^{2}\theta/(1+\cos\theta)$ gives
\[g_{k}(\theta)=(1+a)\,\frac{\sum_{j=2}^{k}\binom kj a^{k-j}b^{j-2}}{\sum_{i=0}^{k-1}a^{i}},\qquad a=\cos\theta,\ b=\sin\theta,\]
a quotient of two polynomials with non-negative coefficients. Also $g_{k}(\pi-\theta)\le g_{k}(\theta)$
on~$(0,\pi/2]$: the numerator is the same and the denominator $1+\cos^{k}\theta$ is larger, as $k$ is odd. So
$\Upsilon_{k}$ is the supremum of $g_{k}$ over $(0,\pi/2]$. On a cell $[\theta_{0},\theta_{1}]\subset[0,\pi/2]$,
$\cos\theta_{1}\le a\le\cos\theta_{0}$ and $b\le\sin\theta_{1}$. So $g_{k}$ is at most the value of the quotient
with $a=\cos\theta_{0}$ and $b=\sin\theta_{1}$ in the factor $1+a$ and in the numerator, and $a=\cos\theta_{1}$ in
the denominator. The largest of these values over a subdivision of $[0,\pi/2]$ is an upper bound $\Upsilon^{+}$ of
$\Upsilon_{k}$, and each value of $g_{k}$ is a lower bound $\Upsilon^{-}$. For $v>\ell_{k}$ let $C(v)>v$
solve~(\ref*{eq:sqrtd-C}) of the main paper with $v$ in place of $v_{k}$, so that $C_{k}=C(v_{k})$. The equation is
$\tfrac12(C^{2}-v^{2})-\log C+\log(v-\ell_{k})=0$. Its left side increases in $C>1$, and it decreases in $v$ when
$v(v-\ell_{k})>1$. So $C(v)$ increases on the half-line of the $v>\ell_{k}$ with $v(v-\ell_{k})>1$. The proof
of~\eqref{eq:Ck-bracket} uses only $v_{k}>\ell_{k}>1$, so~\eqref{eq:Ck-bracket} holds for $v$ and $C(v)$. Let
$v^{\pm}=(\Upsilon^{\pm}+k-1)\sqrt k$. If $v^{-}(v^{-}-\ell_{k})>1$, then
\[C_{k}\ \le\ C(v^{+})\ \le\ v^{+}+\frac{v^{+}}{(v^{+})^{2}-1}\,\log\frac{v^{+}}{v^{+}-\ell_{k}} .\]
The script \texttt{numerics/constant\_Ck\_certified.py}, in the folder \texttt{paper2} of the repository cited in the
main paper, does this computation with outward rounding. It checks that $v^{-}(v^{-}-\ell_{k})>1$, and it finds
$C_{k}/(2^{k/2}\sqrt k)\le1.72$ for every odd $k$ from $3$ to $123$, and at most $1.02$ for every odd $k$ from $21$
to $123$. The largest value is at $k=5$, in $[1.71865,1.71869]$.
\end{proof}

\section{Normalisations in the literature}\label{app:norm}

The noise tensor of \citet{perry2020statistical} is $\sqrt{2/d}\,W$ in law, so their signal is our
$\lambda/(\sigma\sqrt{d/2})$. At fixed $k$, their Theorem~3.2 makes weak detection impossible below the largest signal
$\lambda_{\rm PWB}$ with $\frac12\lambda_{\rm PWB}^{2}t^{k}/(1+t^{k})\le-\frac12\log(1-t^{2})$ for all $t\in[0,1)$,
which is~$1.5847$ for $k=3$, that is $\lambda/\sigma=1.1205\sqrt d$. \citet{jagannath2020statistical}
observe~$\lambda\sqrt d\,x^{\otimes k}+W'$, with $x$ uniform on the sphere and $W'$ a tensor of independent standard Gaussian
entries. The projection of $W'$ on the symmetric tensors has the law of $W$, and the rest is independent of it and
does not depend on the signal, so their $\lambda$ is our $\lambda/(\sigma\sqrt d)$, and their critical value, $1.2066$
for $k=3$, reads $\lambda/\sigma=1.2066\sqrt d$. \citet{lesieur2017statistical} work at fixed
signal-to-noise ratio as $d\to\infty$, with a spike whose entries are independent draws from a known prior, which
is $\lambda\asymp\sqrt d$ in our normalisation.

\section{Closed forms of $G_{r}$ and $H_{r}$}\label{app:closed}

The functions $G_{r}$ and $H_{r}$ of Section~\ref*{sub:closed} of the main paper, with the notation of that
paper, are finite sums of incomplete moments.

\begin{proposition}\label{prop:closed}
Write $\gamma=\sqrt{\kappa/(m-1)}$ and $\nu=m-1$, and let $e_{j}(r)$ be the elementary symmetric functions of the
eigenvalues of $r$. Then
\begin{subequations}
\begin{equation}
\begin{gathered}G_{r}(\ell)=\sum_{j=0}^{N}(-1)^{j}e_{j}(r)\,I_{N-j}(\ell),\\
I_{0}=1-\Phi,\quad I_{1}=\varphi,\quad I_{p}(\ell)=\ell^{p-1}\varphi(\ell)+(p-1)I_{p-2}(\ell),\end{gathered}
\end{equation}
and
\begin{equation}
\begin{gathered}H_{r}(\ell)=\gamma\sum_{j=0}^{N}(-1)^{j}e_{j}(r)\,\gamma^{\,N-j}\,K_{N-j}(\ell/\gamma),\\
K_{p}(a)=\frac{c_{\nu}\nu^{\frac{p+1}{2}}}{2}B\Big(\tfrac{\nu-p}{2},\tfrac{p+1}{2}\Big)I_{w_{a}}\Big(\tfrac{\nu-p}{2},\tfrac{p+1}{2}\Big)\end{gathered}
\end{equation}
\end{subequations}
for $a\ge0$, with $w_{a}=\nu/(\nu+a^{2})$, $c_{\nu}=\Gamma(\frac{\nu+1}{2})/(\sqrt{\nu\pi}\,\Gamma(\frac\nu2))$ and
$I_{w}$ the regularised incomplete Beta function, the integrals converging because $N<\nu$. Both $G_{r}$ and
$H_{r}$ are positive and strictly decreasing on $(\lmax(r),\infty)$, with limit $0$ at $+\infty$. The prefactor
$\gamma$ and the constant $c_{\nu}$ cancel in $\mathbf T$.
\end{proposition}

\begin{proof}
Expanding $\det(u\,\mathrm{Id}-r)=\sum_{j}(-1)^{j}e_{j}(r)u^{N-j}$ reduces both statements to incomplete moments.
For $I_{p}$, integrate by parts: $\int_{\ell}^{\infty}u^{p-1}\,(u\varphi(u))\,\mathrm{d}u=\ell^{p-1}\varphi(\ell)+(p-1)I_{p-2}(\ell)$.
For $K_{p}$, substitute $u=\gamma y$ and then $w=\nu/(\nu+y^{2})$ in
$\int_{a}^{\infty}y^{p}(1+y^{2}/\nu)^{-(\nu+1)/2}\,\mathrm{d}y$, which becomes
$\tfrac12\nu^{(p+1)/2}\int_{0}^{w_{a}}w^{\frac{\nu-p}{2}-1}(1-w)^{\frac{p+1}{2}-1}\,\mathrm{d}w$, finite for~$p<\nu$.
Monotonicity and positivity hold because the integrand $\det(u\,\mathrm{Id}-r)$ is positive for~$u>\lmax(r)$.
\end{proof}


\begin{thebibliography}{25}
\providecommand{\natexlab}[1]{#1}
\providecommand{\url}[1]{\texttt{#1}}
\expandafter\ifx\csname urlstyle\endcsname\relax
  \providecommand{\doi}[1]{doi: #1}\else
  \providecommand{\doi}{doi: \begingroup \urlstyle{rm}\Url}\fi

\bibitem[Anandkumar et~al.(2014)Anandkumar, Ge, Hsu, Kakade, and
  Telgarsky]{anandkumar2014tensor}
Animashree Anandkumar, Rong Ge, Daniel Hsu, Sham~M. Kakade, and Matus
  Telgarsky.
\newblock Tensor decompositions for learning latent variable models.
\newblock \emph{Journal of Machine Learning Research}, 15\penalty0
  (1):\penalty0 2773--2832, 2014.

\bibitem[Armentano et~al.(2025)Armentano, Aza{\"\i}s, and
  Le\'on]{armentano2023general}
Diego Armentano, Jean-Marc Aza{\"\i}s, and Jos\'e~R. Le\'on.
\newblock On a general {K}ac--{R}ice formula for the measure of a level set.
\newblock \emph{The Annals of Applied Probability}, 35\penalty0 (3):\penalty0
  1828--1851, 2025.
\newblock \doi{10.1214/25-AAP2158}.

\bibitem[Auffinger et~al.(2013)Auffinger, {Ben Arous}, and
  {\v{C}}ern{\'y}]{auffinger2013}
Antonio Auffinger, G{\'e}rard {Ben Arous}, and Ji{\v{r}}{\'\i} {\v{C}}ern{\'y}.
\newblock Random matrices and complexity of spin glasses.
\newblock \emph{Communications on Pure and Applied Mathematics}, 66\penalty0
  (2):\penalty0 165--201, 2013.
\newblock \doi{10.1002/cpa.21422}.
\newblock arXiv:1003.1129v2.

\bibitem[Aza{\"\i}s and Wschebor(2009)]{azais2009level}
Jean-Marc Aza{\"\i}s and Mario Wschebor.
\newblock \emph{Level sets and extrema of random processes and fields}.
\newblock John Wiley \& Sons, 2009.
\newblock \doi{10.1002/9780470434642}.

\bibitem[Aza{\"\i}s et~al.(2018)Aza{\"\i}s, De~Castro, and
  Mourareau]{azais2018power}
Jean-Marc Aza{\"\i}s, Yohann De~Castro, and St{\'e}phane Mourareau.
\newblock Power of the spacing test for least-angle regression.
\newblock \emph{Bernoulli}, 24\penalty0 (1):\penalty0 465--492, 2018.
\newblock \doi{10.3150/16-BEJ885}.

\bibitem[Aza{\"\i}s et~al.(2024)Aza{\"\i}s, Dalmao, and
  De~Castro]{azais2024second}
Jean-Marc Aza{\"\i}s, Federico Dalmao, and Yohann De~Castro.
\newblock Second maximum of a {G}aussian random field and exact ($t$-)spacing
  test.
\newblock \emph{arXiv preprint arXiv:2406.18397v3}, 2024.

\bibitem[Aza{\"\i}s et~al.(2026)Aza{\"\i}s, Dalmao, and De~Castro]{kss}
Jean-Marc Aza{\"\i}s, Federico Dalmao, and Yohann De~Castro.
\newblock Non-asymptotic first- and second-order tail brackets for the
  {K}ostlan--{S}hub--{S}male field and spherical $k$-spin complexity.
\newblock \emph{arXiv preprint arXiv:2611.09459}, 2026.

\bibitem[Baik et~al.(2005)Baik, {Ben Arous}, and P{\'e}ch{\'e}]{baik2005phase}
Jinho Baik, G{\'e}rard {Ben Arous}, and Sandrine P{\'e}ch{\'e}.
\newblock Phase transition of the largest eigenvalue for nonnull complex sample
  covariance matrices.
\newblock \emph{The Annals of Probability}, 33\penalty0 (5):\penalty0
  1643--1697, 2005.
\newblock \doi{10.1214/009117905000000233}.

\bibitem[Banach(1938)]{banach1938homogene}
S.~Banach.
\newblock {\"U}ber homogene {P}olynome in ({$L^{2}$}).
\newblock \emph{Studia Mathematica}, 7\penalty0 (1):\penalty0 36--44, 1938.
\newblock \doi{10.4064/sm-7-1-36-44}.

\bibitem[{Ben Arous} et~al.(2019){Ben Arous}, Mei, Montanari, and
  Nica]{arous2019landscape}
G{\'e}rard {Ben Arous}, Song Mei, Andrea Montanari, and Mihai Nica.
\newblock The landscape of the spiked tensor model.
\newblock \emph{Communications on Pure and Applied Mathematics}, 72\penalty0
  (11):\penalty0 2282--2330, 2019.
\newblock \doi{10.1002/cpa.21861}.
\newblock arXiv:1711.05424v2.

\bibitem[Ben~Arous et~al.(2024{\natexlab{a}})Ben~Arous, Gerbelot, and
  Piccolo]{benarous2024langevin}
G{\'e}rard Ben~Arous, C{\'e}dric Gerbelot, and Vanessa Piccolo.
\newblock Langevin dynamics for high-dimensional optimization: the case of
  multi-spiked tensor {PCA}.
\newblock \emph{arXiv preprint arXiv:2408.06401}, 2024{\natexlab{a}}.

\bibitem[Ben~Arous et~al.(2024{\natexlab{b}})Ben~Arous, Gerbelot, and
  Piccolo]{benarous2024sgd}
G{\'e}rard Ben~Arous, C{\'e}dric Gerbelot, and Vanessa Piccolo.
\newblock Stochastic gradient descent in high dimensions for multi-spiked
  tensor {PCA}.
\newblock \emph{arXiv preprint arXiv:2410.18162}, 2024{\natexlab{b}}.

\bibitem[Friedland and Lim(2018)]{friedland2018nuclear}
Shmuel Friedland and Lek-Heng Lim.
\newblock Nuclear norm of higher-order tensors.
\newblock \emph{Mathematics of Computation}, 87\penalty0 (311):\penalty0
  1255--1281, 2018.
\newblock \doi{10.1090/mcom/3239}.
\newblock arXiv:1410.6072v3.

\bibitem[Jagannath et~al.(2020)Jagannath, Lopatto, and
  Miolane]{jagannath2020statistical}
Aukosh Jagannath, Patrick Lopatto, and L{\'e}o Miolane.
\newblock Statistical thresholds for tensor {PCA}.
\newblock \emph{The Annals of Applied Probability}, 30\penalty0 (4):\penalty0
  1910--1933, 2020.
\newblock \doi{10.1214/19-AAP1547}.
\newblock arXiv:1812.03403v2.

\bibitem[Lehmann and Romano(2005)]{lehmann2005testing}
Erich~L. Lehmann and Joseph~P. Romano.
\newblock \emph{Testing Statistical Hypotheses}.
\newblock Springer, 3rd edition, 2005.

\bibitem[{Leiden Declaration}(2026)]{leiden2026declaration}
{Leiden Declaration}.
\newblock The {Leiden} {Declaration} on {Artificial} {Intelligence} and
  {Mathematics}.
\newblock \url{https://leidendeclaration.ai/}, 2026.
\newblock Endorsed by the International Mathematical Union.

\bibitem[Lesieur et~al.(2017)Lesieur, Miolane, Lelarge, Krzakala, and
  Zdeborov{\'a}]{lesieur2017statistical}
Thibault Lesieur, L{\'e}o Miolane, Marc Lelarge, Florent Krzakala, and Lenka
  Zdeborov{\'a}.
\newblock Statistical and computational phase transitions in spiked tensor
  estimation.
\newblock \emph{Proceedings of the 2017 IEEE International Symposium on
  Information Theory (ISIT)}, pages 511--515, 2017.
\newblock \doi{10.1109/ISIT.2017.8006580}.

\bibitem[Lifshits(1983)]{lifshits1983absolute}
M.~A. Lifshits.
\newblock {On the absolute continuity of distributions of functionals of random
  processes}.
\newblock \emph{Theory of Probability \& Its Applications}, 27\penalty0
  (3):\penalty0 600--607, 1983.
\newblock \doi{10.1137/1127066}.

\bibitem[Perry et~al.(2020)Perry, Wein, and Bandeira]{perry2020statistical}
Amelia Perry, Alexander~S. Wein, and Afonso~S. Bandeira.
\newblock Statistical limits of spiked tensor models.
\newblock \emph{Annales de l'Institut Henri Poincar{\'e}. Probabilit{\'e}s et
  Statistiques}, 56\penalty0 (1):\penalty0 230--264, 2020.
\newblock \doi{10.1214/19-AIHP960}.

\bibitem[Richard and Montanari(2014)]{montanari2014statistical}
Emile Richard and Andrea Montanari.
\newblock A statistical model for tensor {PCA}.
\newblock In \emph{Advances in Neural Information Processing Systems 27}, pages
  2897--2905, 2014.

\bibitem[Subag(2017)]{subag2017complexity}
Eliran Subag.
\newblock The complexity of spherical $p$-spin models---a second moment
  approach.
\newblock \emph{The Annals of Probability}, 45\penalty0 (5):\penalty0
  3385--3450, 2017.
\newblock \doi{10.1214/16-AOP1139}.

\bibitem[Taylor et~al.(2016)Taylor, Loftus, and
  Tibshirani]{taylor2016inference}
Jonathan~E. Taylor, Joshua~R. Loftus, and Ryan~J. Tibshirani.
\newblock Inference in adaptive regression via the {K}ac--{R}ice formula.
\newblock \emph{The Annals of Statistics}, 44\penalty0 (2):\penalty0 743--770,
  2016.
\newblock \doi{10.1214/15-AOS1386}.

\bibitem[Wijsman(1967)]{wijsman1967cross}
Robert~A. Wijsman.
\newblock Cross-sections of orbits and their application to densities of
  maximal invariants.
\newblock In \emph{Proceedings of the Fifth Berkeley Symposium on Mathematical
  Statistics and Probability}, volume~1, pages 389--400. University of
  California Press, 1967.

\bibitem[Xia et~al.(2022)Xia, Zhang, and Zhou]{xia2022inference}
Dong Xia, Anru~R. Zhang, and Yuchen Zhou.
\newblock Inference for low-rank tensors---no need to debias.
\newblock \emph{The Annals of Statistics}, 50\penalty0 (2):\penalty0
  1220--1245, 2022.
\newblock \doi{10.1214/21-AOS2146}.

\bibitem[Zhang and Xia(2018)]{zhang2018tensor}
Anru Zhang and Dong Xia.
\newblock Tensor {SVD}: statistical and computational limits.
\newblock \emph{IEEE Transactions on Information Theory}, 64\penalty0
  (11):\penalty0 7311--7338, 2018.
\newblock \doi{10.1109/TIT.2018.2841377}.

\end{thebibliography}
\end{document}